\documentclass[twocolumn]{autart}    % Enable this line and disable the 
\usepackage[utf8]{inputenc}   % This two packages are used for \v and \'
\usepackage[T1]{fontenc}      % same as previous

\usepackage{subfigure, amsmath, rotating}
\usepackage{array, multirow}
\usepackage{amsbsy, amsfonts, graphics}
\usepackage{algorithmic}
\usepackage{algorithm}
\usepackage{epstopdf}
\usepackage{mathrsfs}
\usepackage{graphicx}
\usepackage{caption}
\usepackage{amssymb}
\usepackage{amsmath,bm}
\usepackage{subfiles}
\usepackage{mathtools}

\usepackage{tikz}
\usepackage{float}
\usetikzlibrary{calc,patterns,decorations.pathmorphing,decorations.markings}
\usepackage{soul}
\usepackage{cancel}
\usepackage{todonotes}

\usepackage{subfig}
\usepackage[normalem]{ulem} % 给删除线准备的，haing [normalem] makes \emph{} to be italic not underline.
\usepackage{units}
\newtheorem{sassum}{Standing Assumption}

\usepackage{wrapfig}
\usepackage{graphicx}                              % document contains figures,
\makeatletter
\renewcommand{\maketag@@@}[1]{\hbox{\m@th\normalsize\normalfont#1}}%
\makeatother

\begin{document}

\begin{frontmatter}

\title{Stabilization of singularly perturbed networked control systems over a single channel\thanksref{footnoteinfo}} % Title, preferably not more 

\thanks[footnoteinfo]{This work was supported by the Australian Research Council under the Discovery Project DP200101303, the France Australia collaboration project IRP-ARS CNRS and the ANR COMMITS ANR-23-CE25-0005.}

\author[Melbourne]{Weixuan Wang}\ead{weixuanw@student.unimelb.edu.au},
\author[Chile]{Alejandro I. Maass}\ead{alejandro.maass@uc.cl},
\author[Melbourne]{Dragan Ne\v{s}i\'{c}}\ead{dnesic@unimelb.edu.au},
\author[Melbourne]{Ying Tan}\ead{yingt@unimelb.edu.au},
\author[France]{Romain Postoyan}\ead{romain.postoyan@univ-lorraine.fr},
\author[Netherlands]{W.P.M.H. Heemels}\ead{w.p.m.h.heemels@tue.nl}

\address[Melbourne]{School of Electrical, Mechanical and Infrastructure Engineering, The University of Melbourne, Parkville, 3010, Victoria, Australia}
\address[Chile]{Department of Electrical Engineering, Pontificia Universidad Cat\'olica de Chile, Santiago, 7820436, Chile}
\address[France]{Universit\'e de Lorraine, CNRS, CRAN, F-54000 Nancy, France}
\address[Netherlands]{Department of Mechanical Engineering, Eindhoven University of Technology, The Netherlands}
% \address[Paestum]{Buckingham Palace, Paestum}  % Please supply                                              
% \address[Rome]{Senate House, Rome}             % full addresses
% \address[Baiae]{The White House, Baiae}        % here.

                             % keyword list or with the 
                                          % help of the Automatica 
                                          % keyword wizard

\begin{abstract}                          % Abstract of not more than 200 words.
% This paper studies the emulation-based stabilization of (nonlinear) networked control systems with two time scales. 
% We consider the scenario where only a single communication channel is used to transmit both fast and slow variables between the plant and the controller.
% %
% The challenge is then to appropriately schedule transmissions to (approximately) preserve the stability properties of the closed-loop system in case of perfect, continuous communications. We present for this purpose a novel dual clock mechanism. The networked control system is modeled as a hybrid singularly perturbed dynamical system. 
% %
% Singular perturbation-based analysis is used to obtain individual maximum allowable transmission intervals for the transmission of the fast and slow variables, under which semi-global practical asymptotic stability properties hold. Stronger stability guarantees are also derived by strengthening the made assumptions. 
% %
% We illustrate the results via a numerical example.
%

This paper studies the emulation-based stabilization of nonlinear networked control systems with two time scales. We address the challenge of using a single communication channel for transmitting both fast and slow variables between the plant and the controller. A novel dual clock mechanism is proposed to schedule transmissions for this purpose. The system is modeled as a hybrid singularly perturbed dynamical system, and singular perturbation analysis is employed to determine individual maximum allowable transmission intervals for both fast and slow variables, ensuring semi-global practical asymptotic stability. Enhanced stability guarantees are also provided under stronger assumptions. The efficacy of the proposed method is illustrated through a numerical example.
\end{abstract}

\end{frontmatter}

% \section{Introduction}
% Video, patres conscripti, in me omnium vestrum ora atque oculos esse 
% conversos, video vos non solunn de vestro ac rei publicae, verum 
% etiam, si id depulsum sit, de meo periculo esse sollicitos. Est mihi 
% iucunda in malis et grata in dolore vestra erga me voluntas, sed eam, 
% per deos inmortales, deponite atque obliti salutis meae de vobis ac 
% de vestris liberis cogitate. Mihi si haec condicio consulatus data 
% est, ut omnis acerbitates, onunis dolores cruciatusque perferrem, 
% feram non solum fortiter, verum etiam lubenter, dum modo meis 
% laboribus vobis populoque Romano dignitas salusque pariatur.

% \begin{figure}
% \begin{center}
% \includegraphics[height=4cm]{jcaesar.eps}    % The printed column  
% \caption{Gaius Julius Caesar, 100--44 B.C.}  % width is 8.4 cm.
% \label{fig1}                                 % Size the figures 
% \end{center}                                 % accordingly.
% \end{figure}

% OR

%\begin{figure}
%\begin{center}
%\epsfig{file=jcaesar,width=7cm}
%\caption{Gaius Julius Caesar, 100--44 B.C.}
%\label{fig1}
%\end{center}
%\end{figure}

\section{Introduction}
Networked control systems (NCSs) integrate feedback control loops with real-time communication networks \cite{definition}. The rapid evolution of network technologies has expanded the applications of NCSs across various sectors, such as industrial automation, smart transportation, telemedicine, and both space and terrestrial exploration \cite{xu2018survey}. These applications often involve dynamic variables evolving in multiple time scales \cite{kokotovic_singular_book}. 
Most state-of-the-art research in NCSs, e.g., \cite{dragan_stability,carnevale_stability,heijmans2017computing}, overlook this multi-scale structure, leading to designs demanding high data transmission rates to maintain system stability or robustness. This is especially problematic for Internet of Things devices, which are generally wireless, battery-powered, and have limited bandwidth, processing power, and data storage capacity.
A relevant approach in this context is to extend the singularly perturbed method \cite{nonlinear_systems_Khalil} to NCS, thereby introducing the concept of singularly perturbed NCS (SPNCS).

SPNCSs have gained substantial attention for their applicability in various engineering disciplines. 
For instance, \cite{wang2021observer,song2019dynamic,lei2022event} have introduced a range of control and analysis techniques for linear SPNCSs. 
For nonlinear SPNCSs, \cite{Romain_ETC} formulated stabilizing event-triggered feedback laws based on the reduced model (also known as the quasi-steady-state model), assuming stable fast dynamics. 
Meanwhile, \cite{SPNCS} established sufficient conditions for stability in time-triggered, two-time-scale nonlinear SPNCSs under the scenario that the controller is co-located with the actuators, as well as slow and fast plant states are transmitted over two separate channels, which is not always possible or desirable in practice.
More precisely, \cite{lei2022event}, \cite{Romain_ETC} and \cite{SPNCS} adopted an emulation-based approach for NCS design \cite{dragan_stability}, i.e., the controller is designed to stabilize the plant in absence of communication constraints, and subsequently the event- or time-triggered conditions are determined to maintain stability over networks.
Furthermore, SPNCSs have been framed within the framework of a hybrid singularly perturbed dynamical system using the formalisms outlined in \cite{gosate12}, with stability tools available in \cite{sanfelice2011singular} and \cite{wang2012analysis}.
%framed within the context of a hybrid singularly perturbed dynamical system with the formalisms of \cite{gosate12}, for which stability tools are available in \cite{sanfelice2011singular} and \cite{wang2012analysis}. % and \cite{6161309}.
%
We highlight that, the existing literature on SPNCSs often assumes either stable slow or fast subsystems, perfect transmission of some signals, and dedicated channels for slow and fast signals.

In this paper, we consider a two-time scale nonlinear system stabilized by a dynamical output feedback controller, inspired by both linear \cite{linear1980}, \cite{linear2018}, \cite{linear2010} and nonlinear \cite{output_feedback} research on dynamic controllers for singularly perturbed dynamical systems (SPSs) \cite{nonlinear_systems_Khalil} without network constraints. Unlike previous work \cite{SPNCS} that assumes only the transmission of plant states, we address scenarios where both plant output and control input are transmitted via the network. 
%Additionally, instead of requiring dedicated channels for slow and fast variables as in \cite{SPNCS}, which is not always practical, we consider the scenario where only one channel is involved, which presents challenges in allocating access for slow and fast signals. 
%
Additionally, instead of requiring dedicated channels for slow and fast variables as in \cite{SPNCS}, which is not always practical, we consider a single-channel scenario, which presents significant challenges in allocating access for slow and fast signals.
Therefore, our objective is to provide a general methodology for the design of the controller and transmission mechanism to ensure stability properties for the SPNCS.
%to design a resource-efficient mechanism to manage the transmission of slow and fast signals over a single communication channel, and \red{identify conditions under which stability properties are preserved.}  

We first design a dual clock mechanism to govern the data transmissions.
Then we represent the SPNCS as a hybrid SPS, incorporating jump sets specifically designed to comply with the previously mentioned clock mechanism.
%Then we model the SPNCS as a hybrid SPS\cyan{;} with jump sets designed to satisfy the aforementioned clock mechanism. 
The obtained SPS is more general compared to those in the literature \cite{sanfelice2011singular,wang2012analysis}%,6161309}
, as its flow and jump sets depend on the time scale separation parameter, which is commonly denoted by $\epsilon$.
Following an emulation-based approach, we demonstrate that the if the \emph{reduced system} and \emph{boundary-layer system} are \emph{uniformly globally asymptotically stable} (UGAS) or \emph{uniformly globally exponentially stable} (UGES), and an interconnection condition and some mild conditions are met, the stability properties can be approximately preserved by transmitting both slow and fast variables sufficiently fast.
Specifically, we employ a Lyapunov-based analysis to determine individual \emph{maximum allowable transmission interval} (MATI) \cite{dragan_stability,carnevale_stability} of the slow and fast dynamics.
Finally, we illustrate the benefits of our approach through a numerical case study.

Compared to our preliminary work \cite{Single_channel_NCS_CDC}, we relax a restrictive condition on the \emph{minimum allowable transmission interval} (MIATI). Additionally, we present conditions that guarantee stronger stability properties: UGAS and UGES. Furthermore, while many works, such as \cite{Romain_ETC}, assume that either the slow or fast subsystem is stable without the need of control, 
%and then achieve the stabilisation of the overall system by implementing a stabilizing controller over the network, 
our methodology does not rely on this assumption. 
Nevertheless, our framework accommodates scenarios involving stable slow or fast subsystems as special cases.

The results are novel even for linear time-invariant (LTI) systems. While \cite{wang2021observer,song2019dynamic,lei2022event} assumed periodic transmissions, we allow transmissions to be aperiodic. Compared to \cite{song2019dynamic,lei2022event} that assumed the sampled-data structure, we consider SPNCSs with scheduling protocols. Moreover, we take into account the inter-event continuous behavior, which was ignored by \cite{wang2021observer,song2019dynamic}. 
%
%Finally, we illustrate the benefits of our approach through a numerical case study.

\textbf{Notation:} 
Let $\mathbb{R}\coloneqq (-\infty, \infty)$, $\mathbb{R}_{\geq 0} \coloneqq [0,\infty)$, $\mathbb{Z}_{\geq 0} \coloneqq \{0, 1, 2, \cdots \}$ and $\mathbb{Z}_{\geq 1} \coloneqq \{1, 2, \cdots \}$.
%The sets of real numbers and integers larger than or equal to an integer $n$ are denoted by $\mathbb{R}_{\geq n}$ and $\mathbb{Z}_{\geq n}$, respectively. 
For vectors $v_i\in \mathbb{R}^n$, $i\in \{1,2,\cdots, N\}$, we denote the vector $[v_1^\top \; v_2^\top \; \cdots \; v_N^\top]^\top$ by $(v_1, v_2, \cdots, v_N)$, and inner product by $\left< \cdot , \cdot \right>$. 
Given a vector $x\in \mathbb{R}^{n_x}$ and a non-empty closed set $\mathcal{A} \subseteq \mathbb{R}^{n_x}$, the distance from $x$ to $\mathcal{A}$ is denoted by $|x|_\mathcal{A} \coloneqq \min_{y\in \mathcal{A}}|x-y|$.
We use $U^\circ(x;v)$ to denote the Clarke generalized derivative \cite[Eqn. (20)]{teel2000assigning} of a locally Lipschitz function $U$ at $x$ in the direction of $v$. 
For a real symmetric matrix $P$, we denote its maximum and minimum eigenvalues by $\lambda_{\text{max}}(P)$ and $\lambda_{\text{min}}(P)$ respectively. The logic AND operator is denoted by $\wedge$.

\section{Problem setting} \label{Chapter Problem setting}
\begin{figure}[H]
    \centering
    \includegraphics[width = 0.6\linewidth]{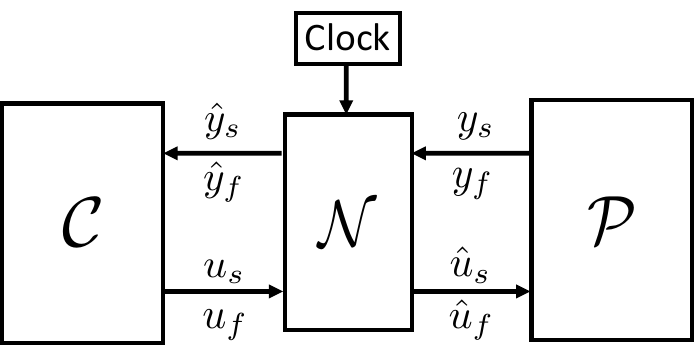}
    \caption{NCS Block Diagram}
    \label{fig: Block Diagram}
\end{figure}

We consider a two-time-scale nonlinear NCS as depicted in Figure \ref{fig: Block Diagram}, designed using emulation techniques \cite{dragan_stability}. Specifically, a dynamic continuous output-feedback controller is developed to ensure robustness for both the \emph{reduced} (slow) system and the \emph{boundary-layer} (fast) system, initially without considering the network. Subsequently, the network is designed by establishing bounds on transmission intervals and selecting an appropriate scheduling protocol \cite{dragan_stability}. The resulting continuous-time controller is then deployed over the network, with the objective of providing conditions under which the stability of the SPNCS is guaranteed. Details on the emulation design framework are provided in Section 5.
Next, we introduce the model of Figure \ref{fig: Block Diagram}.
%
% We consider a two-time-scale nonlinear NCS, shown in Figure \ref{fig: Block Diagram}, which is designed by emulation \cite{dragan_stability}. 
% \red{In particular, a dynamic output-feedback controller is assumed to be designed to ensure robustness properties for both the \emph{reduced system} (slow) and the \emph{boundary-layer system} (fast) without network constraint. Then we design the network by defining the bounds on transmission intervals and selecting the scheduling protocol \cite{dragan_stability}.
% %
% The controller is then deployed over the network, and 
% our aim is to provide condition on the original closed-loop system and the network under which stability of the SPNCS follows. See Section \ref{Section Emulation design framework} for details one the emulation design framework.}
% %appropriate MATIs are selected to ensure the stabilization of the NCS. 
% Next, we introduce the model of Figure \ref{fig: Block Diagram}.

\subsection{Plant ($\mathcal{P}$) and Controller ($\mathcal{C}$)}
%Let $n_{x_p}, n_{z_p}, n_{y_s}, n_{y_f} \in \mathbb{Z}_{\geq 0}$.
We model the plant as the following SPS,
\begin{equation}
    \mathcal{P}:
    \begin{cases}
    \begin{aligned}
    \dot x_p &= f_p(x_p, z_p,\hat u)\\
    \epsilon \dot z_p &= g_p(x_p, z_p, \hat u) \\
    y_p &= \left(y_s, y_f  \right) = \left(k_{p_s}(x_p) , k_{p_f}(x_p, z_p) \right) ,
    \end{aligned}
    \end{cases} 
    \label{eqn:plant}
\end{equation}
where $0 < \epsilon \ll 1$, $x_p \in \mathbb{R}^{n_{x_p}}$, $z_p\in \mathbb{R}^{n_{z_p}}$, $y_s\in \mathbb{R}^{n_{y_s}}$, $y_f \in \mathbb{R}^{n_{y_f}}$ and $n_{x_p}, n_{z_p}, n_{y_s}, n_{y_f} \in \mathbb{Z}_{\geq 0}$.
Here, $x_p$ and $z_p$ denote the slow and fast plant states, respectively, while $y_s$ and $y_f$ represent the slow and fast output, respectively. Additionally, $\hat u = (\hat u_s, \hat u_f)$ refers to the latest received control input $u$ in \eqref{eqn:controller} from the network. It is assumed that $k_{p_s}$ and $k_{p_f}$ are continuously differentiable, and $f_p$ and $g_p$ are locally Lipschitz.
%$f_p(0,0,0) = 0$, $g_p(0,0,0) = 0$, $k_{p_s}(0) = 0$, $k_{p_f}(0,0) = 0$ and $k_{p_f}$ is continuously differentiable.
%
%
%
% In our emulation-based approach, we assume that a dynamic controller has been designed to stabilise plant \eqref{eqn:plant} in the absence of network, \cyan{i.e., $\hat{y}_p \equiv y_p$ and $\hat{u} \equiv u$}, \todo{maybe remove this paragraph}
Similarly, the dynamic controller has the following form,
\begin{equation}
    \mathcal{C}:
    \begin{cases}
    \begin{aligned}
    \dot x_c &= f_c(x_c, z_c, \hat{y}_p)\\
    \epsilon \dot z_c &= g_c(x_c, z_c, \hat y_p) \\
    u &= (u_s, u_f) = \left(k_{c_s}(x_c), k_{c_f}(x_c,z_c) \right) ,
    \end{aligned}
    \end{cases}
    \label{eqn:controller}
\end{equation}
where $\epsilon$ comes from (\ref{eqn:plant}), $x_c \in \mathbb{R}^{n_{x_c}}$, $z_c \in \mathbb{R}^{n_{z_c}}$, $u_s \in \mathbb{R}^{n_{u_s}}$, $u_f \in \mathbb{R}^{n_{u_f}}$ and $n_{x_c}, n_{z_c}, n_{u_s}, n_{u_f} \in \mathbb{Z}_{\geq 0}$. Moreover, $\hat y_p = (\hat y_s, \hat y_f)$ refers to the most recently received output of the plant transmitted via the network. It is assumed that $k_{c_s}$ and $k_{c_f}$ are continuously differentiable, $f_c$ and $g_c$ are locally Lipschitz, and $u_s$, $u_f$, $y_s$, $y_f$ have the dimension as $\hat u_s$, $\hat u_f$, $\hat y_s$, $\hat y_f$, respectively.
%
%\sout{We have $n_{x_p} + n_{x_c} \geq 1$ and $n_{z_p} + n_{z_c} \geq 1$ to guarantee the existence of both the slow and fast states. We also have that $n_{y_s} + n_{u_s} \geq 1$ and $n_{y_f} + n_{u_f} \geq 1$, which ensures that both fast and slow signals are present in the system.}
%

%\sout{Our results generalize those from \cite{SPNCS} by considering the transmission of control inputs and plant outputs via the network, as opposed to transmitting the plant state and assuming the control input is transmitted through a perfect network.}
%
%
% \begin{figure}[H]
%     \centering
%     \includegraphics[width = 0.6\linewidth]{Figures/Block diagram small font.pdf}
%     \caption{NCS Block Diagram}
%     \label{fig: Block Diagram}
% \end{figure}

\subsection{Network ($\mathcal{N}$)}
A channel may consist of multiple \emph{network nodes}, each representing a group of sensors and/or actuators, see \cite{wang2017observer} for more information. In this paper, we consider that each node can only contain either slow (i.e., $y_s$, $u_s$) or fast (i.e., $y_f$, $u_f$) signals, but not both. Only one node can transmit data at any given transmission time, regulated by the channel scheduling protocol. This implies that slow signals are never transmitted simultaneously with fast signals. In particular, at each transmission time allocated to a slow (resp. fast) node, a group of elements in $y_s$ (resp. $y_f$) and $u_s$ (resp. $u_f$) accessible to that node is sampled and transmitted.

%In this context, we define $\mathcal{T} \coloneqq \{t_0, t_1, t_2, \cdots \}$ as a set of all transmission instants. Let $\mathcal{T}^s \coloneqq \{t_0^s, t_1^s, t_2^s, \cdots \} $ be the subsequence of $\mathcal{T}$ such that all the elements of $\mathcal{T}^s$ are the instances that a slow node gets access to the network. Then we define the set of instances that a fast note gets access to the network to be $\mathcal{T}^f \coloneqq \mathcal{T}-\mathcal{T}^s = \{t_0^f, t_1^f, t_2^f, \cdots \}$

In this context, we define $\mathcal{T} \coloneqq \{t_1, t_2, t_3, \cdots \}$ as the set of all transmission instants. Let $\mathcal{T}^s \coloneqq \{t_1^s, t_2^s, t_3^s, \cdots \}$ be the subsequence of $\mathcal{T}$ consisting of the instances that a slow node gains access to the network. We then define the set of instances that a fast node gets access to the network as $\mathcal{T}^f \coloneqq \mathcal{T} \setminus \mathcal{T}^s = \{t_1^f, t_2^f, t_3^f, \cdots \}$.
%
%
%
% Define $\mathcal{T}^s \coloneqq \{t_0^s, t_1^s, t_2^s, \cdots \} $ as the unbounded set of transmission times at which a slow node is transmitted,
% and $\mathcal{T}^f \coloneqq \{t_0^f, t_1^f, t_2^f, \cdots \}$ as the unbounded set of transmission times at which a fast node is transmitted,
% such that $\mathcal{T}^s \cap \mathcal{T}^f = \emptyset$. Then, let $\mathcal{T} \coloneqq \mathcal{T}^s \cup \mathcal{T}^f =  \{t_0, t_1, t_2, \cdots \} $ denote the set of all transmission instances, with its elements arranged in ascending time order.
%
We impose that for any $k \in \mathbb{Z}_{\geq 1}$, the transmission times satisfy
\begin{subequations}
    \begin{align}
    &\tau_{\text{miati}}^s \leq t_{k+1}^s - t_k^s \leq \tau_{\text{mati}}^s, \; \forall t_k^s,t_{k+1}^s\in \mathcal{T}^s,  \label{eqn: timer eqn1}
    \\
    &\tau_{\text{miati}}^f \leq t_{k+1}^f - t_k^f \leq \tau_{\text{mati}}^f ,  \; \forall t_k^f, t_{k+1}^f  \in \mathcal{T}^f,  
    \label{eqn: timer eqn2}
    \\
    &\tau_{\text{miati}}^f \leq t_{k+1} - t_k, \quad \qquad \; \; \ \forall t_k, t_{k + 1} \in \mathcal{T}, \label{eqn: timer eqn3}
    \end{align}
    \label{eqn: Stefan timer}%
\end{subequations}
\noindent where $0<\tau_{\text{miati}}^f\leq \tau_{\text{mati}}^f$ denote, respectively, the MIATI and MATI between any two consecutive fast transmissions. Similarly, $\tau_{\text{miati}}^s$ and $\tau_{\text{mati}}^s$ are the MIATI and MATI between two consecutive slow updates.
We note that since there might be a slow transmission between two consecutive fast transmissions,
\begin{equation}
    \tau_{\text{miati}}^f \leq  \tfrac{1}{2}\tau_{\text{mati}}^f
    \label{eqn: condition on miati^f}
\end{equation}
must hold to satisfy \eqref{eqn: timer eqn2} and \eqref{eqn: timer eqn3}, as in \cite{Stefan_thesis}.

Let the \emph{network-induced errors} be $e_{y_s} \coloneqq \hat{y}_s - y_s$, $e_{y_f} \coloneqq \hat{y}_f - y_f$, $  e_{u_s} \coloneqq \hat{u}_s - u_s$ and $  e_{u_f} \coloneqq \hat{u}_f - u_f $.
For simplicity, $(\hat{y}_s,\hat{y}_f,\hat{u}_s,\hat{u}_f)$ are assumed to be constant between any two successive transmission times, i.e., zero-order hold devices are used.
%Other type of network-processing may be implemented if desired, see, e.g., \cite{dragan_stability}.
Before we present the behaviour of the system at transmission times, we introduce some useful notation regarding the variables: $x\coloneqq (x_p,x_c)\in\mathbb{R}^{n_x}$, $z \coloneqq ( z_p, z_c) \in \mathbb{R}^{n_z}$, $e_s \coloneqq ( e_{y_s} , e_{u_s})\in \mathbb{R}^{n_{e_s}}$ and $e_f \coloneqq (e_{y_f} , e_{u_f}) \in \mathbb{R}^{n_{e_f}}$, with $n_x\coloneqq n_{x_p}+n_{x_c}$,  $n_z\coloneqq n_{z_p}+n_{z_c}$, $n_{e_s}\coloneqq n_{y_s}+n_{u_s}$ and  $n_{e_f}\coloneqq n_{y_f}+n_{u_f}$. 
 
At each transmission time $t_k^s \in \mathcal{T}^s$ for slow updates, the values $(\hat{y}_s,\hat{y}_f,\hat{u}_s,\hat{u}_f) $ are updated according to
$
\big(\hat{y}_s ( {t_k^s}^{+}),\hat{u}_s ( {t_k^s}^{+} )\big)
    =
    \big(
    y_s(t_k^s), u_s(t_k^s)
    \big)+ h_s(k, e_{s}(t_k^s) )
$
and
$
\big(\hat{y}_f ( {t_k^s}^{+} ), 
    \hat{u}_f ( {t_k^s}^{+})
    \big)
    =
    \left(
    \hat{y}_f(t_k^s), \hat{u}_f(t_k^s)
    \right)
$,
%
% \begin{equation*}
%     \begin{aligned}
%     \big(
%     \hat{y}_s ( {t_k^s}^{+}),
%     \hat{u}_s ( {t_k^s}^{+} )
%     \big)
%     =&
%     \big(
%     y_s(t_k^s), u_s(t_k^s)
%     \big)+ h_s(k, e_{s}(t_k^s) ),
%     \\
%     \big(
%     \hat{y}_f ( {t_k^s}^{+} ), 
%     \hat{u}_f ( {t_k^s}^{+})
%     \big)
%     =&
%     \left(
%     \hat{y}_f(t_k^s), \hat{u}_f(t_k^s)
%     \right) ,
%     \end{aligned}
% \end{equation*}
where the function $h_s: \mathbb{Z}_{\geq 0}\times \mathbb{R}^{n_{e_s}}  \rightarrow \mathbb{R}^{n_{e_s}}$ models the scheduling protocol \cite{dragan_stability} for the slow updates.
Similarly, for each $t_k^f \in \mathcal{T}^f$, we have 
$
\big(
    \hat{y}_s ( {t_k^f}^{+} ), 
    \hat{u}_s ( {t_k^f}^{+})
    \big)
    =
    \big(
    \hat{y}_s(t_k^f), \hat{u}_s(t_k^f)
    \big)
$
and
$\big(
    \hat{y}_f ( {t_k^f}^{+} ), 
    \hat{u}_f ( {t_k^f}^{+} )
    \big)
    =
    \big(
    y_f(t_k^f), u_f(t_k^f)
    \big) 
    + h_f\big(k, e_{f}(t_k^f) \big)$,
%
% \begin{align*}
%     \begin{aligned}
%     \big(
%     \hat{y}_s ( {t_k^f}^{+} ), 
%     \hat{u}_s ( {t_k^f}^{+})
%     \big)
%     =&
%     \big(
%     \hat{y}_s(t_k^f), \hat{u}_s(t_k^f)
%     \big),
%     \\
%     \big(
%     \hat{y}_f ( {t_k^f}^{+} ), 
%     \hat{u}_f ( {t_k^f}^{+} )
%     \big)
%     =&
%     \big(
%     y_f(t_k^f), u_f(t_k^f)
%     \big) 
%     + h_f\big(k, e_{f}(t_k^f) \big) ,
%     \end{aligned}
%    % \label{eqn: fast update}
% \end{align*}
where the function $h_f: \mathbb{Z}_{\geq 0}\times \mathbb{R}^{n_{e_f}} \rightarrow \mathbb{R}^{n_{e_f}} $ is the scheduling protocol for the update of fast components. 
If a SPNCS has $\ell$ slow nodes, then $e_s$ can be partitioned as $e_s = [e_{s,1}^\top \; e_{s,2}^\top \; \cdots \; e_{s,\ell}^\top]$. If the slow scheduling protocol $h_s$ grants the $i$th slow node access to the network at a transmission instance $t_k^s \in \mathcal{T}^s$, then $e_{s,i}$ experiences a jump. For protocols such as round robin (RR) and try-one-discard (TOD) \cite{dragan_stability}, $e_{s,i}({t_k^s}^+) = 0$ and $e_{s,j}({t_k^s}^+) = e_{s,j}({t_k^s})$ for all $j \neq i$, although this assumption is not generally necessary. The same rule applies to the fast nodes.

% A variable useful for analysis is the so-called \emph{network-induced error}, which we define as $e_{y_s} \coloneqq \hat{y}_s - y_s$, $e_{y_f} \coloneqq \hat{y}_f - y_f$, $  e_{u_s} \coloneqq \hat{u}_s - u_s$ and $  e_{u_f} \coloneqq \hat{u}_f - u_f $.
% For simplicity, $(\hat{y}_s,\hat{y}_f,\hat{u}_s,\hat{u}_f)$ are assumed to be constant between any two successive transmission times (i.e. zero-order hold behaviour). Other type of network-processing may be implemented if desired, see, e.g., \cite{dragan_stability}.
% Define $x\coloneqq (x_p,x_c)\in\mathbb{R}^{n_x}$, $z \coloneqq ( z_p, z_c) \in \mathbb{R}^{n_z}$, $e_s \coloneqq ( e_{y_s} , e_{u_s})\in \mathbb{R}^{n_{e_s}}$ and $e_f \coloneqq (e_{y_f} , e_{u_f}) \in \mathbb{R}^{n_{e_f}}$, with $n_x\coloneqq n_{x_p}+n_{x_c}$,  $n_z\coloneqq n_{z_p}+n_{z_c}$,  $n_{e_s}\coloneqq n_{y_s}+n_{u_s}$ and  $n_{e_f}\coloneqq n_{y_f}+n_{u_f}$. 

\section{A hybrid model for the SPNCS}
In this section, we present a hybrid system model for the SPNCS described in Section \ref{Chapter Problem setting} in the formalism of \cite{gosate12}, and it is more general than the hybrid SPSs in the literature such as \cite{sanfelice2011singular} and \cite{wang2012analysis}, as its flow and jump sets depend on $\epsilon$.
Firstly, we design a clock mechanism to satisfy \eqref{eqn: timer eqn1}-\eqref{eqn: timer eqn3}, and then we present the model of the overall SPNCS.

\subsection{Clock Mechanism}
We introduce two clocks and two counters, namely $\tau_s, \tau_f \in \mathbb{R}_{\geq 0}$ and $\kappa_s, \kappa_f \in \mathbb{Z}_{\geq 0}$. In particular, $\tau_s$ and $\epsilon \tau_f$ record the time elapsed since the last slow and fast transmission, respectively.
%, and we have $\dot{\tau}_s = 1$, $\epsilon \dot{\tau}_f = 1$ during flow.
Meanwhile, $\kappa_s$ and $\kappa_f$ count the number of slow and fast transmissions, respectively, and are useful for implementing some commonly used protocols, such as RR. 

\begin{figure}[H]
    \centering
    \includegraphics[width = \linewidth]{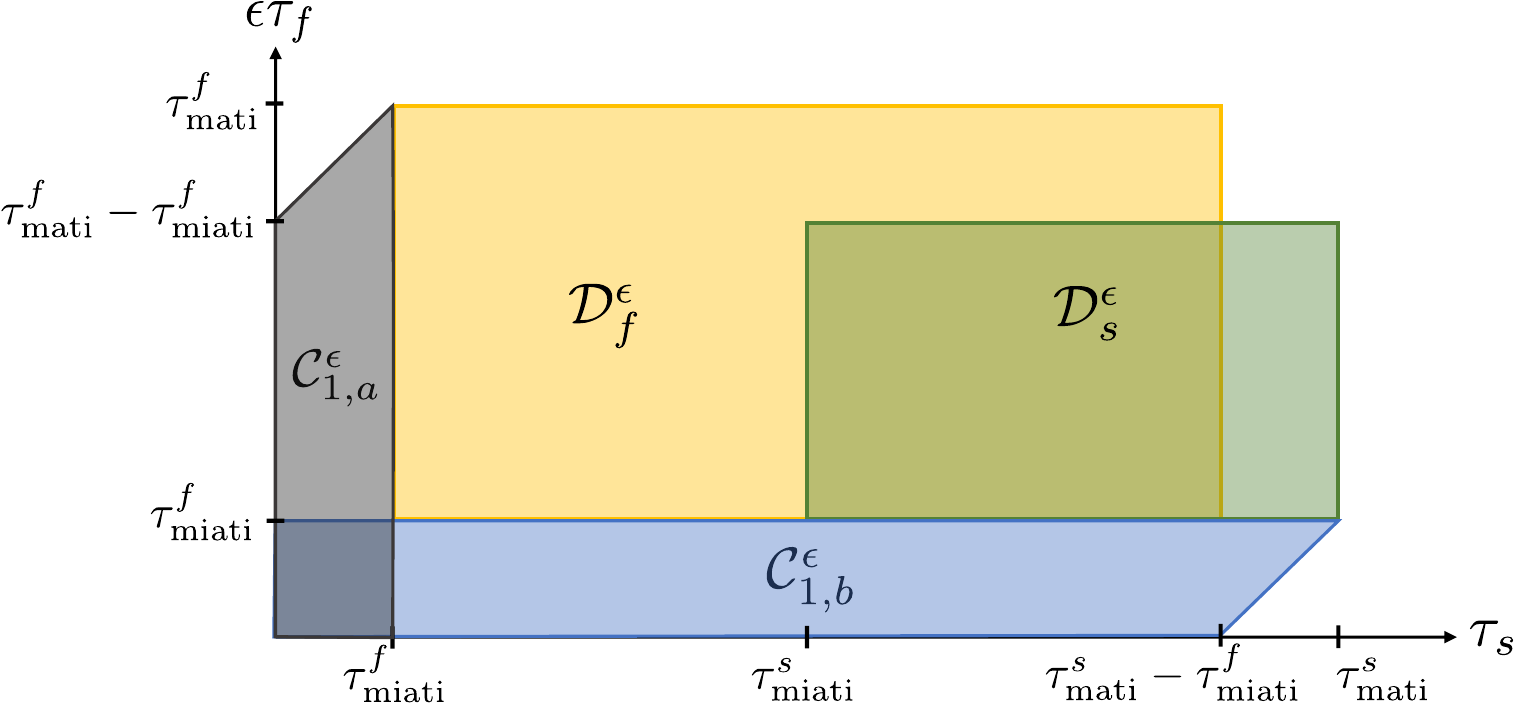}
    \caption{Flow set and jump set}
    \label{fig: Stefan timer}
\end{figure}

Let $\xi \coloneqq (x,e_s, \tau_s, \kappa_s, z,e_f, \tau_f,  \kappa_f)\in \mathbb{X}$,
with $\mathbb{X}\coloneqq \mathbb{R}^{n_x}\times \mathbb{R}^{n_{e_s}}\times  \mathbb{R}_{\geq 0} \times \mathbb{Z}_{\geq 0} \times \mathbb{R}^{n_z}\times \mathbb{R}^{n_{e_f}}\times \mathbb{R}_{\geq 0} \times \mathbb{Z}_{\geq 0}$, 
denote the full state of the hybrid system. We define the jump sets $\mathcal{D}_s^\epsilon$, $\mathcal{D}_f^\epsilon$ and the flow set $\mathcal{C}_1^\epsilon$ as
$\mathcal{D}_s^\epsilon \coloneqq  \{\xi \in \mathbb{X} \; | \; \tau_s \in [\tau_{\text{miati}}^s, \tau_{\text{mati}}^s] \wedge \epsilon \tau_f \in  [\tau_{\text{miati}}^f,  \tau_{\text{mati}}^f - \tau_{\text{miati}}^f]  \}$,
$\mathcal{D}_f^\epsilon \coloneqq \{\xi \in \mathbb{X} \; | \; \tau_s \in [\tau_{\text{miati}}^f, \tau_{\text{mati}}^s-\tau_{\text{miati}}^f]    \wedge\; \epsilon \tau_f \in  [\tau_{\text{miati}}^f, \tau_{\text{mati}}^f]   \}$,
and
$\mathcal{C}_1^\epsilon \coloneqq 
        \mathcal{D}_s^\epsilon \cup \mathcal{D}_f^\epsilon \cup \mathcal{C}_{1,a}^\epsilon \cup \mathcal{C}_{1,b}^\epsilon$,
%
% \begin{align*}
%     \mathcal{D}_s^\epsilon \coloneqq & \Big\{\xi \in \mathbb{X} \; | \; \tau_s \in [\tau_{\text{miati}}^s, \tau_{\text{mati}}^s]  \\
%         & \qquad \qquad \qquad \qquad \wedge \epsilon \tau_f \in  [\tau_{\text{miati}}^f,  \tau_{\text{mati}}^f - \tau_{\text{miati}}^f]  \Big\},
%     \\
%     \mathcal{D}_f^\epsilon \coloneqq &\Big\{\xi \in \mathbb{X} \; | \; \tau_s \in [\tau_{\text{miati}}^f, \tau_{\text{mati}}^s-\tau_{\text{miati}}^f]  \\
%     & \qquad \qquad \qquad \qquad \qquad \quad   \wedge\; \epsilon \tau_f \in  [\tau_{\text{miati}}^f, \tau_{\text{mati}}^f]  \Big \},
%     \\
%     \mathcal{C}_1^\epsilon \coloneqq & 
%         \mathcal{D}_s^\epsilon \cup \mathcal{D}_f^\epsilon \cup \mathcal{C}_{1,a}^\epsilon \cup \mathcal{C}_{1,b}^\epsilon  
% \end{align*}
%
with $\mathcal{C}_{1,a}^\epsilon \coloneqq \{ \xi \in \mathbb{X} \ | \ \tau_s \in [0, \tau_{\text{miati}}^f]  \wedge \epsilon \tau_f \in [0,\tau_s + \tau_{\text{mati}}^f - \tau_{\text{miati}}^f] \} $ and 
$\mathcal{C}_{1,b}^\epsilon \coloneqq \{ \xi \in \mathbb{X} \;|\; \tau_s \in [\tau_{\text{miati}}^f, \epsilon \tau_f + \tau_{\text{mati}}^s - \tau_{\text{miati}}^f]  \wedge \epsilon \tau_f \in  [0, \tau_{\text{miati}}^f]  \}$. 
A transmission of slow (resp. fast) signals is allowed in the set $\mathcal{D}_s^\epsilon$ (resp. $\mathcal{D}_f^\epsilon$), and at the transmission instance, $\tau_s$ (resp. $\tau_f$) is reset to zero.
The sets $\mathcal{C}_1^\epsilon$, $\mathcal{D}_s^\epsilon$ and $\mathcal{D}_f^\epsilon$ are defined to ensure the
satisfaction of \eqref{eqn: Stefan timer}, which can be deduced by visual inspection from Fig. \ref{fig: Stefan timer}. The jump sets $\mathcal{D}_s^\epsilon$ and $\mathcal{D}_f^\epsilon$ are indicated by the orange and green regions, respectively. Additionally, $\mathcal{C}_{1,a}^\epsilon$ and $\mathcal{C}_{1,b}^\epsilon$ are the regions where a jump is not allowed due to a recent transmission of slow and fast signals, respectively.

\subsection{Hybrid model}
Let $f_x,g_z,f_{e_s}$ and $g_{e_f}$ be defined in \eqref{eq:functions} in the next page, where we use $f_{x,\iota}$ and $g_{z,\iota}$, $\iota\in\{1,2\}$, to denote the $\iota$--th component of $f_x  $ and $g_z$, respectively.
%
% Let  
% $\xi \coloneqq (x,e_s, \tau_s, \kappa_s, z,e_f, \tau_f,  \kappa_f)\in \mathbb{X}$, 
% with $\mathbb{X}\coloneqq \mathbb{R}^{n_x}\times \mathbb{R}^{n_{e_s}}\times  \mathbb{R}_{\geq 0} \times \mathbb{Z}_{\geq 0} \times \mathbb{R}^{n_z}\times \mathbb{R}^{n_{e_f}}\times \mathbb{R}_{\geq 0} \times \mathbb{Z}_{\geq 0}$, 
% denote the full state of the hybrid system.
%
Then the SPNCS can now be expressed as the following hybrid model
\begin{equation}
    \mathcal{H}_1:\left\{
\begin{aligned}
    \dot{\xi} &= F(\xi, \epsilon), &&\xi \in \mathcal{C}_1^\epsilon, \\
    \xi^+ &\in G(\xi),  &&\xi\in \mathcal{D}_s^\epsilon \cup \mathcal{D}_f^\epsilon,
\end{aligned}
    \right.
    \label{eqn:full system}
\end{equation}
where 
%$F(\xi) \coloneqq  \big(f_x(x,z,e_s,e_f), \tfrac{1}{\epsilon}g_z(x,z,e_s,e_f),$ $f_{e_s}(x,z,e_s,e_f),\tfrac{1}{\epsilon} g_{e_f}(x,z,e_s,e_f, \epsilon), 1, \frac{1}{\epsilon},0,0\big)$
$F(\xi, \epsilon) \coloneqq  \big(f_x(x,z,e_s,e_f),f_{e_s}(x,z,e_s,e_f),1,0, $ $\tfrac{1}{\epsilon}g_z(x,z,e_s,e_f), \tfrac{1}{\epsilon} g_{e_f}(x,z,e_s,e_f, \epsilon),  \frac{1}{\epsilon},0\big)$, and 
\begin{align*}
    G(\xi) \coloneqq \left\{ 
    \begin{aligned}
    &G_s(\xi), \quad \xi\in\mathcal{D}_s^\epsilon \setminus \mathcal{D}_f^\epsilon , \\
    &G_f(\xi), \quad \xi\in\mathcal{D}_f^\epsilon \setminus \mathcal{D}_s^\epsilon ,\\
    &\{G_s(\xi),G_f(\xi)\},\quad \xi\in \mathcal{D}_s^\epsilon\cap\mathcal{D}_f^\epsilon .
    \end{aligned}
    \right. 
\end{align*}
The jump maps are defined as $G_s(\xi) \coloneqq (x,h_s(\kappa_s, e_s),0,$ $ \kappa_s + 1, z, e_f,  \tau_f,  \kappa_f)$ and $G_f(\xi) \coloneqq(x, e_s, $ $\tau_s,\kappa_s,  z,  h_f(\kappa_f,$ $ e_f),  0,  \kappa_f + 1 )$, where $G_s$ and $G_f$ corresponds to the transmission of slow and fast signals, respectively. 
%The jump map $G$ is defined such that, at any transmission instance where both the slow and fast transmissions are allowed, i.e. $\xi \in \mathcal{D}_s^\epsilon \cap \mathcal{D}_f^\epsilon$, the trajectory experiences a single jump according to either $G_s$ or $G_f$. Similar to the approach in \cite{abdelrahim2017robust}, this \red{design}\todo{or modelling choice?} choice ensures that the jump map $G$ is outer semicontinuous (OSC) \cite[Definition 5.9]{gosate12}, which is one of the hybrid basic conditions \cite[Assumption 6.5]{gosate12}. 
%The jump map $G$ would not be OSC if $G(\xi)$ were defined as $\{G_s(\xi) \}$ or $\{G_f(\xi) \}$ when $\xi \in \mathcal{D}_s^\epsilon \cap \mathcal{D}_f^\epsilon$.
The set-valued map in the definition of $G$
%, i.e., when $\xi \in \mathcal{D}_s^\epsilon \cap \mathcal{D}_f^\epsilon$, 
is introduced to ensure that $\mathcal{H}_1$ satisfies the hybrid basic conditions \cite[Assumption 6.5]{gosate12}, providing well-posedness of the system. This approach is commonly used when modeling the NCS as hybrid dynamical systems, see \cite{abdelrahim2017robust,wang2015emulation} for more details.

\begin{figure*}[!htp]
	\hrule
    \scriptsize % You can change this to \small \footnotesize, \scriptsize, or \tiny
	\begin{equation}\label{eq:functions}
    \begin{aligned}
		f_x(x,z,e_s,e_f) &\coloneqq 
    \big(
    f_p(x_p,z_p,(k_{c_s}(x_c)+e_{u_s},k_{c_f}(x_c,z_c)+e_{u_f}) ),
    f_c(x_c,z_c,(k_{p_s}(x_p)+e_{y_s},k_{p_f}(x_p,z_p)+e_{y_f}) )
    \big)  \\
    g_z(x,z,e_s,e_f) &\coloneqq 
    \big(
    g_p(x_p,z_p,(k_{c_s}(x_c)+e_{u_s},k_{c_f}(x_c,z_c)+e_{u_f}) ) ,
    g_c(x_c,z_c,(k_{p_s}(x_p)+e_{y_s},k_{p_f}(x_p,z_p)+e_{y_f} ) )
    \big)  \\
    f_{e_s}(x,z,e_s,e_f) &\coloneqq
    \Big(- \tfrac{\partial k_{p_s}(x_p)}{\partial x_p} 
        f_{x,1}(x,z,e_s,e_f), 
    - \tfrac{\partial k_{c_s}(x_c)}{\partial x_c} 
        f_{x,2}(x,z,e_s,e_f)\Big)  \\
        g_{e_f}(x,z,e_s,e_f,\epsilon) &\coloneqq \Big(  -\epsilon \tfrac{\partial k_{p_f}(x_p,z_p)}{\partial x_p}  f_{x,1}(x,z,e_s,e_f)  - \tfrac{\partial k_{p_f}(x_p,z_p)}{\partial z_p} g_{z,1}(x,z,e_s,e_f) , \\
        &\hspace{4.7cm}  - \epsilon \tfrac{\partial k_{c_f}(x_c,z_c)}{\partial x_c}  f_{x,2}(x,z,e_s,e_f) -\tfrac{\partial k_{c_f}(x_c,z_c)}{\partial z_c} g_{z,2}(x,z,e_s,e_f)\Big).  
    \end{aligned}
    \end{equation}
    \begin{equation}\label{eq:functions 2}
    \begin{aligned}
    F_s^y(x,y,e_s,e_f) &\coloneqq \big(f_x(x,y+\overline{H}(x,e_s),e_s, e_f ), f_{e_s}(x,y+\overline{H}(x,e_s),e_s, e_f ),1,0\big) 
    \\
    F_f^y(x,y,e_s,e_f,\epsilon) &\coloneqq \big(
      g_z(x,y+\overline{H}(x,e_s),e_s,e_f)- \epsilon \tfrac{\partial \overline{H}}{\partial \xi_s} F_s^y(x,y,e_s,e_f ),
      g_{e_f}(x,y+\overline{H}(x,e_s),e_s, e_f , \epsilon),1,0 \big)
    \end{aligned}
    \end{equation}
    \normalsize
	\hrule
\end{figure*}

\section{Auxiliary systems}

% To facilitate the forthcoming analysis, we introduce $\mathcal{H}_1$ as the hybrid system with dynamics as per (\ref{eqn:full system}), but with the "patched" flow set defined as 
% $ \mathcal{C}_2^\epsilon \coloneqq \{ \xi \in \mathbb{X} \;|\; \tau_s \in [0, \tau_{\text{mati}}^s] \;\wedge\ \epsilon \tau_f \in  [0, \tau_{\text{mati}}^f]  \}$.
%
%
% \begin{equation*}
%      \mathcal{C}_2^\epsilon \coloneqq \left\{ \xi \in \mathbb{X} \;|\; \tau_s \in [0, \tau_{\text{mati}}^s] \;\wedge\ \epsilon \tau_f \in  [0, \tau_{\text{mati}}^f]  \right\}.
% \end{equation*}
% We note that $\mathcal{H}_1$ \emph{contains} $\mathcal{H}_1$ in the sense that all solutions of $\mathcal{H}_1$ are also solutions to $\mathcal{H}_1$.
%, since $\mathcal{C}_1^\epsilon  \subseteq \mathcal{C}_2^\epsilon$ and they have identical flow map, jump map and jump set. 
% Therefore, using \cite[Proposition 3.32]{gosate12}, we can conclude the stability properties of $\mathcal{H}_1$ by analysing the stability of $\mathcal{H}_1$.
%We also note that if $\mathcal{H}_1$ is initialized at  some $\xi_0 \in \mathcal{C}_1^\epsilon$, its maximal solution will be complete, otherwise it will have a non-complete maximal solution.
%
We adopt a similar approach to the standard singularly perturbed method \cite[Section 11.5]{nonlinear_systems_Khalil} to establish stability properties for $\mathcal{H}_1$, but generalised to hybrid systems. Particularly, we first derive a system $\mathcal{H}_1^y$ by changing the $z$--coordinate of $\mathcal{H}_1$ to $y$--coordinate, where $y$ is defined in \eqref{eqn: map between y and z}, and determine its stability through a \emph{boundary layer} and \emph{reduced system}. 
\subsection{Change of coordinates}

We first derive the \emph{quasi-steady-state} of $\mathcal{H}_1$, under the following assumption.

% \textbf{Standing Assumption 1}\hspace{5pt}\rm\textbf{(SA1)} 
% \textit{For any $\overline{x}\in \mathbb{R}^{n_x}$, $\overline{e}_s\in \mathbb{R}^{n_{e_s}}$ and $\overline{z}\in \mathbb{R}^{n_z}$, equation $ 0 = g_z\left(\bar x,\bar z, \bar e_{s},0\right)$ has a unique real solution $\bar z = \overline{H}(\bar x,  \bar e_{s})$, where $\overline{H}$ is continuously differentiable and $0 =\overline{H}(0, 0)$.}

\begin{sassum}\label{assum:standing-ss} \rm 
\textbf{(SA1)} \it
For any $\overline{x}\in \mathbb{R}^{n_x}$, $\overline{e}_s\in \mathbb{R}^{n_{e_s}}$ and $\overline{z}\in \mathbb{R}^{n_z}$, equation $ 0 = g_z\left(\bar x,\bar z, \bar e_{s},0\right)$ has a unique real solution $\bar z = \overline{H}(\bar x,  \bar e_{s})$, where $\overline{H}$ is continuously differentiable and $0 =\overline{H}(0, 0)$.
\end{sassum}

The \emph{quasi-steady-states} $\bar z$ and $\bar e_{f}$, referring to the equilibrium of the fast states as $\epsilon$ approaches zero, are obtained as follows:
$\bar{e}_f$ is equal to zero, as  for sufficiently high frequency of fast-output transmissions, $e_f$ converges to zero; and $\bar{z}$ corresponds to the unique solution $\bar z = \overline{H}(\bar x,  \bar e_{s})$ as per SA\ref{assum:standing-ss}.
We define the variable $y$ as
\begin{equation}
    y\coloneqq z - \overline{H}(x, e_{s}).
    \label{eqn: map between y and z}
\end{equation}
Then similar to the assumptions in the continuous-time SPSs literature such as \cite{nonlinear_systems_Khalil,christofides1996singular}, SA\ref{assum:standing-ss} guarantees the map \eqref{eqn: map between y and z} to be stability preserving, which means the origin of the $x$-$z$ coordinate is asymptotically stable if and only if the origin of the $x$-$y$ system is asymptotically stable, see \cite[Section 11.5]{nonlinear_systems_Khalil} for more detail.
% \begin{equation}
%     y\coloneqq z - \overline{H}(x, e_{s})
%     \label{eqn: map between y and z}
% \end{equation}
%
% The \emph{quasi-steady-states} $\bar z$ and $\bar e_{f}$, referring to the equilibrium of the fast states as $\epsilon$ approaches zero, are obtained as follows:
% $\bar{e}_f$ is equal to zero, as  for sufficiently high frequency of fast-output transmissions, $e_f$ converges to zero; and
%   $\bar{z}$ corresponds to the unique solution $\bar z = \overline{H}(\bar x,  \bar e_{s})$ as per SA\ref{assum:standing-ss}.
%
Next, to derive $\mathcal{H}_1^y$, we define the full state of $\mathcal{H}_1^y$, namely 
\begin{equation}
    \xi^y \coloneqq (\xi_s, \xi_f) \coloneqq \big((x,e_s,\tau_s, \kappa_s), (y,e_f, \tau_f, \kappa_f)\big),
    \label{eqn: definition of xi_s and xi_f}
\end{equation}
where $\xi^y \in \mathbb{X}$, $\xi_s \in \mathbb{X}^{s} \coloneqq \mathbb{R}^{n_x} \times \mathbb{R}^{n_{e_s}} \times \mathbb{R}_{\geq 0} \times \mathbb{Z}_{\geq 0}$ and $\xi_f \in  \mathbb{X}^{f} \coloneqq \mathbb{R}^{n_z} \times \mathbb{R}^{n_{e_f}} \times \mathbb{R}_{\geq 0} \times \mathbb{Z}_{\geq 0}$. 
When a slow variable is transmitted at $t_k^s \in \mathcal{T}^s$, $e_s$ updates according to $h_s$, then by the definition of $y$ in \eqref{eqn: map between y and z}, we know at each slow transmission, the value of $y$ updates according to
\begin{equation}
    \begin{aligned}
        y^+ &= z^+ - \overline{H}(x^+, e_{s}^+)= z - \overline{H}(x, h_s(\kappa_s, e_s)) \\
        %&= z - \overline{H}(x, h_s(\kappa_s, e_s)) \\
        &= y + \overline{H}(x, e_s) - \overline{H}(x, h_s(\kappa_s, e_s)) \\
        & \eqqcolon h_y(\kappa_s,x,e_s,y). 
    \end{aligned}
    \label{eqn: Jump of y at slow transmission}
\end{equation}
Then, $\mathcal{H}_1^y$ is given by
\begin{equation}
    \mathcal{H}_1^y:\left\{
\begin{aligned}
    \dot{\xi}^y &= F^y(\xi^y, \epsilon),\ \xi^y \in \mathcal{C}_2^{y,\epsilon}, \\
    {\xi^y}^+ &\in G^y(\xi^y), \ \xi^y\in \mathcal{D}_s^{y,\epsilon} \cup \mathcal{D}_f^{y,\epsilon},
\end{aligned}
    \right.
    \label{eqn: H_2^y}
\end{equation}
where $F^y(\xi^y, \epsilon) = \big(F_s^y(x,y,e_s,e_f), \tfrac{1}{\epsilon}F_f^y(x,y, $     $e_s,e_f,\epsilon)\big)$, with $F_s^y$ and $F_f^y$ from \eqref{eq:functions 2}. 
% $F_s^y(x,y,e_s,e_f) \coloneqq 
% \big(f_x(x,y+\overline{H}(x,e_s),e_s, e_f ), f_{e_s}(x,y+\overline{H}(x,e_s),e_s, e_f ),1,0\big) $, 
% $ F_f^y(x,y,e_s,e_f,\epsilon) \coloneqq \big(
%       \epsilon \tfrac{\partial y}{\partial t},
%       g_{e_f}(x,y+\overline{H}(x,e_s),e_s, e_f , \epsilon),1,0 \big)$
% and 
% $\epsilon \tfrac{\partial y}{\partial t} = g_z(x,y+\overline{H}(x,e_s),e_s,e_f)- \epsilon \tfrac{\partial \overline{H}}{\partial \xi_s} F_s^y(x,y,e_s,e_f ) $. 
The jump map $G^y$ is given by
\begin{equation}
\begin{aligned}
    G^y(\xi^y) \coloneqq \left\{ 
    \begin{aligned}
    &G_s^y(\xi^y),  \;\xi^y\in\mathcal{D}_s^{y,\epsilon} \setminus \mathcal{D}_f^{y,\epsilon} , \\
    &G_f^y(\xi^y),  \;  \xi^y \in\mathcal{D}_f^{y,\epsilon} \setminus \mathcal{D}_s^{y,\epsilon} ,\\
    &\{G_s^y(\xi^y),G_f^y(\xi^y)\}, \; \xi^y\in \mathcal{D}_s^{y,\epsilon} \cap \mathcal{D}_f^{y,\epsilon},
    \end{aligned}
    \right. 
\end{aligned}
\label{eqn: G^y}
\end{equation}
with $G_s^y(\xi_y) \coloneqq \big(x, h_s(\kappa_s, e_s), 0, \kappa_s + 1, h_y(\kappa_s,x,e_s,y),  e_f,$ $ \tau_f, \kappa_f \big)$; $G_f^y(\xi_y) \coloneqq \big(x,e_s, \tau_s, \kappa_s , y, h_f(\kappa_f, e_f), 0, \kappa_f + 1 \big)$.

For analysis purposes, we write $\tau_{\text{mati}}^f = \epsilon T^*$ with $T^* \in \mathbb{R}_{>0}$ independent of $\epsilon$. We also write $\tau_{\text{miati}}^f = a\tau_{\text{mati}}^f$ for some $a \in(0,\tfrac{1}{2}] $, which satisfies the inequality \eqref{eqn: condition on miati^f}. Since $\epsilon > 0$, $\epsilon \tau_f \in [\tau_{\text{miati}}^f, \tau_{\text{mati}}^f]$ is equivalent to $\tau_f \in [aT^*,T^*]$. Then the jump and flow sets in \eqref{eqn: H_2^y} are defined by
$\mathcal{D}_s^{y,\epsilon} \coloneqq  \{\xi^y \in \mathbb{X} \; | \; \tau_s \in [\tau_{\text{miati}}^s, \tau_{\text{mati}}^s] \wedge \tau_f \in  [aT^*, (1-a)T^*] \}$, 
$\mathcal{D}_f^{y,\epsilon} \coloneqq \{\xi^y \in \mathbb{X} \; | \; \tau_s \in [\epsilon aT^*, \tau_{\text{mati}}^s-\epsilon aT^*]  \wedge  \tau_f \in  [aT^*, T^*]  \}$ 
and
%$\mathcal{C}_2^{y,\epsilon} \coloneqq  \{\xi^y \in \mathbb{X} \; | \; \tau_s \in [0, \tau_{\text{mati}}^s] \wedge  \tau_f \in  [0, T^*] \}$.
%
$\mathcal{C}_1^{y,\epsilon} \coloneqq 
        \mathcal{D}_s^{y,\epsilon} \cup \mathcal{D}_f^{y,\epsilon} \cup \mathcal{C}_{1,a}^{y,\epsilon} \cup \mathcal{C}_{1,b}^{y,\epsilon}$,
with $\mathcal{C}_{1,a}^{y,\epsilon} \coloneqq \{ \xi^y \in \mathbb{X} \ | \ \tau_s \in [0, \epsilon a T^*]  \wedge \epsilon \tau_f \in [0,\tau_s + \epsilon T^* - \epsilon a T^*] \} $ and 
$\mathcal{C}_{1,b}^{y,\epsilon} \coloneqq \{ \xi^y \in \mathbb{X} \;|\; \tau_s \in [\epsilon a T^*, \epsilon \tau_f + \tau_{\text{mati}}^s - \epsilon a T^*]  \wedge \epsilon \tau_f \in  [0, \epsilon a T^*]  \}$.

%
% \begin{align*}
%     \mathcal{D}_s^{y,\epsilon} \coloneqq & \{\xi^y \in \mathbb{X} \; | \; \tau_s \in [\tau_{\text{miati}}^s, \tau_{\text{mati}}^s] 
%         \\ & \qquad \qquad \qquad \qquad \quad \wedge  \tau_f \in  [aT^*, (1-a)T^*] \},
%     \\ 
%     \mathcal{D}_f^{y,\epsilon} \coloneqq &\{\xi^y \in \mathbb{X} \; | \; \tau_s \in [\epsilon aT^*, \tau_{\text{mati}}^s-\epsilon aT^*]  \\
%     & \qquad \qquad \qquad \qquad \qquad \qquad   \wedge  \tau_f \in  [aT^*, T^*]  \},
%     \\
%     \mathcal{C}_2^{y,\epsilon} \coloneqq & 
%         \{\xi^y \in \mathbb{X} \; | \; \tau_s \in [0, \tau_{\text{mati}}^s] \wedge  \tau_f \in  [0, T^*] \}.
% \end{align*}
%
We have changed the coordinate from $z$ to $y$, and we are now ready to derive the reduced system $\mathcal{H}_r$ and boundary layer system $\mathcal{H}_{bl}$ associated with $\mathcal{H}_1^y$.

\subsection{Boundary layer system and reduced system of $\mathcal{H}_1$}
 We define the fast time scale $\sigma \coloneqq \tfrac{t-t_0}{\epsilon}$, where we can assume $t_0 = 0$ as the system is time invariant. Then we have $\tfrac{\partial}{\partial \sigma} = \epsilon \tfrac{\partial}{\partial t}$. We set $\epsilon = 0$ for system \eqref{eqn: H_2^y}, then $\mathcal{C}_1^{y,0}$, which corresponds to $\mathcal{C}_1^{y,\epsilon}$ with $\epsilon = 0$, is given by $\mathcal{C}_1^{y,0} \coloneqq \{ \xi^y \in \mathbb{X} \ | \ \tau_s \in [0, \tau_\text{mati}^s]  \wedge \tau_f \in [0,  T^*] \}  $, and $\mathcal{D}_s^{y,0}$, $\mathcal{D}_f^{y,0}$ are derived accordingly. In the perspective of fast dynamics, the slow dynamics are now frozen. Meanwhile, the jump and flow sets of $\mathcal{H}_{bl}$ contain the condition $\tau_{s}\in [0, \tau_{\text{mati}}^s]$, which is always satisfied. Therefore, the jumps and flows of $\mathcal{H}_{bl}$ are only determined by $\tau_{f}$. We thus write
\begin{equation}
    \mathcal{H}_{bl}\! : \! \left\{
\begin{aligned}
    (\tfrac{\partial \xi_s}{\partial \sigma}, \tfrac{\partial \xi_f}{\partial \sigma} ) &= (\mathbf{0}_{n_{\xi_s}\! \times 1}, F_f^y(x,y,e_s,e_f,0) ), \xi^y \! \in \mathcal{C}_{1,bl}^{y,0}, \\
    {\xi^y}^+  &=   G_f^y(\xi^y), \qquad \qquad \qquad \qquad \, \xi^y \! \in \mathcal{D}_f^{y,0},
\end{aligned}
    \right.
    \label{eqn: H_bl}
\end{equation}
where $\mathcal{C}_{2,bl}^{y,0} \coloneqq \{\xi^y \in  \mathbb{X} \ | \ \tau_f \in [0, T^*]\}$ and $\mathcal{D}_f^{y,0}\coloneqq \{\xi^y \in  \mathbb{X} \ | \ \tau_f \in [aT^*, T^*] \}$.

From the perspective of $\mathcal{H}_r$ (i.e., slow dynamics), the fast dynamics evolve infinitely fast. Therefore, for any $\tau_s \in [0, \tau_{\text{mati}}^s] $, the waiting time for the condition $\tau_f \in [aT^*, T^*]$ in the jump set to be satisfied approaches to zero, and the flows and jumps of $\mathcal{H}_r$ are essentially determined only by $\tau_s$. 
%
%We assume $\mathcal{H}_{bl}$ satisfies an asymptotic stability property at its quasi-steady state, which we formalise in the sequel. 
Moreover, we have $y=0$ and $e_f = 0$ in $\mathcal{H}_r$, that is 
\begin{equation}
    \mathcal{H}_{r}:\left\{
\begin{aligned}
    \dot \xi_s &= F_s^y(x,0,e_s, 0) , \quad \xi^y \in \mathcal{C}_{1,r}^{y,0}, \\
    \xi_s^+  &=   (x, h_s(\kappa_s, e_s), 0, \kappa_s + 1), \  \xi^y\in \mathcal{D}_s^{y,0},
\end{aligned}
    \right.
    \label{eqn: H_r}
\end{equation}
where $\mathcal{C}_{1,r}^{y,0} \coloneqq \{\xi^y \in  \mathbb{X} \ | \ \tau_s \in [0, \tau_{\text{mati}}^s]\}$ and $\mathcal{D}_s^{y,0}\coloneqq \{\xi^y \in  \mathbb{X} \ | \ \tau_s \in [\tau_{\text{miati}}^s, \tau_{\text{mati}}^s] \}$.

\section{Emulation design framework} \label{Section Emulation design framework}
% In this section, we first outline the assumptions and preliminaries used to ensure semi-global practical asymptotic stability of $\mathcal{H}_1$, and we will then present the conditions that guarantee UGAS and UGES of $\mathcal{H}_1$.

This section presents the main results that provide the framework for emulation design. The first step is to design a controller that making the reduced system and boundary-layer system robust with respect to network induced error by satisfying \eqref{eqn: NCS Vs flow} and \eqref{eqn: NCS Vf flow} in Assumptions \ref{Assumption reduced model} and \ref{Assumption boundary layer system}, respectively. Next, we select UGAS protocols for slow and fast transmissions and verify the growth conditions on error dynamics, i.e., \eqref{eqn: NCS Ws dot} and \eqref{eqn: NCS Wf dot} in Assumptions \ref{Assumption reduced model} and \ref{Assumption boundary layer system}. Then, by checking interconnection condition during flow (Assumption \ref{Assumption interconnection}) and at slow transmissions (Assumption \ref{Assumption Vf at slow transmission}), as well as a mild assumption, we guarantee semi-global practical asymptotic stability given $\tau_{\text{mati}}^s$, $\tau_{\text{mati}}^f$ and $\epsilon$ are sufficiently small.
Finally, we present additional conditions that guarantee UGAS and UGES of $\mathcal{H}_1$ in section \ref{Section UGES and UGAS}.

\subsection{Semi-global practical asymptotic stability}
Assumptions \ref{Assumption reduced model} and \ref{Assumption boundary layer system} below provide sufficient conditions to guarantee asymptotic stability properties for $\mathcal{H}_r$ and $\mathcal{H}_{bl}$, respectively, which align with those commonly encountered in the NCS literature, see \cite{carnevale_stability,SPNCS}.
\begin{assum}
 There exist a function $W_s: \mathbb{Z}_{\geq 0}\times \mathbb{R}^{n_{e_s}} \to \mathbb{R}_{\geq 0}$ that is locally Lipschitz in its second argument uniformly in its first argument, a continuous function $H_s : \mathbb{R}^{n_x}\times\mathbb{R}^{n_{e_s}} \rightarrow \mathbb{R}_{\geq 0} $, $\mathcal{K}_{\infty}$-functions $ \underline{\alpha}_{W_s},\overline{\alpha }_{W_s}  $, constants $\lambda_s \in [0,1)$ and $ L_s \geq 0$ such that, for all $ \kappa_s \in \mathbb{Z}_{\geq 0}$ and $e_s \in \mathbb{R}^{n_{e_s}}$, the following properties hold:
\begin{align}
    \underline{\alpha}_{W_s}\left(\left|  {e_s}  \right|\right) \leq {W_s}(k_s, {e_s}) \leq \overline{\alpha }_{W_s}\left(\left|   {e_s}  \right|\right) ,\label{eqn: NCS assumption Ws sandwich bound}
    \\
    {W_s}(\kappa_s + 1, h_s(\kappa_s, e_s)) \leq \lambda_s {W_s}(\kappa_s, {e_s}). \label{eqn: NCS assumption Ws jump}
\end{align}
For all $x \in \mathbb{R}^{n_x}, \kappa_s \in \mathbb{Z}_{\geq 0}$ and almost all ${e_s} \in \mathbb{R}^{n_{e_s}}$,
% \begin{equation}
%     \begin{aligned}
%     &\left< \tfrac{\partial {W_s}(\kappa_s,{e_s})}{\partial {e_s}}, f_{e_s}(x,\overline{H}(x,e_s),e_s, 0)\right> \
%     \\
%     & \phantom{aaaaaaaaaaaaaaa}   \leq  L_s {W_s}(\kappa_s, e_s)  + H_s(x,e_s).
%     \end{aligned}
%     \label{eqn: NCS Ws dot}
% \end{equation} 
\begin{multline}
    \left< \tfrac{\partial {W_s}(\kappa_s,{e_s})}{\partial {e_s}}, f_{e_s}(x,\overline{H}(x,e_s),e_s, 0)\right> \
    \\
    \phantom{aaaaaaaaaaaaaaa}   \leq  L_s {W_s}(\kappa_s, e_s)  + H_s(x,e_s).
    \label{eqn: NCS Ws dot}
\end{multline}
Moreover, there exist a locally Lipschitz, positive definite and radially unbounded function ${V_s}: \mathbb{R}^{n_x} \to \mathbb{R}_{\geq 0}$, positive definite function $\rho_s $, and real number $\gamma_s > 0$, such that for all $e_s \in \mathbb{R}^{n_{e_s}}$, all $\kappa_s \in \mathbb{Z}_{\geq 0}$, and almost all $x\in \mathbb{R}^{n_x}$, the following inequality holds
\begin{multline}
    \left< \tfrac{\partial {V_s}(x)}{\partial x},f_x(x,\overline{H}(x,e_s),e_s, 0) \right> \leq - \rho_s(|x|) 
    \\
     - \rho_s\left(W_s(\kappa_s, e_s)\right) - H_s^2(x,e_s) + \gamma_s^2 W_s^2(\kappa_s, e_s). 
    \label{eqn: NCS Vs flow}
\end{multline}
\label{Assumption reduced model}
\end{assum}
\vspace{-0.5cm}
\begin{assum}
There exist a function ${W_f}: \mathbb{Z}_{\geq 0}\times \mathbb{R}^{n_{e_f}} \to \mathbb{R}_{\geq 0}$ that is locally Lipschitz in its second argument uniformly in its first argument, a continuous function $H_f : \mathbb{R}^{n_x}\times\mathbb{R}^{n_{e_f}} \rightarrow \mathbb{R}_{\geq 0} $, $\mathcal{K}_{\infty}$-functions $ \underline{\alpha}_{W_f},\overline{\alpha }_{W_f}  $, constants $\lambda_f \in [0,1)$ and $ L_f \geq 0$ such that, for all $ \kappa_f \in \mathbb{Z}_{\geq 0}$ and $e_f \in \mathbb{R}^{n_{e_f}}$, the following properties hold:
\begin{align}
   \underline{\alpha}_{W_f}\left(\left|  {e_f}  \right|\right) \leq {W_f}(k_f, {e_f}) \leq \overline{\alpha }_{W_f}\left(\left|   {e_f}  \right|\right) ,\label{eqn: NCS assumption Wf sandwich bound}
    \\
    {W_f}(\kappa_f + 1, h_f(\kappa_f, e_f)) \leq \lambda_f {W_f}(\kappa_f, {e_f}). \label{eqn: NCS assumption Wf jump}
\end{align}
For all $x \in \mathbb{R}^{n_x}, \kappa_f \in \mathbb{Z}_{\geq 0}$ and almost all ${e_f} \in \mathbb{R}^{n_{e_f}}$,
\begin{multline}
   \left< \tfrac{\partial {W_f}(\kappa_f,{e_f})}{\partial {e_f}}, g_{e_f}(x,y+  \overline{H}(x, e_s),e_s,e_f, 0)\right>
    \\
     \leq  L_f {W_f}(\kappa_f, e_f) + H_f(y,e_f).
    \label{eqn: NCS Wf dot}
\end{multline}  
Moreover, there exist a locally Lipschitz function ${V_f}: \mathbb{R}^{n_x}\times \mathbb{R}^{n_z} \to \mathbb{R}_{\geq 0}$, $\mathcal{K}_\infty$-functions $\underline{\alpha}_{V_f}$, $\overline{\alpha}_{V_f}$, such that for all $x \in \mathbb{R}^{n_x}$ and $y\in \mathbb{R}^{n_z}$, the following inequality holds. 
\begin{equation}
    \underline{\alpha}_{V_f}\left(\left|  y \right|\right) \leq {V_f}(x, y) \leq \overline{\alpha }_{V_f}\left(\left|   y  \right|\right). \label{eqn: Vf sandwich bound}
\end{equation}
At the same time, there exist positive definite function $\rho_f $, and real number $\gamma_f > 0$ such that for all $e_s \in \mathbb{R}^{n_{e_s}}$, $e_f \in \mathbb{R}^{n_{e_f}}$, all $\kappa_f \in \mathbb{Z}_{\geq 0}$, all $x\in \mathbb{R}^{n_x}$, and almost all $y\in \mathbb{R}^{n_z}$, the following inequality holds.
\begin{multline}
        \left< \tfrac{\partial {V_f}(x,y)}{\partial y},g_z(x,y+ \overline{H}(x, e_s),e_s,e_f)  \right>\leq - \rho_f(|y|) 
            \\
             - \rho_f\left(W_f(\kappa_f, e_f)\right) - H_f^2(y,e_f) 
             + \gamma_f^2 W_f^2(\kappa_f,e_f). \!\!\!\!
             \label{eqn: NCS Vf flow}
\end{multline}
\vspace{-0.5cm}
\label{Assumption boundary layer system}
\end{assum}
\vspace{-0.5cm}
In Assumption \ref{Assumption reduced model} (similarly with Assumption \ref{Assumption boundary layer system}), conditions (\ref{eqn: NCS assumption Ws sandwich bound}) and (\ref{eqn: NCS assumption Ws jump}) relate to UGAS protocols and are satisfied by sampled-data systems and NCSs with RR, TOD, etc; for more details, see \cite{dragan_stability}. Inequality \eqref{eqn: NCS Ws dot} bounds the growth of $e_s$ during flow, and (\ref{eqn: NCS Vs flow}) relates to the $\mathcal{L}_2$-stability of $\mathcal{H}_r$ from $W_s$ to $H_s$, which is typically ensured at the first stage of emulation.
According to \cite{carnevale_stability}, Assumption \ref{Assumption reduced model} implies there exists a $\tau_{\text{mati}}^{s,*} > 0$ such that for all $0<\tau_{\text{mati}}^{s}\leq \tau_{\text{mati}}^{s,*}$, the set $\{\xi^y \in \mathbb{X} | x = 0 \wedge e_s = 0 \}$ is UGAS for $\mathcal{H}_r$.
See \cite{dragan_stability} for more details on finding Lyapunov functions to satisfy Assumptions \ref{Assumption reduced model} and \ref{Assumption boundary layer system}. 
We next provide a lemma as a preliminary to the main result. %We will use $U^\circ$ to denote the Clarke generalized derivative of a locally Lipschitz function $U$ \cyan{\cite[Eqn. (20)]{teel2000assigning}.}

Recalling that $\xi^s = (x, e_s, \tau_s, \kappa_s)$ and $\xi^f = (y, e_f, \tau_f, \kappa_f)$, we define Lyapunov functions ${U_s}:  \mathbb{X}^{s} \to \mathbb{R}_{\geq 0}$ and $U_f :\mathbb{X}^{s} \times\mathbb{X}^{f} \to \mathbb{R}_{\geq 0}$ as \cite[Eqn. (25)]{carnevale_stability} 
\begin{subequations}
    \begin{align}
        U_s(\xi_s) &= V_s(x) + \gamma_s \phi_s(\tau_s) W_s^2(\kappa_s, e_s), \label{eqn: definition of U_s}\\
        U_f(\xi_s,\xi_f) &=  V_f(x,y) + \gamma_f \phi_f(\tau_f) W_f^2(\kappa_f, e_f), \label{eqn: definition of U_f}%
    \end{align}
    \label{eqn: Us and Uf}%
\end{subequations}
where $\dot \phi_\star = -2L_\star \phi_\star - \gamma_\star (\phi_\star^2 + 1)$, $\phi_\star(0) = 1/\lambda_\star^*$, $\lambda_\star^* \in (\lambda_\star, 1)$, $\star \in \{s,f\}$. Note that by abuse of notation, we write $\dot{\phi}_s = \tfrac{d \phi_s}{d \tau_s}$, $\dot{\phi}_f = \tfrac{d \phi_f}{d \tau_f}$ and $U_f(\xi^y) = U_f(\xi_s, \xi_f)$.

We define the nonlinear mapping $T: \mathbb{R}_{\geq0}\times(0,1)\times\mathbb{R}_{>0} \rightarrow \mathbb{R}$ for the upcoming lemma. For any $L > 0$, $\lambda \in (0,1)$ and $\gamma > 0$, 
\begin{equation*}
    T(L,\gamma,\lambda) \coloneqq
    \begin{cases}
    \tfrac{1}{Lr}\tan^{-1}\bigg(\tfrac{r(1-\lambda)}{2\tfrac{\lambda}{1+\lambda}\big(\tfrac{\gamma}{L}-1\big)+1+\lambda} \bigg), \;\;\;  \gamma > L
    \\
    \tfrac{1}{L} \tfrac{1-\lambda}{1+\lambda}, \qquad\qquad\qquad\qquad\qquad \;\;\; \gamma = L
    \\
    \tfrac{1}{Lr}\tanh^{-1}\bigg(\tfrac{r(1-\lambda)}{2\tfrac{\lambda}{1+\lambda}\left(\tfrac{\gamma}{L}-1\right)+1+\lambda} \bigg), \; \gamma < L ,
    \end{cases}   
\end{equation*}
where $r \coloneqq \sqrt{|\left(\tfrac{\gamma}{L}\right)^2 -1|}$. For $L = 0$ and any $\lambda \in (0,1)$ and $\gamma > 0$, this nonlinear mapping becomes
$$T(0,\gamma, \lambda) = \tfrac{1}{\gamma} \big(\tan^{-1}(\tfrac{1}{\lambda}) - \tan^{-1}(\lambda) \big).$$ 
\begin{lem}\label{Lemma MATI}
Suppose Assumptions \ref{Assumption reduced model} and \ref{Assumption boundary layer system} hold. 
%
% For any $L \geq 0$, $\lambda \in (0,1)$ and $\gamma > 0$, we define the following nonlinear mapping:
% \begin{equation*}
%     T(L,\gamma,\lambda) \coloneqq
%     \begin{cases}
%     \tfrac{1}{Lr}\tan^{-1}\bigg(\tfrac{r(1-\lambda)}{2\tfrac{\lambda}{1+\lambda}\big(\tfrac{\gamma}{L}-1\big)+1+\lambda} \bigg), \;\;\;  \gamma > L
%     \\
%     \tfrac{1}{L} \tfrac{1-\lambda}{1+\lambda}, \qquad\qquad\qquad\qquad\qquad \;\;\; \gamma = L
%     \\
%     \tfrac{1}{Lr}\tanh^{-1}\bigg(\tfrac{r(1-\lambda)}{2\tfrac{\lambda}{1+\lambda}\left(\tfrac{\gamma}{L}-1\right)+1+\lambda} \bigg), \; \gamma < L ,
%     \end{cases}   
% \end{equation*} 
% where $r \coloneqq \sqrt{|\left(\tfrac{\gamma}{L}\right)^2 -1|}$. 
%
Let $(L_s, \gamma_s, \lambda_s)$ and $(L_f, \gamma_f, \lambda_f)$ come from Assumption \ref{Assumption reduced model} and \ref{Assumption boundary layer system}, respectively, and $U_s$ and $U_f$ come from \eqref{eqn: Us and Uf} with some $\lambda_s^*\in (\lambda_s,1)$ and $\lambda_f^*\in (\lambda_f,1)$.
For all $\tau_{\text{mati}}^s \leq T(L_s,\gamma_s, \lambda_s^*)$ and $T^* \leq T(L_f, \gamma_f, \lambda_f^*)$, there exist %Lyapunov functions ${U_s}: \mathbb{R}^{n_{\xi_s}} \to \mathbb{R}_{\geq 0} , U_f : \mathbb{R}^{n_{\xi_s}} \times \mathbb{R}^{n_{\xi_f}} \to \mathbb{R}_{\geq 0}$,
$\mathcal{K}_\infty$-functions $\underline{\alpha}_{U_s}, \overline{\alpha}_{U_s},\underline{\alpha}_{U_f}, \overline{\alpha}_{U_f}$, continuous positive definite functions $\psi_s, \psi_f$ and positive constants $a_s,a_f$ such that 
(\ref{eqn: Us sandwich bound}) holds for all $\xi_s \in \mathcal{C}_{2,r}^{y,0} \cup \mathcal{D}_s^{y,0}$, (\ref{eqn: Us flow}) holds for all $\xi_s \in \mathcal{C}_{2,r}^{y,0}$, (\ref{eqn: Us jump}) holds for all $\xi_s \in \mathcal{D}_s^{y,0}$, (\ref{eqn: Uf sandwich bound}) holds for all $\xi_f \in \mathcal{C}_{2,bl}^{y,0} \cup \mathcal{D}_f^{y,0}$, (\ref{eqn: Uf flow}) holds for all $\xi_f \in \mathcal{C}_{2,bl}^{y,0}$ and (\ref{eqn: Uf jump}) holds for all $\xi_f \in \mathcal{D}_f^{y,0}$,
\begin{subequations}
     \begin{align}
        \underline{\alpha}_{U_s}\left(\left| ( x , e_s ) \right|\right) \leq {U_s}(\xi_s) &\leq \overline{\alpha}_{U_s}\left(\left| (x, e_s) \right|\right),
        \label{eqn: Us sandwich bound}
        % \\
        % \tfrac{\partial {U_s}(\xi_s)}{\partial \xi_s} F_s(x,0,e_s, 0) &\leq -a_s \psi_s^2\left(\left| (x, e_s) \right|\right) , \label{eqn: Us flow}
        \\
        U_s^\circ(\xi_s; F_s^y(x,0,e_s, 0)) &\leq -a_s \psi_s^2\left(\left| (x, e_s) \right|\right) ,\label{eqn: Us flow}
        \\
        {U_s}((x, h_s(\kappa_s, e_s), 0, & \kappa_s + 1))  \leq {U_s}(\xi_s),   \label{eqn: Us jump}
     \end{align}
     \label{eqn: Us}%
\end{subequations}
\vspace{-0.7cm}
\begin{subequations}
     \begin{gather}
        \underline{\alpha}_{U_f}\left(\left| (y, e_f) \right|\right)\leq {U_f}(\xi_s,\xi_f) \leq \overline{\alpha}_{U_f}\left(\left| (y, e_f) \right|\right) ,
        \label{eqn: Uf sandwich bound}%
        % \\
        % \tfrac{\partial {U_f}(\xi_s, \xi_f)}{\partial \xi_f} F_f(x,y,e_s, e_f, 0) &\leq -a_f \psi_f^2 \left(\left| (y, e_f) \right|\right),
        % \label{eqn: Uf flow}
        \\
        \begin{aligned}
        U_f^\circ \big((\xi_s,\xi_f); (\mathbf{0}_{n_{\xi_s \times 1}}, F_f^y(x,&y,e_s,e_f, 0))\big) 
            \\
            &\leq -a_f \psi_f^2 \left(\left| (y, e_f) \right|\right),
        \end{aligned}
        \label{eqn: Uf flow}%
        \\
        {U_f}(G_f^y(\xi^y))  \leq {U_f}(\xi_s,\xi_f).
        \label{eqn: Uf jump}
     \end{gather}
     \label{eqn: Uf}%
\end{subequations}
\end{lem}
\vspace{-0.5cm}
\textbf{Proof:} The proof of Lemma \ref{Lemma MATI} follows similarly to \cite[Theorem 1]{sampled_data_system} and is therefore omitted.

% 
% By abuse of notation, we write $U_f(\xi^y) = U_f(\xi_s, \xi_f)$. The functions $U_s$ and $U_f$ typically have the form \cite[Eqn. (25)]{carnevale_stability}
% \begin{subequations}
%     \begin{align}
%         U_s(\xi_s) &= V_s(x) + \gamma_s \phi_s(\tau_s) W_s^2(\kappa_s, e_s) \\
%         U_f(\xi_s,\xi_f) &=  V_f(x,y) + \gamma_f \phi_f(\tau) W_f^2(\kappa_f, e_f) 
%     \end{align}
% \end{subequations}
% where $\dot \phi_\star = -2L_\star \phi_\star - \gamma_\star (\phi_\star^2 + 1)$, $\phi_\star(0) = 1/\lambda_\star^*$,  $\star \in \{s,f\}$.
%
%\cyan{From \eqref{eqn: Us} and the fact $t \geq \tau_{\text{miati}}^s j_k^s -  \tau_{\text{miati}}^s$, we can show $\mathcal{H}_r$ is Uniformly Globally pre-Asymptotically Stable (UGpAS) by utilizing \cite[Proposition 3.27]{_systems}. Similarly, we can also show $\mathcal{H}_{bl}$ is UGpAS.} 
%
Lemma \ref{Lemma MATI} asserts that, under Assumptions  \ref{Assumption reduced model} and \ref{Assumption boundary layer system}, we can establish upper bounds on $\tau_{\text{mati}}^s$ and $T^*$ in a manner such that, when both bounds are met, we can construct Lyapunov functions $U_s$ and $U_f$ that guarantee stability properties for $\mathcal{H}_r$ and $\mathcal{H}_{bl}$, respectively. These Lyapunov functions will play a crucial role in the proof of our main result (namely Theorem \ref{Theorem H_1} below), since we will conclude stability property of $\mathcal{H}_1^y$ by considering $\mathcal{H}_r$, $\mathcal{H}_{bl}$, and their interconnection induced by nonzero $\epsilon$. 

Assumption \ref{Assumption interconnection} specifies the \emph{interconnection condition} between the 
slow and fast dynamics during flow, analogous to the continuous-time case as described in \cite[pp. 451]{nonlinear_systems_Khalil}.% and \cite{khorasani1985asymptotic}.% On the other hand, Assumption \ref{Assumption U_f slow jump} considers the interconnection during jumps.
\begin{assum}
    Given a set $\widetilde {\mathcal{C}} \in \mathbb{X}$, for any $\Delta_1$, $\nu_1>0$, there exist $b_1$, $b_2$, $b_3 \geq 0$, such that $|\xi^y|_{\mathcal{E}^y} \leq \Delta_1$ implies
    \begin{align*}
        &\begin{aligned}
            &\left < \tfrac{\partial {U_s}}{\partial \xi_s}, F_s^y(x,y,e_s,e_f) - F_s^y(x,0,e_s,0)  \right> \leq
            \\
            &\phantom{aaaaaaaaaaaaaaa} b_1 \psi_s\left(\left| (x, e_s) \right|\right) \psi_f\left(\left| (y, e_f) \right|\right) + \nu_1, 
        \end{aligned}
        \\
        &\begin{aligned}
            &\Big< \tfrac{\partial {U_f}}{\partial \xi_s} - \tfrac{\partial {U_f}}{\partial y} \tfrac{\partial \overline{H}}{\partial \xi_s} - \tfrac{\partial {U_f}}{\partial e_f} \tfrac{\partial \tilde k}{\partial \xi_s} ,  F_s^y(x,y,e_s,e_f) \Big> \leq
            \\
            & \phantom{aaaaa} b_2 \psi_s\left(\left| (x, e_s) \right|\right) \psi_f\left(\left| (y, e_f) \right|\right) + b_3 \psi_f^2\left(\left| (y, e_f) \right|\right) + \nu_1
        \end{aligned}
    \end{align*}
%  
% \begin{subequations}
%     \begin{align}
%         &\begin{aligned}
%             &\left < \tfrac{\partial {U_s}}{\partial \xi_s}, F_s^y(x,y,e_s,e_f) - F_s^y(x,0,e_s,0)  \right> \leq
%             \\
%             &\phantom{aaaaaaaaaaaa} b_1 \psi_s\left(\left| (x, e_s) \right|\right) \psi_f\left(\left| (y, e_f) \right|\right) + \nu_1, 
%         \end{aligned}\label{eqn: SPNCS interconnection 1}
%         \\
%         &\begin{aligned}
%             &\Big< \tfrac{\partial {U_f}}{\partial \xi_s} - \tfrac{\partial {U_f}}{\partial y} \tfrac{\partial \overline{H}}{\partial \xi_s} - \tfrac{\partial {U_f}}{\partial e_f} \tfrac{\partial \tilde k}{\partial \xi_s} ,  F_s^y(x,y,e_s,e_f) \Big> \leq
%             \\
%             & b_2 \psi_s\left(\left| (x, e_s) \right|\right) \psi_f\left(\left| (y, e_f) \right|\right) + b_3 \psi_f^2\left(\left| (y, e_f) \right|\right) + \nu_1
%         \end{aligned}\label{eqn: SPNCS interconnection 2}
%         % \\
%         % &\begin{aligned}
%         %     &\Big< \tfrac{\partial {U_f}}{\partial \xi_f}, F_f^y(x,y,e_s, e_f,\epsilon) - F_f^y(x,y,e_s,e_f,0)\Big> \leq
%         %     \\
%         %      & \phantom{aa}\; \epsilon b_4 \psi_s\left(\left| (x, e_s) \right|\right) \psi_f\left(\left| (y, e_f) \right|\right) 
%         %          + \epsilon b_5 \psi_f^2\left(\left| (y, e_f) \right|\right)
%         % \end{aligned}\label{eqn: SPNCS interconnection 3}
%     \end{align}
%     \label{eqn: SPNCS interconnections}%
% \end{subequations}
hold for almost all $\xi^y \in \widetilde {\mathcal{C}}$, where $\tilde k(x,z) = (k_{pf}(x_p,z_p), $ $ k_{cf}(x_c,z_c))$.
    \label{Assumption interconnection}
\end{assum}

At each slow transmission, there is a potential increase in $V_f$ due to \eqref{eqn: Jump of y at slow transmission}, we bound this jump of $V_f$ using the following assumption, which is adapted from \cite[Assumption 5]{Romain_ETC}.
\begin{assum}
    There exist $\lambda_1$, $\lambda_2\geq0$ such that for all $\xi^y \in \mathbb{X}$, we have
    $V_f(x,h_y(x,e_s,y)) \leq  V_f(x,y) +  \lambda_1 W_s^2(\kappa_s,e_s) 
            + \lambda_2 \sqrt{W_s^2(\kappa_s,e_s) V_f(x,y)}$.
    %
    % \begin{equation} \begin{aligned}
    %     V_f(x,h_y(x,e_s,y)) \leq & V_f(x,y) +  \lambda_1 W_s^2(\kappa_s,e_s) 
    %         \\ &+ \lambda_2 \sqrt{W_s^2(\kappa_s,e_s) V_f(x,y)}.
    %     \end{aligned}
    %     \label{eqn: Vf at slow transmission}
    % \end{equation}
    \label{Assumption Vf at slow transmission}%
\end{assum}
Finally, we introduce the next assumption, which is required to guarantee the exponential decay of the composite Lyapunov function $U$ defined in \eqref{eqn: U} during flow, and naturally holds for linear-time-invariant (LTI) SPNCSs as we will see in section \ref{Section LMI}.
\begin{assum}
    Let $\psi_s$ and $\psi_f$ come from Assumption \ref{Assumption interconnection}. There exist $a_{\psi_s}$, $a_{\psi_f}>0$ such that
    $\psi_s(|(x, e_s)|) \leq a_{\psi_s} \sqrt{U_s(\xi_s)}$ and $\psi_f(|(y, e_f)|) \leq a_{\psi_f} \sqrt{U_f(\xi_s,\xi_f)}$.
%    
    % \begin{align*}
    %     \psi_s(|(x, e_s)|) &\leq a_{\psi_s} \sqrt{U_s(\xi_s)} \\
    %     \psi_f(|(y, e_f)|) &\leq a_{\psi_f} \sqrt{U_f(\xi_s,\xi_f)}.
    % \end{align*}
    \label{Assumption Extra 2}%
\end{assum}
We note that $\xi_s $ and $\xi_f$ are defined in \eqref{eqn: definition of xi_s and xi_f}, and contain $(x,e_s)$ and $(y,e_f)$, respectively. Assumption \ref{Assumption Extra 2} naturally holds in LTI systems with UGES protocols \cite{dragan_stability}.

%Before stating our main result, we introduce the attractor $ \mathcal{E} \coloneqq \{ \xi \in \mathbb{X}\ | \ x=0 \wedge e_s = 0 \wedge z=0 \wedge e_f = 0 \} .$
%
% \begin{definition} \red{Semiglobal practical?}
%     We say that set $\mathcal{E}$ is uniformly globally pre-asymptotically stable (UGpAS) for system $\mathcal{H}$ if there exists a class $\mathcal{KL}$ function $\beta$ such that any solution $\xi$ to $\mathcal{H}$ satisfies $|\xi(t,j)|_{\mathcal{E}} \leq \beta(|\xi(0,0)|_{\mathcal{E}}, t+j)$ for all $(t,j)\in dom(\xi)$.
% \end{definition}
% %
% \begin{definition}
%     Suppose the set $\mathcal{E}$ is UGpAS for system $\mathcal{H}$. If every maximal solution to $\mathcal{H}$ is complete, then we say $\mathcal{E}$ is uniformly globally asymptotically stable (UGAS) for system $\mathcal{H}$.
% \end{definition}
%
%\subsection{Stability guarantees}
%------------------------ Theorem ----------------------------
By introducing the attractor $ \mathcal{E} \coloneqq \{ \xi \in \mathbb{X}\ | \ x=0 \wedge e_s = 0 \wedge z=0 \wedge e_f = 0 \}$, we can now state our main results, whose proofs are postponed to the appendix. 
\begin{thm}
Consider system $\mathcal{H}_1$ in \eqref{eqn:full system} and suppose Assumptions \ref{Assumption reduced model}, \ref{Assumption boundary layer system}, \ref{Assumption Vf at slow transmission} and \ref{Assumption Extra 2} hold, and Assumption \ref{Assumption interconnection} holds with $\widetilde {\mathcal{C}} = \mathcal{C}_2^{y,\epsilon}$. 
Let $L_s$ and $\gamma_s$ come from Assumption \ref{Assumption reduced model}, $L_f$ and $\gamma_f$ come from Assumption \ref{Assumption boundary layer system}, and $\lambda_s^*$ and $\lambda_f^*$ come from Lemma \ref{Lemma MATI}.
Then for any $\tau_{\text{miati}}^s \leq \tau_{\text{mati}}^s \leq T(L_s, \gamma_s, \lambda_s^*)$ and $2\tau_{\text{miati}}^f \leq \tau_{\text{mati}}^f \leq \epsilon T(L_f, \gamma_f,\lambda_f^*)$, the following statement holds:

There exists a $\mathcal{KL}$-function $\beta$, such that for all $\Delta, \nu > 0$, there exists an $\epsilon^* >0 $ such that for all $0<\epsilon<\epsilon^*$, any solution $\xi$ with $ |\xi(0,0)|_{\mathcal{E}}<\Delta$ satisfies $|\xi(t,j)|_\mathcal{E} \leq \beta(|\xi(0,0)|_\mathcal{E}, t+j) + \nu$ for any $(t,j)\in \text{dom} \, \xi$.
\label{Theorem H_1}
\end{thm}
%\textbf{Proof:} The proof of Theorem \ref{Theorem H_1} is given in the Appendix.

%

% \begin{remark}
%     \cyan{We note our result is applicable to the double channel SPNCS, with plant and controller in the form of this paper or \cite{SPNCS}. }
%         \todo[inline]{I may explain this remark with more detail, and put this in the thesis not this paper.}
% \end{remark}
Theorem \ref{Theorem H_1} establishes that for any bounded initial condition and ultimate bound, if the condition in Theorem \ref{Theorem H_1} is satisfied, and $\tau_{\text{mati}}^s$, $\tau_{\text{mati}}^f$ and $\epsilon$ are sufficiently small, then the trajectory of system \eqref{eqn:full system} asymptotically approach the ultimate bound.
%
%system \eqref{eqn:full system} satisfies a semiglobal practical asymptotic stability when slow and fast variables are transmitted via a single channel according to the clock mechanism \eqref{eqn: Stefan timer}, with sufficiently small $\tau_{\text{mati}}^s$, $\tau_{\text{mati}}^f$ and $\epsilon$.
% 
%
In the proof of Theorem \ref{Theorem H_1}, it is observed that $\epsilon^*$ approaches zero when $\tau_{\text{miati}}^s$ decreases. This is because the Lyapunov function needs to decrease during flow for some time to compensate the potential increase at slow transmissions. From the perspective of the fast time scale, the transmission interval of slow signals is lower bounded by $\tau_{\text{miati}}^s / \epsilon$. Thus, a smaller $\epsilon$ is required for smaller $\tau_{\text{miati}}^s$ to ensure sufficient flow between any two consecutive slow transmissions.
\subsection{Uniform global asymptotic/exponential stability}\label{Section UGES and UGAS}

Next we are going to present global results such as UGAS and UGES. 
In order to obtain global stability, we need global assumptions. Therefore, we state a global version of Assumption \ref{Assumption interconnection}, which is Assumption \ref{Assumption interconnection Exponential} below.
\begin{assum}
    There exist  $b_1$, $b_2$, $b_3 \geq 0$, such that
\begin{subequations}
    \begin{align}
        &\begin{aligned}
            &\Big < \tfrac{\partial {U_s}}{\partial \xi_s}, F_s^y(x,y,e_s,e_f) - F_s^y(x,0,e_s,0)  \Big> \leq
            \\
            &\phantom{aaaaaaaaaaaa} b_1 \psi_s\left(\left| (x, e_s) \right|\right) \psi_f\left(\left| (y, e_f) \right|\right), 
        \end{aligned}\label{eqn: SPNCS interconnection Exponential 1}
        \\
        &\begin{aligned}
            &\Big< \tfrac{\partial {U_f}}{\partial \xi_s} - \tfrac{\partial {U_f}}{\partial y} \tfrac{\partial \overline{H}}{\partial \xi_s} - \tfrac{\partial {U_f}}{\partial e_f} \tfrac{\partial \tilde k}{\partial \xi_s} ,  F_s^y(x,y,e_s,e_f) \Big> \leq
            \\
            & b_2 \psi_s\left(\left| (x, e_s) \right|\right) \psi_f\left(\left| (y, e_f) \right|\right) + b_3 \psi_f^2\left(\left| (y, e_f) \right|\right)
        \end{aligned}\label{eqn: SPNCS interconnection Exponential 2}
    \end{align}
    \label{eqn: SPNCS interconnections Exponential}%
\end{subequations}
hold for almost all $\xi^y \in \mathcal{C}_2^{y,\epsilon}$, where $\tilde k(x,z) = (k_{pf}(x_p,z_p), k_{cf}(x_c,z_c))$.
    \label{Assumption interconnection Exponential}
\end{assum}

% The following assumption guarantees exponential decay of the composite Lyapunov function $U$, see the proof of Corollary \ref{Corollary UGAS} for details.

% \begin{assum}
%     There exist $a_{\psi_s}$, $a_{\psi_f}>0$ such that
%     \begin{align*}
%         \psi_s(|(x, e_s)|) &\leq a_{\psi_s} \sqrt{U_s(\xi_s)} \\
%         \psi_f(|(y, e_f)|) &\leq a_{\psi_f} \sqrt{U_f(\xi_s,\xi_f)}.
%     \end{align*}
%     \label{Assumption Extra 2}
% \end{assum}

\begin{cor}
Considering system $\mathcal{H}_1$ in \eqref{eqn:full system} and suppose Assumptions \ref{Assumption reduced model}, \ref{Assumption boundary layer system}, \ref{Assumption Vf at slow transmission}, \ref{Assumption Extra 2} and \ref{Assumption interconnection Exponential} hold.
Let $L_s$ and $\gamma_s$ come from Assumption \ref{Assumption reduced model}, and $L_f$ and $\gamma_f$ come from Assumption \ref{Assumption boundary layer system}.
Let $b_1$, $b_2$, and $b_3$ come from Assumption \ref{Assumption interconnection Exponential} and $a_s$, $a_f$, $\lambda_s^*$ and $\lambda_f^*$ come from Lemma \ref{Lemma MATI}.
Then for any $\tau_{\text{miati}}^s \leq \tau_{\text{mati}}^s \leq T(L_s, \gamma_s, \lambda_s^*)$ and $ 2\tau_{\text{miati}}^f \leq \tau_{\text{mati}}^f \leq \epsilon T(L_f, \gamma_f,\lambda_f^*)$, there exists $\epsilon^* >0 $ and $\beta \in \mathcal{KL}$, such that for all $0<\epsilon<\epsilon^*$ and $ |\xi(0,0)|_{\mathcal{E}} \in \mathbb{R}$, we have $|\xi(t,j)|_\mathcal{E} \leq \beta(|\xi(0,0)|_\mathcal{E}, t+j)$ for any $(t,j)\in \text{dom} \, \xi$.
% and  $\mathcal{H}_1$ is UGAS w.r.t $\mathcal{E}$.
\label{Corollary UGAS}
\end{cor}
Corollary \ref{Corollary UGAS} is proved similarly to Theorem \ref{Theorem H_1} by setting $\Delta$ and $\nu$ to infinity and zero, respectively, so its proof is omitted.
We can also state UGES of $\mathcal{H}_1$ if $\mathcal{H}_{bl}$ and $\mathcal{H}_r$ are UGES, and if $h_s$ and $h_f$ are UGES protocols.
The following assumption guarantees these conditions when combined with Assumptions \ref{Assumption reduced model} and \ref{Assumption boundary layer system}.
\begin{assum}
    Let $W_s, W_f, V_s, V_f, \rho_s$ and $\rho_f$ come from Assumptions \ref{Assumption reduced model} and \ref{Assumption boundary layer system}. There exist positive real numbers $\underline{a}_{W_s}$, $\overline{a}_{W_s}$, $\underline{a}_{V_s}$, $\overline{a}_{V_s}$, $\underline{a}_{W_f}$, $\overline{a}_{W_f}$, $\underline{a}_{V_f}$, $\overline{a}_{V_f}$, $a_{\rho_s}$ and $a_{\rho_f}$ such that 
    \begin{subequations}
        \begin{align}
            \underline{a}_{W_s} |e_s| \leq & W_s(\kappa_s, e_s) \leq \overline{a}_{W_s} |e_s|, \label{eqn: Ws exponential sandwich bound} \\
            \underline{a}_{W_f} |e_f| \leq & W_f(\kappa_f, e_f) \leq \overline{a}_{W_f} |e_f|, \label{eqn: Wf exponential sandwich bound} \\
            \underline{a}_{V_s} |x|^2 \leq &V_s(x) \leq \overline{a}_{V_s} |x|^2, \label{eqn: Vs exponential sandwich bound}\\
            \underline{a}_{V_f} |y|^2 \leq &V_f(x,y) \leq \overline{a}_{V_f} |y|^2, \label{eqn: Vf exponential sandwich bound}%
        \end{align}
    \end{subequations}
 and \vspace{-0.5cm}
\begin{subequations}
\begin{align}
    a_{\rho_s}s^2 &\leq \rho_s(s), \label{eqn: a_{rho_s}}
    \\
    a_{\rho_f}s^2 &\leq \rho_f(s)\label{eqn: a_{rho_f}}
    \end{align}
\end{subequations}
for all $s\in \mathbb{R}$. Additionally, $\overline{H}$ is globally Lipschitz in both arguments.
\label{Assumption Exponential}
\end{assum}
We note that Assumption \ref{Assumption Exponential} implies Assumption \ref{Assumption Extra 2} holds, this can be see in the proof of Theorem \ref{Theorem Exponential decay}. Additionally, \eqref{eqn: Ws exponential sandwich bound} and \eqref{eqn: NCS assumption Ws jump} imply $W_s$ is an UGES protocol, and the similar statement holds for $W_f$.

\begin{thm}
Consider system $\mathcal{H}_1$ in \eqref{eqn:full system} and suppose Assumptions \ref{Assumption reduced model}, \ref{Assumption boundary layer system}, \ref{Assumption Vf at slow transmission}, \ref{Assumption interconnection Exponential} and \ref{Assumption Exponential} hold.
%Let $b_1$, $b_2$, and $b_3$ come from Assumption \ref{Assumption interconnection} and $a_s$ and $a_f$ come from Lemma \ref{Lemma MATI}. 
%%%%%%%%%%%%%%%%%%%%%%
Let $L_s$ and $\gamma_s$ come from Assumption \ref{Assumption reduced model}, and $L_f$ and $\gamma_f$ come from Assumption \ref{Assumption boundary layer system}.
Let 
%$b_1$, $b_2$, and $b_3$ come from Assumption \ref{Assumption interconnection Exponential} and $a_s$, $a_f$, 
$\lambda_s^*$ and $\lambda_f^*$ come from Lemma \ref{Lemma MATI}.
Then for any $\tau_{\text{miati}}^s \leq \tau_{\text{mati}}^s \leq T(L_s, \gamma_s, \lambda_s^*)$ and $ 2\tau_{\text{miati}}^f \leq \tau_{\text{mati}}^f \leq \epsilon T(L_f, \gamma_f,\lambda_f^*)$, there exist $\epsilon^*, c_1, c_2 >0 $ such that for all $0<\epsilon<\epsilon^*$, the solution $\xi$ satisfies $|\xi(t,j)|_\mathcal{E} \leq c_1 |\xi(0,0)|_\mathcal{E} \text{exp}( -c_2(t+j))$ for any $(t,j)\in \text{dom} \, \xi$.
%$\mathcal{H}_1$ is UGES w.r.t $\mathcal{E}$.
\label{Theorem Exponential decay}
\end{thm}
%\textbf{Proof:} The proof of Theorem \ref{Theorem Exponential decay} is given in the Appendix.
%
%
\begin{rem}
    We state only global assumptions and results (i.e., UGAS and UGES) to simplify the presentation. However, for local results, all assumptions can be rephrased. For example, the requirement for a unique real solution in SA\ref{assum:standing-ss} can be relaxed to isolated real solutions, and the interconnection condition in Assumption \ref{Assumption interconnection Exponential} needs to hold only on a compact set of $\xi$ containing the set $\mathcal{E}$. 
\end{rem}

\section{Further stability results}
%\input{Chapters/7_Stable_subsystem}
% This section will study two specific configurations of our single-channel SPNCSs, both of them are special cases of $\mathcal{H}_1$, with lower dimension.

% There are two typical configurations usually adopted in the literature: assuming either the slow or fast subsystems are stable and thus do not require stabilization via the network. In this paper, we do not rely on these assumptions, though they are special cases of our main results. 

In this section, we consider the scenario when the plant $\mathcal{P}$ has a stable fast subsystem (i.e., boundary layer system) and the controller is designed to stabilise the slow subsystem. 
%This scenario typically arises when we have sensors or actuators operating on a faster time scale. 
%
Due to the stable fast subsystem, the plant only needs to transmit  $y_s$ and the controller only generates $u_s$. The plant and controller below are in the form of \eqref{eqn:plant} and \eqref{eqn:controller}, where $k_{pf}$ and $k_{cf}$ are removed, that is, $y_p = y_s = k_{p_s}(x_p) $ and $u = u_s = k_{c_s}(x_c)$.
% \begin{equation*}
%     \mathcal{P}\!:\!
%     \begin{cases}
%     \dot x_p \!\!\!\!\!\!&= f_p(x_p, z_p,\hat u)\\
%     \epsilon \dot z_p\!\!\!\!\!\! &= g_p(x_p, z_p, \hat u) \\
%     y_p \!\!\!\!\!\!&= y_s = k_{p_s}(x_p)  , 
%     \end{cases} \qquad
%     \mathcal{C}\!:\!
%     \begin{cases}
%     \dot x_c \!\!\!\!\!\!&= f_c(x_c, z_c, \hat{y}_p)\\
%     \epsilon \dot z_c \!\!\!\!\!\!&= g_c(x_c, z_c, \hat y_p) \\
%     u \!\!\!\!\!\!&= u_s = k_{c_s}(x_c).
%     \end{cases}
% \end{equation*}
%
Given the exclusive presence of slow transmissions within the communication channel, a simpler version of clock mechanism \eqref{eqn: Stefan timer} is used to govern the system dynamics, that is $\tau_{\text{miati}}^s \leq t_{k+1}^s - t_k^s \leq \tau_{\text{mati}}^s$ for all $t_k^s, t_{k+1}^s\in \mathcal{T}^s$ and $k \in \mathbb{Z}_{\geq 1}$.
% \begin{equation*}
%     \tau_{\text{miati}}^s \leq t_{k+1}^s - t_k^s \leq \tau_{\text{mati}}^s, \ \forall t_k^s, t_{k+1}^s\in \mathcal{T}^s, k \in \mathbb{Z}_{\geq 1}.
% \end{equation*}
Moreover, the networked induced error becomes $e_s =  (\hat{y}_s - y_s, \hat{u}_s - u_s)$, and $e_f$ no longer exist. We define the state $\xi_2 \coloneqq (x,e_s, \tau_s, \kappa_s, z)\in \mathbb{X}_2$, with $\mathbb{X}_2\coloneqq \mathbb{R}^{n_x}\times \mathbb{R}^{n_{e_s}}\times  \mathbb{R}_{\geq 0} \times \mathbb{Z}_{\geq 0} \times \mathbb{R}^{n_z}$, and we define our system by the hybrid model $\mathcal{H}_2$, and we omit its expression here.

% \begin{equation*}
%     \mathcal{H}_2:\left\{
% \begin{aligned}
%     \dot{\xi}_2 &= F_2(\xi_2),\ \xi_2 \in \mathcal{C}_2, \\
%     \xi^+ &= G_{2,s}(\xi_2), \ \xi_2 \in \mathcal{D}_{2,s},
% \end{aligned}
%     \right.
% \end{equation*}
% where $F_2(\xi_2) \coloneqq \big(f_x(x,z,e_s,e_f), f_{e_s}(x,z,e_s,e_f),1,0,$ $\tfrac{1}{\epsilon} g_z(x,z,e_s,e_f) \big)$, $C_2 \coloneqq \{ \xi_2 \in \mathbb{X}_2 | \tau_s \in [0, \tau_{\text{mati}}^s]\}$ , $G_{2,s}(\xi_2) \coloneqq (x, h_s(\kappa_s, e_s), 0, \kappa_s + 1, z)$ and $D_{2,s} \coloneqq \{ \xi_2 \in \mathbb{X}_2 | \tau_s \in [\tau_{\text{miati}}^s, \tau_{\text{mati}}^s]\}$. We emphasize again that $k_{pf} \equiv 0$, $k_{cf} \equiv 0$ and $e_f \equiv 0$ in $\mathcal{H}_2$.

We again define $y$ as in \eqref{eqn: map between y and z}.
%$y \coloneqq z - \overline{H}(x,e_s)$, where $\overline{H}$ comes from SA\ref{assum:standing-ss}. 
Let $\xi_2^y \coloneqq (x,e_s,\tau_s, \kappa_s,y) \in \mathbb{X}_2$,
then we can define $\mathcal{H}_2^y$, which is $\mathcal{H}_2$ after changing the coordinates. Since $\mathcal{H}_2$ is a special case of $\mathcal{H}_1$, the functions we used to define $\mathcal{H}_2^y$ have the same form as $\mathcal{H}_1^y$ (e.g., \eqref{eq:functions}), but no longer depend on the state $e_f$, nor functions $k_{pf}$ and $k_{cf}$. For example, we will have that $\dot x$ equals to $f_x\big(x,y+\overline{H}(x,e_s), e_s\big)$, not $f_x\big(x,y+\overline{H}(x,e_s), e_s, e_f\big)$. Then we can write
\begin{equation*}
    \mathcal{H}_2^y:\left\{
\begin{aligned}
    \dot{\xi}_2^y &= F_2^y(\xi_2^y, \epsilon),\ \xi_2^y \in \mathcal{C}_2^y, \\
    {\xi_2^y}^+ &= G_{2,s}(\xi_2^y), \ \xi_2^y \in \mathcal{D}_{2,s}^y,
\end{aligned}
    \right.
\end{equation*}
where $F_2^y(\xi_2^y, \epsilon) \coloneqq \big(F_s^y(x,y,e_s), \tfrac{1}{\epsilon} (\epsilon \tfrac{\partial y}{\partial t})  \big)$, with $F_s^y$ and $\epsilon \tfrac{\partial y}{\partial t}$ come from system \eqref{eqn: H_2^y}, $C_2^y \coloneqq \{ \xi_2^y \in \mathbb{X}_2 | \tau_s \in [0, \tau_{\text{mati}}^s]\}$, $G_{2,s}^y(\xi_2^y) \coloneqq \big(x, h_s(\kappa_s, e_s), 0, \kappa_s + 1, h_y(x,e_s,y)\big)$ and $D_{2,s}^y \coloneqq \{ \xi_2^y \in \mathbb{X}_2^y | \tau_s \in [\tau_{\text{miati}}^s, \tau_{\text{mati}}^s]\}$.
Furthermore, we can derive a continuous time boundary-layer system $\mathcal{H}_{2,bl} : \{
        \tfrac{\partial y}{\partial \sigma} = g_z(x,y+\overline{H}(x,e_s),e_s)
        $,
% \begin{equation*}
%     \mathcal{H}_{2,bl} : \begin{cases}
%         \tfrac{\partial y}{\partial \sigma} = g_z(x,y+\overline{H}(x,e_s),e_s),
%     \end{cases}
% \end{equation*}
and a NCS $\mathcal{H}_{2,r}$ 
\begin{equation*}
    \mathcal{H}_{2,r}:\left\{
\begin{aligned}
    \dot{\xi_s} &= F_s^y(\xi_s),\ \xi_2^y \in \mathcal{C}_2^y, \\
    \xi_s^+ &= \big(x,h_s(\kappa_s, e_s), 0, \kappa_s + 1\big), \ \xi_2^y\in \mathcal{D}_{2,s}^y,
\end{aligned}
    \right.    
\end{equation*}
where we recall that $\xi_s \coloneqq (x,e_s,\tau_s, \kappa_s)$.
We note that only $\mathcal{H}_{2,r}$ is a hybrid system, while $\mathcal{H}_{2,bl}$ is not, as it is already stable and does not need to be stabilized through network transmissions. Consequently, we need to modify relevant assumptions. In this case, we adjust Assumption \ref{Assumption boundary layer system}, resulting in Assumption \ref{Assumption Stable fast subsystem} below. 

\begin{assum}
    Consider $\mathcal{H}_2^y$, there exist locally Lipschitz function $V_f: \mathbb{R}^{n_x} \times \mathbb{R}^{n_z}\rightarrow \mathbb{R}_{\geq 0}$, class $\mathcal{K}_\infty$ functions $\underline{\alpha}_{V_f}$ and $\overline{\alpha}_{V_f}$, $a_f>0$ and positive definite function $\psi_f$ such that for all $\xi_2 \in \mathbb{X}_2$, we have 
    $\underline{\alpha}_{V_f}\left(\left| y \right|\right)\leq {V_f}(x,y) \leq \overline{\alpha}_{V_f}\left(\left| y\right|\right)$ and $\big< \tfrac{\partial {V_f}(x,y)}{\partial y},g_z(x,y+ \overline{H}(x, e_s),e_s)  \big> \leq -a_f \psi_f^2 \left(| y |\right)$.  
    % \begin{align}
    % \underline{\alpha}_{V_f}\left(\left| y \right|\right)\leq {V_f}(x,y) & \leq \overline{\alpha}_{V_f}\left(\left| y\right|\right) \\
    %  \left< \tfrac{\partial {V_f}(x,y)}{\partial y},g_z(x,y+ \overline{H}(x, e_s),e_s)  \right> &\leq -a_f \psi_f^2 \left(| y |\right).
    %  \label{eqn: Assumption stable fast subsystem, Uf flow}
    % \end{align}
    \label{Assumption Stable fast subsystem}
\end{assum}
Assumptions \ref{Assumption reduced model}, \ref{Assumption interconnection} and \ref{Assumption Extra 2} are written for a more general model (i.e., $\mathcal{H}_1^y$), when we apply them to a specialized model with a lower dimension, we ignore the states that do not exist in the specialized model. %\cyan{For instance, $\mathcal{H}_2$ does not have the state $e_f$, then when we apply Assumption \ref{Assumption reduced model} to it, $f_x(x, \overline{H}(x,e_s),e_s,0)$ is replaced by $f_x(x, \overline{H}(x,e_s),e_s)$.} 
At the same time, Lemma \ref{Lemma MATI} is applicable to any reduced and boundary layer system in the form of NCS.
%with the same form as $\mathcal{H}_r$ and $\mathcal{H}_{bl}$. 
Since $\mathcal{H}_{2,r}$ is a NCS and Assumption \ref{Assumption reduced model} holds, we can conclude inequalities in \eqref{eqn: Us} hold, with $\mathcal{C}_{2,r}^{y,0}$ and $\mathcal{D}_{2,r}^{y,0}$ replaced by $\mathcal{C}_2^y$ and $\mathcal{D}_{2,s}^y$, respectively. Since there is no fast transmission, we only have $V_f$ but not $U_f$, and all $U_f$ in Assumption \ref{Assumption interconnection} and \ref{Assumption Extra 2} should be replace by $V_f$.
We define the set $\mathcal{E}_2 \coloneqq \{\xi_2 \in \mathbb{X}_2 | x=0 \wedge e_s = 0 \wedge z=0 \}$.
\begin{cor}
    Consider system $\mathcal{H}_2$ and suppose Assumptions \ref{Assumption reduced model}, \ref{Assumption Vf at slow transmission}, \ref{Assumption Extra 2} and \ref{Assumption Stable fast subsystem} hold, and Assumption \ref{Assumption interconnection} holds with $\widetilde{\mathcal{C}} =  \mathcal{C}_2^y$.
    Let $b_1$, $b_2$, and $b_3$ come from Assumption \ref{Assumption interconnection}, $a_s$ comes from Lemma \ref{Lemma MATI} and $a_f$ comes from Assumption \ref{Assumption Stable fast subsystem}. Then for any $\tau_{\text{miati}}^s \leq \tau_{\text{mati}}^s \leq T(L_s, \gamma_s, \lambda_s^*)$, the following statement holds:

    There exists a $\mathcal{KL}$-function $\beta$, such that for all $\Delta, \nu>0$, there exist $\epsilon^* >0$ such that for all $0<\epsilon<\epsilon^*$, any solution $\xi_2$ with $ |\xi_2(0,0)|_{\mathcal{E}_2}<\Delta$ satisfies $|\xi_2(t,j)|_{\mathcal{E}_2} \leq \beta(|\xi_2(0,0)|_{\mathcal{E}_2}, t+j) + \nu$ for any $(t,j)\in \text{dom} \, \xi_2$.
    \label{Corollary Stable fast subsystem}
\end{cor}
Corollary \ref{Corollary Stable fast subsystem} is proved similarly to Theorem \ref{Theorem H_1} by defining the composite Lyapunov function as $U_2(\xi_2^y)\coloneqq {U_s}(\xi_s) + d{V_f}(x,y)$. Its proof is therefore omitted.
\begin{rem}
    The stability analysis for system with a stable slow subsystem can be conducted similarly.
\end{rem}

\section{Special case: linear time invariant systems} \label{Section LMI}
We show in this section how to apply the result seen so far to a LTI plant and a LTI controller with RR or TOD protocols. Consider systems \eqref{eqn:plant} and \eqref{eqn:controller} as
%
% In this section, we consider an LTI plant and controller, with UGES protocols, then we illustrate conditions in Theorem \ref{Theorem Exponential decay} holds if we satisfy two linear matrix inequalities (LMIs). 
% Let the plant and controller be defined as
\begin{equation}
\begin{aligned}
    \left[ \begin{smallmatrix}
        \dot{x}_p \\ \epsilon \dot{z}_p
    \end{smallmatrix} \right]
    &=
    \left[ \begin{smallmatrix}
        A_{11}^p & A_{12}^p \\ A_{21}^p & A_{22}^p
    \end{smallmatrix} \right]
    \left[ \begin{smallmatrix}
        x_p \\ z_p
    \end{smallmatrix} \right] 
    + 
    \left[ \begin{smallmatrix}
        A_{13}^p & A_{14}^p \\ A_{23}^p & A_{24}^p
    \end{smallmatrix} \right]
    \left[ \begin{smallmatrix}
        \hat{u}_s \\ \hat{u}_f
    \end{smallmatrix} \right],
    \\
    \left[ \begin{smallmatrix}
        y_s \\ y_f
    \end{smallmatrix} \right]
    &=
    \left[ \begin{smallmatrix}
        A_x^{p_s} & 0 \\ A_x^{p_f} & A_z^{p_f}
    \end{smallmatrix} \right]
    \left[ \begin{smallmatrix}
        x_p \\ z_p
    \end{smallmatrix} \right],
    \\ 
    \left[ \begin{smallmatrix}
        \dot{x}_c \\ \epsilon \dot{z}_c
    \end{smallmatrix} \right]
    &=
    \left[ \begin{smallmatrix}
        A_{11}^c & A_{12}^c \\ A_{21}^c & A_{22}^c
    \end{smallmatrix} \right]
    \left[ \begin{smallmatrix}
        x_c \\ z_c
    \end{smallmatrix} \right] 
    + 
    \left[ \begin{smallmatrix}
        A_{13}^c & A_{14}^c \\ A_{23}^c & A_{24}^c
    \end{smallmatrix} \right]
    \left[ \begin{smallmatrix}
        \hat{y}_s \\ \hat{y}_f
    \end{smallmatrix} \right],
    \\
    \left[ \begin{smallmatrix}
        u_s \\ u_f
    \end{smallmatrix} \right]
    &=
    \left[ \begin{smallmatrix}
        A_x^{c_s} & 0 \\ A_x^{c_f} & A_z^{c_f}
    \end{smallmatrix} \right]
    \left[ \begin{smallmatrix}
        x_c \\ z_c
    \end{smallmatrix} \right].
\end{aligned}
\label{eqn: linear plant and controller}
\end{equation}
The hybrid model that describes our SPNCS is given by \eqref{eqn:full system}, with $F(\xi, \epsilon) =  \big(f_x,f_{e_s},1,0,\tfrac{1}{\epsilon}g_z, \tfrac{1}{\epsilon} g_{e_f},  \frac{1}{\epsilon},0\big)$, where
\begin{equation*}
    \left[\begin{smallmatrix}
        f_x \\ f_{e_s} \\ g_z \\ g_{e_f}
    \end{smallmatrix} \right]
    =
    \left[\begin{smallmatrix}
        A_{11} & A_{12} & A_{13} & A_{14} \\
        A_{21} & A_{22} & A_{23} & A_{24} \\
        A_{31} & A_{32} & A_{33} & A_{34} \\
        \epsilon A_{41}^\epsilon + A_{41} & \epsilon A_{42}^\epsilon + A_{42} & \epsilon A_{43}^\epsilon + A_{43} & \epsilon A_{44}^\epsilon + A_{44} \\
    \end{smallmatrix}\right]
    \left[\begin{smallmatrix}
        x \\ e_s \\  z \\  e_f
    \end{smallmatrix}\right],
\end{equation*}
$A_{11} = \left[\begin{smallmatrix}A_{11}^p & A_{13}^p A_x^{c_s} + A_{14}^p A_x^{c_f} \\ A_{13}^c A_x^{p_s} + A_{14}^c A_x^{p_f} & A_{11}^c \end{smallmatrix} \right]$,
$A_{12} = \left[\begin{smallmatrix} 0 & A_{13}^p \\ A_{13}^c & 0\end{smallmatrix}\right]$,
$A_{13} = \left[\begin{smallmatrix} A_{12}^p & A_{14}^p A_{z}^{c_f} \\ A_{14}^c A_z^{p_f} & A_{12}^c \end{smallmatrix} \right]$,
$A_{14} = \left[ \begin{smallmatrix}0 & A_{14}^p \\ A_{14}^c & 0\end{smallmatrix} \right]$,
$A_{21} = A_x^s A_{11}$,
$A_{22} = A_x^s A_{12}$,
$A_{23} = A_x^s A_{13}$,
$A_{24} = A_x^s A_{14}$,
$A_{31} = \left[\begin{smallmatrix}
    A_{21}^p & A_{23}^p A_x^{c_s} + A_{24}^p A_x^{c_f} \\
    A_{23}^c A_x^{p_s} + A_{24}^c A_x^{p_f} & A_{21}^c
\end{smallmatrix}\right]$,
$A_{32} = \left[\begin{smallmatrix}
    0 & A_{23}^p \\ A_{23}^c & 0
\end{smallmatrix} \right]$,
$A_{33} = \left[\begin{smallmatrix}
    A_{22}^p & A_{24}^p A_z^{c_f} \\ A_{24}^c A_z^{p_f} & A_{22}^c
\end{smallmatrix} \right]$,
$A_{34} = \left[\begin{smallmatrix}
    0 & A_{24}^p \\ A_{24}^c & 0
\end{smallmatrix} \right]$,
$A_{41}^\epsilon = A_x^f A_{11}$, 
$A_{42}^\epsilon = A_x^f A_{12}$,
$A_{43}^\epsilon = A_x^f A_{13}$, 
$A_{44}^\epsilon = A_x^f A_{14}$, 
$A_{41} = A_z^f A_{31}$,
$A_{42} = A_z^f A_{32}$,
$A_{43} = A_z^f A_{33}$,
$A_{44} = A_z^f A_{34}$,
$A_x^s = \left[\begin{smallmatrix}
    -A_x^{p_s} & 0 \\ 0 & -A_x^{c_s}
\end{smallmatrix} \right]$,
$A_x^f = \left[\begin{smallmatrix}
    -A_x^{p_f} & 0 \\ 0 & -A_x^{c_f}
\end{smallmatrix} \right]$ and 
$A_z^f = \left[\begin{smallmatrix}
    -A_z^{p_f} & 0 \\ 0 & -A_z^{c_f}
\end{smallmatrix} \right]$.

The quasi-steady state of $z$, which is denoted by $\overline{H}(x,e_s)$, is given by
\begin{equation}
    \overline{H}(x,e_s) = - A_{33}^{-1} A_{31} x - A_{33}^{-1} A_{32} e_s.
    \label{eqn: H bar linear}
\end{equation}
Recall that $y$ is defined in \eqref{eqn: map between y and z}, then by setting $\epsilon$ to zero, the boundary-layer system $\mathcal{H}_{bl}$ is given by \eqref{eqn: H_bl}, where $F_f^y(x,y,e_s,e_f,0)$ is specified in $\eqref{eqn: linear functions}$. The reduced system $\mathcal{H}_{r}$ is given by \eqref{eqn: H_r}, where $F_s^y(x,0,e_s, 0)$ is given in \eqref{eqn: linear functions}.
\begin{equation}
    \begin{aligned}
        &F_f^y(x,y,e_s,e_f,0) = (A_{11}^f y + A_{12}^f e_f, A_{21}^f y + A_{22}^f e_f, 1, 0), \\
        &F_s^y(x,0,e_s, 0) = (A_{11}^s x + A_{12}^s e_s, A_{21}^s x + A_{22}^s e_s, 1, 0), \\
        &A_{11}^f = A_{33}, \ A_{12}^f = A_{34},\ A_{21}^f = A_z^f A_{33},\ A_{22}^f = A_z^f A_{34}, \\
        &A_{11}^s = A_{11} - A_{13}A_{33}^{-1}A_{31},\ A_{12}^s = A_{12} - A_{13}A_{33}^{-1}A_{32}, \\
        &A_{21}^s = A_{21} - A_{23}A_{33}^{-1}A_{31},\ A_{22}^s = A_{22} - A_{23}A_{33}^{-1}A_{32}.   
    \end{aligned}
    \label{eqn: linear functions}
\end{equation}
%
% By Propositions 4 and 5 in \cite{dragan_stability}, there exist positive definite function $W_s$, positive constants $\underline{a}_{W_s}, \overline{a}_{W_s}$ and $\lambda_s \in [0, 1)$ such that \eqref{eqn: NCS assumption Ws sandwich bound}, \eqref{eqn: NCS assumption Ws jump} and \eqref{eqn: Ws exponential sandwich bound} hold.
% %
% Moreover, Examples 3 and 4 in \cite{dragan_stability} show that there exist $L_s \geq 0$ and a matrix $A_{H_s}$, such that \eqref{eqn: NCS Ws dot} holds with $H_s(x,e_s) = \left| A_{H_s} x \right|$.
% %
% Similarly, we can show there exist a $W_f$ locally lipschitz function $W_f$, positive constants $\underline{a}_{W_f}, \overline{a}_{W_f}$, $\lambda_f \in [0, 1)$, $L_f \geq 0$ and a matrix $A_{H_f}$, such that \eqref{eqn: NCS assumption Wf sandwich bound}-\eqref{eqn: NCS Wf dot} and \eqref{eqn: Ws exponential sandwich bound} are satisfied, with $H_f(y,e_f) = \left| A_{H_f} y \right|$.

\begin{claim}
    For LTI plant and controller given by \eqref{eqn: linear plant and controller}, with RR or TOD protocols, there exist positive definite function $W_s$, positive constants $\underline{a}_{W_s}, \overline{a}_{W_s}$ and $\lambda_s \in [0, 1)$ such that \eqref{eqn: NCS assumption Ws sandwich bound}, \eqref{eqn: NCS assumption Ws jump} and \eqref{eqn: Ws exponential sandwich bound} hold. there exist $L_s \geq 0$ and a matrix $A_{H_s}$, such that \eqref{eqn: NCS Ws dot} holds with $H_s(x,e_s) = \left| A_{H_s} x \right|$. Similarly, there exist a locally lipschitz function $W_f$, positive constants $\underline{a}_{W_f}, \overline{a}_{W_f}$, $\lambda_f \in [0, 1)$, $L_f \geq 0$ and a matrix $A_{H_f}$, such that \eqref{eqn: NCS assumption Wf sandwich bound}-\eqref{eqn: NCS Wf dot} and \eqref{eqn: Ws exponential sandwich bound} are satisfied, with $H_f(y,e_f) = \left| A_{H_f} y \right|$.
    \label{Claim for LTI section}
\end{claim}
\textbf{Proof:} Claim \ref{Claim for LTI section} is obtained by inspecting Propositions 4 and 5, as well as Examples 3 and 4 in \cite{dragan_stability}.

\begin{prop}
    Consider system \eqref{eqn:full system}, with the LTI plant and controller specified in \eqref{eqn: linear plant and controller}, as well as RR or TOD protocols. Let $\underline{a}_{W_s}$, $\underline{a}_{W_f}$, $A_{H_s}$ and $A_{H_f}$ come from Claim \ref{Claim for LTI section}. Suppose there exist $a_{\rho_s}$, $a_{\rho_f}$, $\gamma_s$, $\gamma_f > 0$ and positive definite symmetric real matrices $P_s$ and $P_f$, such that the following LMI holds for $\ell \in \{s, f\}$.
%
    % \begin{subequations}
    %     \begin{align}
    %     \left[\begin{smallmatrix}
    %     A_{11}^s P_s + P_s A_{11}^{s\top} + a_{\rho_s} I + A_{H_s}^\top A_{H_s} &  \bigstar  \\
    %     A_{12}^{s\top} P_s & a_{\rho_s} I - \gamma_s^2 \underline{a}_{W_s}^2 I
    %     \end{smallmatrix}\right]
    %     \leq 0,
    %     \label{eqn: LMI slow}
    %     \\
    %     \left[\begin{smallmatrix}
    %     A_{11}^f P_f + P_f A_{11}^{f\top} + a_{\rho_f} I + A_{H_f}^\top A_{H_f} & \bigstar \\
    %     A_{12}^{f\top} P_f & a_{\rho_f} I - \gamma_f^2 \underline{a}_{W_f}^2 I
    %     \end{smallmatrix}\right] \leq 0.
    %     \label{eqn: LMI fast}
    %     \end{align}
    %     \label{eqn: LMIs}%
    % \end{subequations}
    \begin{equation}
        \left[\begin{smallmatrix}
        A_{11}^{\ell} P_\ell + P_\ell A_{11}^{\ell\top} + a_{\rho_\ell} I + A_{H_\ell}^\top A_{H_\ell} &  \bigstar  \\
        A_{12}^{\ell\top} P_\ell & a_{\rho_\ell} I - \gamma_\ell^2 \underline{a}_{W_\ell}^2 I
        \end{smallmatrix}\right]
        \leq 0 .
        \label{eqn: LMIs}
    \end{equation}
    Then conditions in Theorem \ref{Theorem Exponential decay} hold with $V_s = x^\top P_s x$, $V_f = y^\top P_f y$, $\gamma_s$ and $\gamma_f$ from \eqref{eqn: LMIs}, as well as $\lambda_s^* \in (\lambda_s, 1)$, $\lambda_f^* \in (\lambda_f,1)$ from Lemma \ref{Lemma MATI}, with $\lambda_s$ and $\lambda_f$ come from Claim \ref{Claim for LTI section}.
    %$\underline{a}_{V_s} = \lambda_{\text{min}}(P_s)$ and $\overline{a}_{V_s} = \lambda_{\text{max}}(P_s)$, $\underline{a}_{V_f} = \lambda_{\text{min}}(P_f)$ and $\overline{a}_{V_f} = \lambda_{\text{max}}(P_f)$,.
    \label{Proposition LTI}  
\end{prop}
\textbf{Proof:} The proof of Proposition \ref{Proposition LTI} is given in the Appendix.
\begin{rem}
    Proposition \ref{Proposition LTI} can be easily extended to other UGES protocols as long as \eqref{eqn: NCS Ws dot} and \eqref{eqn: NCS Wf dot} hold with $H_s(x,e_s) = \left| A_{H_s} x \right|$ and $H_f(y,e_f) = \left| A_{H_f} y \right|$. See illustrative example and \cite{nesic2009unified} for more details.
\end{rem}
Proposition \ref{Proposition LTI} implies that for SPNCS with an LTI plant, an LTI controller, and RR or TOD protocols, the satisfaction of the LMI \eqref{eqn: LMIs} guarantees that we can always find sufficiently small $\tau_{\text{mati}}^s$, $\tau_{\text{mati}}^f$ and $\epsilon^*$, such that if $\epsilon < \epsilon^*$, system \eqref{eqn:full system} considered in this section is UGES. Two necessary conditions to guarantee feasibility of \eqref{eqn: LMIs} are $A_{11}^{\ell} P_\ell + P_\ell A_{11}^{\ell\top} + a_{\rho_\ell} I + A_{H_\ell}^\top A_{H_\ell} < 0 $ and $a_{\rho_\ell} I - \gamma_\ell^2 \underline{a}_{W_\ell}^2 I <0$, where the first condition can be satisfied if $A_{11}^{\ell}$ is Hurwitz, and the second condition can always be verified by selecting $a_{\rho_\ell}$ sufficiently small and $\gamma_{\ell}$ sufficiently large.

\section{An illustrative example}
This section provides a numerical example of the result of section \ref{Section LMI}.
%an example to show how to determine stability of the system using Theorem \ref{Theorem Exponential decay} and Section VII. 
%
% Let the plant and controller be defined as
% \begin{equation*}
% \begin{aligned}
%     \left[ \begin{smallmatrix} 
%         \dot{x}_p \\ \epsilon \dot{z}_p
%     \end{smallmatrix} \right]
%     &=
%     \left[ \begin{smallmatrix} 
%         A_{11}^p & A_{12}^p \\
%         A_{21}^p & A_{22}^p
%     \end{smallmatrix} \right]
%     \left[ \begin{smallmatrix} 
%         x_p \\ z_p
%     \end{smallmatrix} \right] 
%     + 
%     \left[ \begin{smallmatrix} 
%         A_{13}^p \\ A_{23}^p
%     \end{smallmatrix} \right]
%         \hat{u}_s ,
%     \\
%      y_f
%     &=
%     \left[ \begin{smallmatrix} 
%          A_x^{p_f} & A_z^{p_f}
%     \end{smallmatrix} \right]
%     \left[ \begin{smallmatrix} 
%         x_p \\ z_p
%     \end{smallmatrix} \right],
%     \\ 
%     \dot{x}_c 
%     &=
%     A_{11}^c 
%      x_c 
%     + 
%     A_{14}^c 
%    \hat{y}_f,
%     \quad 
%     u_s
%     =
%     A_x^{c_s}
%     x_c .
% \end{aligned}
% \end{equation*}
% %
Consider system \eqref{eqn: linear plant and controller} with
where $A_{11}^p = a_1$, 
$A_{12}^p = \left[ \begin{smallmatrix}
    a_2 & 0
\end{smallmatrix} \right]$,
$A_{21}^p = \left[ \begin{smallmatrix}
    0 \\ a_3
\end{smallmatrix} \right]$,
$A_{22}^p = \left[ \begin{smallmatrix}
    -a_2 & 0 \\ -a_2 & -a_4
\end{smallmatrix} \right]$,
$A_{13}^p = n_1$,
$A_{23}^p = \left[ \begin{smallmatrix}
    -n_2 \\ -n_2
\end{smallmatrix} \right]$,
$A_{x}^{p_f} = 1$,
$A_{z}^{p_f} = \left[ \begin{smallmatrix}
    0 & 1
\end{smallmatrix} \right]$,
$A_{11}^c = -a_5$,
$A_{14}^c = a_6$,
$A_{x}^{c_s}= -k$,
$a_1 =10^{-4}$, $a_2 = 0.2$, $a_3 = 0.6$, $a_4 = 0.73$, $a_5 = 1.11$, $a_6 = 0.37$, $k = 1.5$, $n_1 = 0.02$ and $n_2 = 0.0018$ are designed such that the controller stabilizes the plant under perfect communication. 

%Let $x = (x_p, x_c)$, $z = z_p$, $e_s = e_{u_s}\coloneqq \hat{u}_s - u_s$ and $e_f = e_{y_f}  \coloneqq \hat{y}_f - y_f$, then 
The hybrid model $\mathcal{H}_1$ is given by \eqref{eqn:full system}. We note that our plant and controller are simpler compared to \eqref{eqn: linear plant and controller}, since $u_f$, $y_s$ and $z_c$ does not exist in the system, nor do matrices such as $A_{14}^p$, $A_x^{p_s}$, $A_{21}^c$, etc. Consequently, the flow map $F$ in $\mathcal{H}_1$ has to be modified accordingly. Specifically, we have 
$A_{11} = \left[\begin{smallmatrix}A_{11}^p & A_{13}^p A_x^{c_s}\\ A_{14}^c A_x^{p_f} & A_{11}^c \end{smallmatrix} \right]$, 
$A_{12} = \left[\begin{smallmatrix}  A_{13}^p \\ 0\end{smallmatrix}\right]$,
$A_{13} = \left[\begin{smallmatrix} A_{12}^p  \\ A_{14}^c A_z^{p_f} \end{smallmatrix} \right]$,
$A_{14} = \left[ \begin{smallmatrix}0 \\ A_{14}^c\end{smallmatrix} \right]$,
$A_x^s = \left[\begin{smallmatrix}
     0 & -A_x^{c_s}
\end{smallmatrix} \right]$,
$A_x^f = \left[\begin{smallmatrix}
    -A_x^{p_f} & 0
\end{smallmatrix} \right]$, 
$A_z^f = -A_z^{p_f}$, 
$A_{31} = \left[\begin{smallmatrix}
    A_{21}^p & A_{23}^p A_x^{c_s}
\end{smallmatrix}\right]$,
$A_{32} = A_{23}^p$,
$A_{33} = 
    A_{22}^p$,
$A_{34} = 0$.
%with $f_x(x,z,e_s,e_f) = \big(a_1 x_p  + a_2 z_1 -n_1 k x_c + n_1e_s, -a_5 x_c + a_6 x_p+a_6 z_2+a_6 e_f \big)$, $f_{e_s}(x,z,e_s,e_f) = -a_5 k x_c + a_6 k(x_p + z_2 + e_f)$, $g_z(x,z,e_s,e_f) = (-a_2 z_1 + n_2 k x_c - n_2 e_s, a_3 x_p - a_2 z_1 - a_4 z_2 + n_2 k x_c -n_2 e_s )$ and $g_{e_f}(x,z,e_s,e_f,\epsilon) = -\epsilon (a_1 x_p  + a_2 z_1 - n_1 k x_c + n_1 e_s ) - (a_3 x_p - a_2 z_1 - a_4 z_2 + n_2 k x_c - n_2 e_s )$. 
%
% 
%By \eqref{eqn: H bar linear}, we have $\overline{H}(x,e_s) = \big(\tfrac{n_2}{a_2}(k x_c - e_s), \tfrac{a_3}{a_4}x_p\big)$. We define $y \coloneqq z - \overline{H}(x,e_s) $. 
%
%$h_y(x,e_s,y) = (y_1 - \tfrac{n_2}{a_2} e_s, y_2)$ by \eqref{eqn: h_y linear}. 
%
Then by \eqref{eqn: linear functions} we have
$A_{11}^s = \left[\begin{smallmatrix} a_1 & -\overline{n} k \\ a_6(1+\tfrac{a_3}{a_4}) & -a_5\end{smallmatrix} \right]$, 
$A_{12}^s = \left[\begin{smallmatrix} \overline{n}  \\ 0\end{smallmatrix} \right]$,
$A_{21}^s = \left[\begin{smallmatrix} a_6(1+\tfrac{a_3}{a_4})k & -a_5k \end{smallmatrix} \right]$, 
%$A_{22}^s = 0$,
$A_{11}^f = \left[\begin{smallmatrix} -a_2 & 0 \\ -a_2 & -a_4\end{smallmatrix} \right]$ and
$A_{12}^f = \left[\begin{smallmatrix} 0  \\ 0\end{smallmatrix} \right]$,
$A_{21}^f = \left[\begin{smallmatrix} a_2 & a_4 \end{smallmatrix} \right]$
%and $A_{22}^f = 0$,
where
$\overline{n} \coloneqq n_1 - n_2$.
% \begin{equation*}
%     \mathcal{H}_r \!: \! 
%     \begin{cases}
%     \begin{aligned} % Used to align \right\}
%     \left.
%     \begin{aligned} %Used to align flow map
%     \dot{x}_p &= x_p + (-k x_c + e_s)\\
%     \dot{x}_c &= 2a x_p-a x_c,\; \dot{e}_{s}=2ak x_p - ak x_c  \\
%     \dot{\tau}_s &= 1, \;  \dot{\kappa}_s = 0 \\
%     \end{aligned}
%     \right\}
%     & \begin{aligned}&\text{when } \\ & \xi^y \in \mathcal{C}_{2,r}^{y,0}\end{aligned} 
%     \\[1mm]
%     \left. 
%     \begin{aligned}
%     x^+ &= x,\;  e_{s}^+ = 0
%     \\
%     \tau_s^+ &= \tau_s, \; \kappa_s^+ = \kappa_s + 1
%     \end{aligned}
%     \qquad \qquad \qquad  \right\} 
%     &\begin{aligned}&\text{when } \\ & \xi^y \in \mathcal{D}_s^{y,0}.\end{aligned}
%     \end{aligned}
%     \end{cases}
%     %\label{eqn: example reduced system}
% \end{equation*}

%%%%%%%%%%%%%%%
%
Next, we find Lyapunov functions $W_s$ and $W_f$ in Claim \ref{Claim for LTI section}.
%
% First, we show that Assumption \ref{Assumption reduced model}, along with \eqref{eqn: Ws exponential sandwich bound} and \eqref{eqn: Vs exponential sandwich bound} in Assumption \ref{Assumption Exponential}, hold. 
% We write the flow map of $\mathcal{H}_r$ in the following state-space form:
% \begin{equation*}
%     \left[ \begin{smallmatrix} 
%         \dot{x} \\ \dot{e}_s
%     \end{smallmatrix} \right]
%     =
%     \left[ \begin{smallmatrix} 
%         A_{11}^s & A_{12}^s \\ A_{21}^s & 0
%     \end{smallmatrix} \right]
%     \left[ \begin{smallmatrix} 
%         x \\ e_s
%     \end{smallmatrix} \right],
% \end{equation*}
% where $A_{11}^s = \left[\begin{smallmatrix} a_1 & -\overline{n} k \\ a_6(1+\tfrac{a_3}{a_4}) & -a_5\end{smallmatrix} \right]$, $A_{12}^s = \left[\begin{smallmatrix} \overline{n}  \\ 0\end{smallmatrix} \right]$, $A_{21}^s = \left[\begin{smallmatrix} a_6(1+\tfrac{a_3}{a_4})k & -a_5k \end{smallmatrix} \right]$ and $\overline{n} \coloneqq n_1 - n_2$.
%
Since both $u_s$ and $y_f$ are scalars, the protocols are given by $h_s(\kappa_s, e_s) = 0$ and $h_f(\kappa_f, e_f) = 0$, which are UGES protocols. 
Let $W_s(\kappa_s,e_s) \coloneqq |e_s|$, then \eqref{eqn: NCS assumption Ws sandwich bound} and \eqref{eqn: Ws exponential sandwich bound} hold with $\underline{a}_{W_s}(s) = \overline{a}_{W_s}(s) = s$, \eqref{eqn: NCS assumption Ws jump} and \eqref{eqn: NCS Ws dot} hold for $\lambda_s = 0 $, $L_s = 0$ and $A_{H_s} = A_{21}^s$. 
Let $W_f(\kappa_f,e_f) \coloneqq |e_f|$, then \eqref{eqn: NCS assumption Wf sandwich bound} and \eqref{eqn: Wf exponential sandwich bound} hold with $\underline{a}_{W_f}(s) = \overline{a}_{W_f}(s) = s$, \eqref{eqn: NCS assumption Wf jump} and \eqref{eqn: NCS Wf dot} hold for $\lambda_f = 0 $, $L_f = 0$ and $A_{H_f} = A_{21}^f$. 
Since $W_s(x,e_s) = |e_s|$, we have $\left|\tfrac{\partial W_s(\kappa_s,e_s)}{\partial e_s} \right| \leq L_1$ where $L_1 = 1$. Then by \eqref{eqn: Lambda_b1}, we have 
$
    \Lambda_{b_1} = 2
            \left[ \begin{smallmatrix}
            a_2 |p_{11}^s| & a_6 |p_{12}^s|                    & 0 \\
            a_2 |p_{12}^s| & 0                               & a_6 |p_{22}^s| \\
            0          & \tfrac{\gamma_s}{\lambda_s^*}a_6k & \tfrac{\gamma_s}{\lambda_s^*}a_6k
        \end{smallmatrix} \right]$
%
% such that 
% \begin{equation*}
%     \begin{aligned}
%         &\Big < \tfrac{\partial {U_s}}{\partial \xi_s}, F_s^y(x,y,e_s,e_f) - F_s^y(x,0,e_s,0)  \Big> \\
%         \leq & 
%     \left[ \begin{smallmatrix} 
%         |x_1| \\ |x_2| \\ |e_s|
%     \end{smallmatrix} \right]^\top
%     \Lambda_{b_1}
%     \left[ \begin{smallmatrix} 
%         |y_1| \\ |y_2| \\ |e_f|
%     \end{smallmatrix} \right]
%         \\
%         \leq & b_1 \psi_s(|(x,e_s)|) \psi_f(|(y,e_f)|),
%     \end{aligned}
% \end{equation*}
and $b_1 = \sqrt{\lambda_{\text{max}} (\Lambda_{b_1}^\top \Lambda_{b_1})}$. Similarly, we can show that $b_2 = \sqrt{\lambda_{\text{max}}(\Lambda_{b_2}^\top \Lambda_{b_2}) }$ and $b_3 = \lambda_{\text{max}}(\Lambda_{b_3})$ satisfy \eqref{eqn: Lambda_b2 and Lambda_b3}, where %$\Lambda_{b_2}$ and $\Lambda_{b_3}$ are given below. 
$
        \Lambda_{b_2} = 2 
        \left[ \begin{smallmatrix} 
            a_1 \tfrac{a_3}{a_4} |p_{12}^f|   &  a_1 \tfrac{a_3}{a_4} |p_{22}^f| & a_1 \tfrac{\gamma_f}{\lambda_f^*} \\
            \tfrac{a_3}{a_4} \overline{n} k |p_{12}^f| & \tfrac{a_3}{a_4} \overline{n} k |p_{12}^f|& \overline{n} k \tfrac{\gamma_f}{\lambda_f^*} \\
            \tfrac{a_3}{a_4}\overline{n} |p_{12}^f|   &  \tfrac{a_3}{a_4} \overline{n}|p_{22}^f| & \overline{n}\tfrac{\gamma_f}{\lambda_f^*}
        \end{smallmatrix} \right]$ 
        and
        $\Lambda_{b_3}=
        \left[ \begin{smallmatrix} 
            2  a_2 \tfrac{a_3}{a_4} |p_{12}^f| & \star & \star \\
            a_2 \tfrac{a_3}{a_4} |p_{22}^f|   &     0      &                \star              \\
            a_2 \tfrac{\gamma_f}{\lambda_f^*} & 0 &        0              
        \end{smallmatrix} \right].$
% such that 
% \begin{equation*}
%     \begin{aligned}
%         & \Big< \tfrac{\partial {U_f}}{\partial \xi_s} - \tfrac{\partial {U_f}}{\partial y} \tfrac{\partial \overline{H}}{\partial \xi_s} - \tfrac{\partial {U_f}}{\partial e_f} \tfrac{\partial \tilde k}{\partial \xi_s} ,  F_s^y(x,y,e_s,e_f) \Big> \\
%         \leq &
%         \left[ \begin{smallmatrix} 
%             |x_1| \\ |x_2| \\ |e_s|
%         \end{smallmatrix} \right]^\top
%         \Lambda_{b_2}
%         \left[ \begin{smallmatrix} 
%             |y_1| \\ |y_2| \\ |e_f|
%         \end{smallmatrix} \right]
%         + 
%         \left[ \begin{smallmatrix} 
%         |y_1| \\ |y_2| \\ |e_f|
%         \end{smallmatrix} \right]^\top
%         \Lambda_{b_3}
%         \left[ \begin{smallmatrix} 
%             |y_1| \\ |y_2| \\ |e_f|
%         \end{smallmatrix} \right]
%         \\
%         \leq & b_2 \psi_s(|(x,e_s)|) \psi_f(|(y,e_f)|) + b_3 \psi_f^2(|(y,e_f)|),
%     \end{aligned}
% \end{equation*}
% where $b_2 = \sqrt{\lambda_{\text{max}}(\Lambda_{b_2}^\top \Lambda_{b_2}) }$, $b_3 = \lambda_{\text{max}}(\Lambda_{b_3})$, and we show \eqref{eqn: SPNCS interconnection 2} is satisfied. 
%
Finally, by \eqref{eqn: lambda_1 and lambda_2}, we show that Assumption \ref{Assumption Vf at slow transmission} holds
% Since $h_y(x,e_s,y) = (y_1 - e_s, y_2)$, we have 
% \begin{equation*}
%     \begin{aligned}
%         &V_f(x, h_y(x,e_s,y)) - V_f(x,y) \\
%         =& -2 \tfrac{n_2}{a_2} (p_{11}^f y_1 + p_{12}^f y_2) e_s +  \tfrac{n_2}{a_2}^2p_{11}^f e_s^2 \\
%         \leq & \lambda_1 W_s^2(\kappa_s,e_s) + \lambda_2 \sqrt{V_f(x,y) W_s^2(\kappa_s, e_s)} ,
%     \end{aligned}
% \end{equation*}
with $\lambda_1 = \tfrac{n_2}{a_2}^2|p_{11}^f|$ and $\lambda_2 = \tfrac{2\frac{n_2}{a_2}(|p_{11}^f| + |p_{12}^f|)}{\sqrt{\lambda_\text{min}(P_f)}}$.

%Now we have shown that Assumptions \ref{Assumption Vf at slow transmission} and \ref{Assumption interconnection Exponential} hold if Assumption \ref{Assumption reduced model}, \ref{Assumption boundary layer system} and \ref{Assumption Exponential} hold. Then we start looking for the value of unknown variables that satisfy the above-mentioned Assumptions.

We want to satisfy the LMI \eqref{eqn: LMIs} for all $\ell \in \{s,f\}$ and maximize $T(L_s,\gamma_s, \lambda_s^*)$, $T(L_f,\gamma_f, \lambda_f^*)$, $\epsilon^*$ under the following constraints: $P_s > 0$, $\gamma_\star > 0$, $a_{\rho_\star} > 0$, $\gamma_\star>0$, $\lambda_\star^* \in (0,1)$ and $\Lambda_\star < 0$ for $\star \in \{s,f\}$. We note that $\epsilon^*$ is given by \eqref{eqn: epsilon star} with $d = d^*$ defined by \eqref{eqn: d star exponential}.
We pose this problem as an optimisation problem with constraints (see \cite{Github_SPNCS_illustrative_example} for the problem fomulation and the code), we get $P_s = \left[\begin{smallmatrix} 54.91  & -1.76\\ -1.76 & 1.81\end{smallmatrix} \right]$, $\gamma_s = 2.58$, $\lambda_s^* = 0.33$, $a_{\rho_s} = 1.16$, $P_f = \left[\begin{smallmatrix} 1.12  & 0.018\\ 0.018 & 0.65\end{smallmatrix} \right]$, $\gamma_f = 0.64$, $\lambda_f^* = 0.46$, $a_{\rho_f} = 0.41$, $T(L_s,\gamma_s, \lambda_s^*) = 360.1 \ ms$ and $T(L_f,\gamma_f, \lambda_f^*) = 1.11 \ \tfrac{s}{\epsilon}$ (in fast time scale $\sigma$). By selecting $\tau_{\text{miati}}^{s} = 324.1 \ ms$ and $\mu = 0.66 a_s \underline{a}_{U_s}$, we have $\epsilon^* = 0.0162$, see the proof of Theorem \ref{Theorem Exponential decay} and \cite{Github_SPNCS_illustrative_example} for more detail.
This implies $\mathcal{H}_1$ is UGES if $\epsilon < \epsilon^*$, $\tau_{\text{miati}}^s = 324.1 \ ms$, $\tau_{\text{mati}}^s = 360.1 \ ms$ and $\tau_{\text{mati}}^f \leq 18 \ ms$. Finally, we have $\tau_{\text{miati}}^f \leq 9\ ms$ such that \eqref{eqn: condition on miati^f} is satisfied.

% \begin{figure}[H]
%     \centering
%     \includegraphics[width = 0.8\linewidth]{Figures/Illustrative Example.png}
%     \caption{Example simulation}
%     \label{fig: Example Simulation}
% \end{figure}

% Let $\phi_{\star} = -2L_\star \phi_\star - \gamma_\star(\phi_\star^2 + 1)$, $\phi_\star(0) = \lambda_\star^*$ and $U_\star = V_\star + \gamma_\star \phi_\star(\tau_\star)W_\star^2(e_\star)$ for $\star \in \{s, f\}$, and we let $\lambda_s^* = 0.4$ and $\lambda_f^* = 0.6$. By Assumptions \ref{Assumption reduced model}, \ref{Assumption boundary layer system} and the definition of $\phi_\star$, inequalities (\ref{eqn: Us}) and (\ref{eqn: Uf}) hold with $a_s=1$, $a_f = 0.5$, $\psi_s(s) = s$, $\psi_f(s) = s$.
% Then inequality (\ref{eqn: SPNCS interconnections}) holds with $b_1 = 293.61, b_2 = 11.25, b_3 = 4.36$.

% % \underline{Assumption \ref{Assumption U_f slow jump}}: Since $U_f$ is independent to $\xi_s$, Assumption \ref{Assumption U_f slow jump} immediately holds.

% Now that we have verified all the assumptions, we can compute $\epsilon^*   = 0.000151$, $T(L_s, \gamma_s, \lambda_s^*) = 0.0707$ and $T(L_f, \gamma_f, \lambda_f^*) = 0.6929$, and the required MATIs to stabilise the system are given by $\tau_{\text{mati}}^s \leq 0.0707$ and $\tau_{\text{mati}} \leq 0.6929\epsilon$.

\section{Conclusion}
This paper provided conditions ensuring various stability properties of two-time-scale singularly perturbed networked control systems (SPNCSs) under a single communication channel, including semi-global practical asymptotic stability, uniformly globally asymptotic stability, and uniformly globally exponential stability. 
We introduce a framework for analyzing SPNCSs, extendable to multiple channels and time scales, and other stability properties. Additionally, we propose a resource-aware strategy for stabilizing SPNCSs, laying a foundation for further research.

%We introduced a framework for analyzing SPNCSs that can be extended to other settings, such as multiple channels and time scales, and other stability properties. Additionally, we proposed a resource-aware strategy for stabilizing SPNCSs with a single channel. Our approach enhances the understanding and design of SPNCSs, laying a foundation for further research in more complex networked control scenarios.

% \begin{ack}                               % Place acknowledgements
% Partially supported by the Roman Senate.  % here.
% \end{ack}

\bibliographystyle{plain}        % Include this if you use bibtex 

%\bibliography{autosam}           % and a bib file to produce the 

                                    % bibliography (preferred). The
                                 % correct style is generated by
                                 % Elsevier at the time of printing.

%\begin{thebibliography}{99}     % Otherwise use the  
                                 % thebibliography environment.
                                 % Insert the full references here.
                                 % See a recent issue of Automatica 
                                 % for the style.
%  \bibitem[Heritage, 1992]{Heritage:92}
%     (1992) {\it The American Heritage. 
%     Dictionary of the American Language.}
%     Houghton Mifflin Company.
%  \bibitem[Able, 1956]{Abl:56}
%     B.~C.~Able (1956). Nucleic acid content of macroscope. 
%     {\it Nature 2}, 7--9. 
%  \bibitem[Able {\em et al.}, 1954]{AbTaRu:54}   
%     B.~C. Able, R.~A. Tagg, and M.~Rush (1954).
%     Enzyme-catalyzed cellular transanimations.
%     In A.~F.~Round, editor, 
%     {\it Advances in Enzymology Vol. 2} (125--247). 
%     New York, Academic Press.
%  \bibitem[R.~Keohane, 1958]{Keo:58}
%     R.~Keohane (1958).
%     {\it Power and Interdependence: 
%     World Politics in Transition.}
%     Boston, Little, Brown \& Co.
%  \bibitem[Powers, 1985]{Pow:85}
%     T.~Powers (1985).
%     Is there a way out?
%     {\it Harpers, June 1985}, 35--47.

%\end{thebibliography}

\appendix
%\vspace{-4mm}
\section{Proof of Theorem $\ref{Theorem H_1}$}
The conditions stated in Theorem \ref{Theorem H_1} indicate Lemma \ref{Lemma MATI} holds. 
In the proof, the following variables and functions are used.
Let 
\begin{gather}
    \mu \in (0, a_s a_{\psi_s}^2), \qquad \mu_1 \in (0, \mu),
    \label{eqn: mu}
\end{gather}
where $a_s$ and $a_{\psi_s}$ come from Lemma \ref{Lemma MATI} and Assumption \ref{Assumption Extra 2}, respectively.
Let
\begin{equation}
    \lambda \in (\exp (-\mu_1\tau_{\text{miati}}^s ), 1),
    \label{eqn: lambda}
\end{equation}
where $\tau_{\text{miati}}^{s}$ comes from Theorem \ref{Theorem H_1}. 
We define 
\begin{equation}
    d \coloneqq \tfrac{-b+\sqrt{b^2-4ac}}{2a},
    \label{eqn: d}
\end{equation}
where $a = \tfrac{\lambda_1}{\gamma_s \lambda_s^*}$, $b= \tfrac{1}{2}( \tfrac{\lambda_1}{\gamma_s \lambda_s^*} + \lambda_2)$ and $c = 1 - \lambda e^{\mu_1 \tau_{\text{miati}}^s}$. Note that by definition of $\lambda$ in \eqref{eqn: lambda}, we have $c < 0$, which implies $b^2 - 4 a c > b> 0$.

Next we define the bound when we change between $x-z$ and $x-y$ coordinate.
Since the map $\overline{H}$ between $(x,e_s)$ and the quasi-steady state is continuously differentiable and $0 =\overline{H}(0, 0)$, by \cite[Lemma 4.3]{nonlinear_systems_Khalil}, there exists a class $\mathcal{K}$ function $\zeta_1$ such that $|\overline{H}(x,e_s)| \leq \zeta_1(|(x,e_s)|)$ for all $x \in \mathbb{R}^{n_x}, e_s \in \mathbb{R}^{n_{e_s}}$. 
We define a class $\mathcal{K}_{\infty}$ function 
\begin{equation}
    \zeta_2(s) \coloneqq s + \zeta_1(s).
\end{equation}
 Let $ \mathcal{E}^y \coloneqq \{ \xi^y \in \mathbb{X}\ | \ x=0 \wedge e_s = 0 \wedge y=0 \wedge e_f = 0 \}$, and $\mathcal{E}$ be defined above Theorem \ref{Theorem H_1}. As
$
|\xi^y|_{\mathcal{E}^y} \! = \! \left| \! \big(x,e_s,y + \overline{H}(x,e_s) - \overline{H}(x,e_s),e_f \big)\! \right|
\leq   \left| \! \big(x,e_s,y + \overline{H}(x,e_s),e_f \big)\! \right| + \left|\! \big(0,0,-\overline{H}(x,e_s),0\big)\! \right|
\leq  |\xi|_{\mathcal{E}} + \zeta_1(|\xi|_{\mathcal{E}})
= \zeta_2(|\xi|_{\mathcal{E}})
$,
% \begin{equation*}
%     \begin{aligned}
%         |\xi^y|_{\mathcal{E}^y} \! &= \! \left| \! \big(x,e_s,y + \overline{H}(x,e_s) - \overline{H}(x,e_s),e_f \big)\! \right|   \\
%             &\leq   \left| \! \big(x,e_s,y + \overline{H}(x,e_s),e_f \big)\! \right| + \left|\! \big(0,0,-\overline{H}(x,e_s),0\big)\! \right| \\
%             &\leq  |\xi|_{\mathcal{E}} + \zeta_1(|\xi|_{\mathcal{E}}) \\
%             &= \zeta_2(|\xi|_{\mathcal{E}}),
%     \end{aligned}
% \end{equation*}
we have for any $\Delta>0$, there exists a $\Delta^y = \zeta_2(\Delta)$ such that $|\xi|_{\mathcal{E}} \leq \Delta$ implies $|\xi^y|_{\mathcal{E}^y} \leq \Delta^y$. Similarly, we can show  $|\xi|_{\mathcal{E}} \leq \zeta_2(|\xi^y|_{\mathcal{E}^y})$. 
Let $\nu^y \coloneqq \zeta_2^{-1}(\nu)$, then if $|\xi^y|_{\mathcal{E}^y} \leq \nu^y$, we have $|\xi|_{\mathcal{E}}\leq \zeta_2(|\xi^y|_{\mathcal{E}^y}) \leq \zeta_2(\nu^y) = \nu$.

For given $d$ coming from \eqref{eqn: d}, as well as $U_s$ and $U_f$ coming from \eqref{eqn: Us and Uf}, we define a composite Lyapunov function 
\begin{equation}
    U(\xi^y)\coloneqq {U_s}(\xi_s) + d{U_f}(\xi_s,\xi_f).
    \label{eqn: U}
\end{equation}
From (\ref{eqn: Us sandwich bound}) and (\ref{eqn: Uf sandwich bound}), we know there exist $\underline{\alpha}_U, \overline{\alpha}_U \in \mathcal{K}_{\infty}$, such that for all $\xi^y \in \mathcal{C}_1^{y,\epsilon} \cup \mathcal{D}_s^{y,\epsilon} \cup \mathcal{D}_f^{y,\epsilon} $, 
    \begin{equation}
        \underline{\alpha}_U\left( |\xi^y|_{\mathcal{E}^y} \right) \leq U(\xi^y) \leq \overline{\alpha}_U\left( |\xi^y|_{\mathcal{E}^y} \right),
        \label{eqn: No disturbance U sandwich bound}
    \end{equation} 
where the sets are defined in Sections 3 and 4.    
%where $ \mathcal{E}^y \coloneqq \{ \xi^y \in \mathbb{X}\ | \ x=0 \wedge e_s = 0 \wedge y=0 \wedge e_f = 0 \} .$

We define the $\nu_1$ and $\Delta_1$ in Assumption \ref{Assumption interconnection} as follows. 
Let
\begin{equation}
    \nu_1 \coloneqq \tfrac{\mu - \mu_1}{2 a_d}\underline{\alpha}_U(\zeta_2^{-1}(\nu)),
    \label{eqn: nu_1}
\end{equation}
where 
\begin{equation}
    a_d \coloneqq 1 + \tfrac{\lambda_1}{\gamma_s \lambda_s^*} d + \tfrac{1}{2}\left(\tfrac{\lambda_2}{\gamma_s \lambda_s^*} + \lambda_2 \right)\sqrt{d}.
    \label{eqn: a_d}
\end{equation}
We define
\begin{equation}
    \nu_2 \coloneqq \tfrac{2\nu_1}{\mu - \mu_1},
    \label{eqn: nu_2}
\end{equation}
such that by \eqref{eqn: nu_1} and \eqref{eqn: nu_2}, $U(\xi^y) \leq a_d \nu_2$ implies
% $
% U(\xi^y) \leq  a_d \nu_2 = a_d \tfrac{2\nu_1}{\mu - \mu_1}
% %
% = a_d \tfrac{2}{\mu - \mu_1} \tfrac{\mu - \mu_1}{2 a_d}\underline{\alpha}_U(\zeta_2^{-1}(\nu))
% %
% = \underline{\alpha}_U(\zeta_2^{-1}(\nu))
% $
\begin{equation}
    \begin{aligned}
        U(\xi^y) \leq & a_d \nu_2 = a_d \tfrac{2\nu_1}{\mu - \mu_1} \\
            = & a_d \tfrac{2}{\mu - \mu_1} \tfrac{\mu - \mu_1}{2 a_d}\underline{\alpha}_U(\zeta_2^{-1}(\nu))\\
            =& \underline{\alpha}_U(\zeta_2^{-1}(\nu)),
    \end{aligned}
    \label{eqn: U less than nu2 implies xi less than nu}
\end{equation}
and $|\xi^y|_{\mathcal{E}^y} \leq \zeta_2^{-1}(\nu) = \nu^y$ by \eqref{eqn: No disturbance U sandwich bound}, and $|\xi|_{\mathcal{E}} \leq \nu$.

We remind that $\Delta^y = \zeta_2(\Delta)$, and for any $|\xi|_{\mathcal{E}}\leq \Delta$, we have $|\xi^y|_{\mathcal{E}^y}\leq \Delta^y$ and $U \leq \overline{\alpha}_U(\Delta^y)$. Let 
$\Delta^U \coloneqq a_d \overline{\alpha}_U(\Delta^y)$
and $\Delta_1$ from Assumption \ref{Assumption interconnection} be 
\begin{equation}
    \Delta_1 \coloneqq \underline{\alpha}_U^{-1}(\Delta^U),
    \label{eqn: Delta_1}
\end{equation}
such that $U(\xi^y) \leq \Delta^U$ implies $|\xi^y|_{\mathcal{E}^y} \leq \Delta_1$. Next, we will bound the system trajectories. We start by showing that during flow, the derivative of the $U$ is negative definite when it is bounded within an interval. This is followed by demonstrating the boundedness during jumps.

\noindent\emph{\underline{During flow}:} Let $\xi^y \in \mathcal{C}_1^{y,\epsilon}$, where $\mathcal{C}_1^{y,\epsilon}$ is the flow set of \eqref{eqn: H_2^y}, by Lemma \ref{Lemma MATI} and (23)-(24) in \cite{teel2000assigning}, we have that
$
U^\circ(\xi^y, F^y(\xi^y, \epsilon))
\leq  \Big< \! \tfrac{\partial {U_s}}{\partial \xi_s} ,F_s^y(x,y,e_s, e_f) \!\Big> \! +\! \tfrac{d}{\epsilon} \Big<\! \tfrac{\partial {U_f}}{\partial \xi_f} , F_f^y(x,y,e_s, e_f,\epsilon) \!\Big>
+ d \Big< \tfrac{\partial {U_f}}{\partial \xi_s},F_s^y(x,y,e_s, e_f) \Big>
= \Big< \! \tfrac{\partial {U_s}}{\partial \xi_s}, F_s^y(x,0,e_s, 0) \! \Big> + \tfrac{d}{\epsilon} \Big<\! \tfrac{\partial {U_f}}{\partial \xi_f}, F_f^y(x,y,e_s, e_f,0) \!  \Big>
+ \Big< \tfrac{\partial {U_s}}{\partial \xi_s}, F_s^y(x,y,e_s,e_f) - F_s^y(x,0,e_s,0)  \Big>
+ \frac{d}{\epsilon} \Big<  \tfrac{\partial {U_f}}{\partial \xi_f}, $
$
F_f^y(x,y,e_s,e_f, \epsilon) - F_f^y(x,y,e_s,e_f,0)\Big>
+ d \Big< \tfrac{\partial {U_f}}{\partial \xi_s}, F_s^y(x,y, $ $e_s, e_f) \Big>
$.
% \begin{equation*}
%     \begin{aligned}
%     U&^\circ(\xi^y, F^y(\xi^y, \epsilon)) 
%     \\
%     \leq&  \Big< \! \tfrac{\partial {U_s}}{\partial \xi_s} ,F_s^y(x,y,e_s, e_f) \!\Big> \! +\! \tfrac{d}{\epsilon} \Big<\! \tfrac{\partial {U_f}}{\partial \xi_f} , F_f^y(x,y,e_s, e_f,\epsilon) \!\Big>
%         \\
%     &  + d \Big< \tfrac{\partial {U_f}}{\partial \xi_s},F_s^y(x,y,e_s, e_f) \Big>
%     \\
%     = &  \Big< \! \tfrac{\partial {U_s}}{\partial \xi_s}, F_s^y(x,0,e_s, 0) \! \Big> + \tfrac{d}{\epsilon} \Big<\! \tfrac{\partial {U_f}}{\partial \xi_f}, F_f^y(x,y,e_s, e_f,0) \!  \Big>
%         \\
%         & + \Big< \tfrac{\partial {U_s}}{\partial \xi_s}, F_s^y(x,y,e_s,e_f) - F_s^y(x,0,e_s,0)  \Big>
%         \\
%         & + \frac{d}{\epsilon} \Big<  \tfrac{\partial {U_f}}{\partial \xi_f}, F_f^y(x,y,e_s,e_f, \epsilon) - F_f^y(x,y,e_s,e_f,0)\Big>
%         \\
%         & + d \Big< \tfrac{\partial {U_f}}{\partial \xi_s},F_s^y(x,y,e_s, e_f) \Big>.
%     % \\
%     % \leq& (1-d) \Big< \tfrac{\partial {U_s}}{\partial \xi_s} ,F_s^y(x,y,e_s, e_f) \Big> + d \Big< \tfrac{\partial {U_f}}{\partial \xi_s},F_s^y(x,y,e_s, e_f) \Big>
%     %     \\
%     % & +\! \tfrac{d}{\epsilon} \!\Big<\! \tfrac{\partial {U_f}}{\partial \xi_f} , F_f^y(x,y,e_s, e_f,0) \!\Big>
%     % \\
%     % & - d \Big< \!\tfrac{\partial {U_f}}{\partial y}\tfrac{\partial \overline{H}}{\partial \xi_s} +\tfrac{\partial {U_f}}{\partial e_f} \tfrac{\partial \tilde k}{\partial \xi_s}, F_s^y(x,y,e_s, e_f) \! \Big>
%     \end{aligned}
%     \end{equation*}
We remind the definition of $F_f^y$ is given by \eqref{eq:functions 2},
% \begin{align*}
%     &F_f^y(x,y,e_s,e_f, \epsilon) \\
%     =& \left[ \begin{smallmatrix} 
%         g_z(x,y+\overline{H}(x,e_s),e_s,e_f)- \epsilon \tfrac{\partial \overline{H}}{\partial \xi_s} F_s^y(x,y,e_s,e_f ) 
%         \\
%         g_{e_f}(x,y+\overline{H}(x,e_s),e_s, e_f , \epsilon)
%         \\
%         1
%         \\
%         0
%     \end{smallmatrix} \right],
% \end{align*}
where the definition of $g_{e_f}$ is given in \eqref{eq:functions}.
Then we have 
\begin{equation*}
    \begin{aligned}
         &F_f^y(x,y,e_s,e_f, \epsilon) - F_f^y(x,y,e_s,e_f,0)  \\
         = & \begin{bsmallmatrix}
            - \epsilon \tfrac{\partial \overline{H}}{\partial \xi_s} F_s^y(x,y,e_s,e_f )
             \\
             -\epsilon \tfrac{\partial k_{p_f}(x_p,z_p)}{\partial x_p}  f_{x,1}(x,z,e_s,e_f) - \epsilon \tfrac{\partial k_{c_f}(x_c,z_c)}{\partial x_c}  f_{x,2}(x,z,e_s,e_f) 
             \\
             0
             \\
             0
         \end{bsmallmatrix} \\
         =& -\epsilon\Big(\tfrac{\partial \overline{H}}{\partial \xi_s} F_s^y(x,y,e_s,e_f ),  \tfrac{\partial \tilde k}{\partial \xi_s} F_s^y(x,y,e_s, e_f), 0, 0\Big),
    \end{aligned}
\end{equation*}
where $\tilde k(x,z) = \big(k_{pf}(x_p,z_p), k_{cf}(x_c,z_c)\big)$. Consequently,
$
\frac{d}{\epsilon} \! \big< \!  \tfrac{\partial {U_f}}{\partial \xi_f}, F_f^y(x,y,e_s,e_f, \epsilon) - F_f^y(x,y,e_s,e_f,0)\big> $ $= - d \big< \!\tfrac{\partial {U_f}}{\partial y}\tfrac{\partial \overline{H}}{\partial \xi_s} +\tfrac{\partial {U_f}}{\partial e_f} \tfrac{\partial \tilde k}{\partial \xi_s}, F_s^y(x,y,e_s, e_f) \! \big>
$. 
% \begin{multline*}
%      \frac{d}{\epsilon} \Big<  \tfrac{\partial {U_f}}{\partial \xi_f}, F_f^y(x,y,e_s,e_f, \epsilon) - F_f^y(x,y,e_s,e_f,0)\Big>  
%     \\
%     = - d \Big< \!\tfrac{\partial {U_f}}{\partial y}\tfrac{\partial \overline{H}}{\partial \xi_s} +\tfrac{\partial {U_f}}{\partial e_f} \tfrac{\partial \tilde k}{\partial \xi_s}, F_s^y(x,y,e_s, e_f) \! \Big>.     
% \end{multline*}
Then we have
\begin{small}
    \begin{equation}
\setlength\abovedisplayskip{-4pt}%shrink space
\setlength\belowdisplayskip{3pt}
    \begin{aligned}
        U&^\circ(\xi^y, F^y(\xi^y, \epsilon)) 
        \\
      \leq  &  \Big< \! \tfrac{\partial {U_s}}{\partial \xi_s}, F_s^y(x,0,e_s, 0) \! \Big>\! +\! \tfrac{d}{\epsilon} \Big<\! \tfrac{\partial {U_f}}{\partial \xi_f}, F_f^y(x,y,e_s, e_f,0) \!  \Big>
        \\
        & + \Big< \tfrac{\partial {U_s}}{\partial \xi_s}, F_s^y(x,y,e_s,e_f) - F_s^y(x,0,e_s,0)  \Big>
        \\
        & + d \Big<  \tfrac{\partial {U_f}}{\partial \xi_s} - \tfrac{\partial {U_f}}{\partial y} \tfrac{\partial \overline{H}}{\partial \xi_s} - \tfrac{\partial {U_f}}{\partial e_f} \tfrac{\partial \tilde k}{\partial \xi_s} ,F_s^y(x,y,e_s,e_f)\Big>.
    \end{aligned}
    \label{eqn: U derivative during flow eqn1}
\end{equation}
\end{small}

By Assumption \ref{Assumption interconnection}, (\ref{eqn: Us flow}) and (\ref{eqn: Uf flow}), we know that for $\Delta_1$ and $\nu_1$ defined in \eqref{eqn: Delta_1} and \eqref{eqn: nu_1}, there exist positive constants $b_1$, $b_2$ and $b_3$ such that for all $|\xi^y|_{\mathcal{E}^y} \leq \Delta_1$,
$
U^\circ(\xi^y, F^y(\xi^y, \epsilon))
\leq  -a_s \psi_s^2\left(\left| (x, e_s) \right|\right) - d (\tfrac{a_f}{\epsilon} - b_3) \psi_f^2\left(\left| (y, e_f) \right|\right)
+ \left(b_1 + d b_2 \right)\psi_s\left(\left| (x, e_s) \right|\right)\psi_f\left(\left| (y, e_f) \right|\right) + (1+d)\nu_1
\leq - 
     \left[ \begin{smallmatrix} 
         \psi_s(|(x,e_s)|) \\ \psi_f(|(y,e_f)|) 
     \end{smallmatrix} \right]^T \tilde{\Lambda} \left[ \begin{smallmatrix} 
         \psi_s(|(x,e_s)|) \\ \psi_f(|(y,e_f)|) 
     \end{smallmatrix} \right] + 2\nu_1
$,
% \begin{equation*}
%     \begin{aligned}
%        U&^\circ(\xi^y, F^y(\xi^y, \epsilon)) \\ 
%        \leq & -a_s \psi_s^2\left(\left| (x, e_s) \right|\right) - d (\tfrac{a_f}{\epsilon} - b_3) \psi_f^2\left(\left| (y, e_f) \right|\right)
%         \\
%             & + \left(b_1 + d b_2 \right)\psi_s\left(\left| (x, e_s) \right|\right)\psi_f\left(\left| (y, e_f) \right|\right) + (1+d)\nu_1
%         \\
%         \leq & - 
%      \left[ \begin{smallmatrix} 
%          \psi_s(|(x,e_s)|) \\ \psi_f(|(y,e_f)|) 
%      \end{smallmatrix} \right]^T \tilde{\Lambda} \left[ \begin{smallmatrix} 
%          \psi_s(|(x,e_s)|) \\ \psi_f(|(y,e_f)|) 
%      \end{smallmatrix} \right] + 2\nu_1
%     \end{aligned},
% \end{equation*}
where $\tilde{\Lambda} \coloneqq \left[ \begin{smallmatrix} 
    a_s & -\tfrac{1}{2} b_1 - \tfrac{1}{2}d b_2 \\
    -\tfrac{1}{2} b_1 - \tfrac{1}{2}d b_2 & d (\tfrac{a_f}{\epsilon} - b_3)
\end{smallmatrix} \right]$.
Then by Assumption \ref{Assumption Extra 2}, we have
$
  U^\circ(\xi^y, F^y(\xi^y, \epsilon)) \leq  - 
        \left[ \begin{smallmatrix} 
            \sqrt{U_s(\xi_s)} \\ \sqrt{U_f(\xi_s, \xi_f)}
        \end{smallmatrix} \right]^\top
        \Lambda 
$ 
$
        \left[ \begin{smallmatrix} 
            \sqrt{U_s(\xi_s)} \\ \sqrt{U_f(\xi_s, \xi_f)}
        \end{smallmatrix} \right] + 2\nu_1
$,
% \begin{equation*}
%     \begin{aligned}
%         U^\circ(\xi^y, F^y(\xi^y, \epsilon)) \leq & - 
%         \left[ \begin{smallmatrix} 
%             \sqrt{U_s(\xi_s)} \\ \sqrt{U_f(\xi_s, \xi_f)}
%         \end{smallmatrix} \right]^\top
%         \Lambda
%         \left[ \begin{smallmatrix} 
%             \sqrt{U_s(\xi_s)} \\ \sqrt{U_f(\xi_s, \xi_f)}
%         \end{smallmatrix} \right] \\
%         &+ 2\nu_1,
%     \end{aligned}
% \end{equation*}
where 
\begin{equation}
\Lambda \coloneqq 
    \left[ \begin{smallmatrix} 
        a_s a_{\psi_s}^2 & -\tfrac{1}{2}(b_1 + d b_2)a_{\psi_s}a_{\psi_f} \\
        -\tfrac{1}{2}(b_1 + d b_2)a_{\psi_s}a_{\psi_f} & d (\tfrac{a_f}{\epsilon} - b_3)a_{\psi_f}^2
    \end{smallmatrix} \right].
    \label{eqn: Lambda}
    \end{equation} 
%
%We remind that $\mu \in (0, a_s a_{\psi_s}^2)$. 
In order to satisfy $\Lambda \geq \mu
    \left[ \begin{smallmatrix} 
        1 & 0 \\ 0 & d
    \end{smallmatrix} \right]$, we need 
\begin{tiny}
\begin{subequations}
\begin{align}   
    a_s a_{\psi_s}^2 > \mu \\
    (a_s a_{\psi_s}^2-\mu) \big( d (\tfrac{a_f}{\epsilon} - b_3)a_{\psi_f}^2 - \mu d \big) & \! \geq  \! \tfrac{1}{4}(b_1 + db_2)^2a_{\psi_s}^2a_{\psi_f}^2,
    \label{eqn: inequality of epsilon Exponential}
\end{align}
\end{subequations}
\end{tiny}%
where the first inequality is satisfied by the definition of $\mu$ in \eqref{eqn: mu}, and the second inequality can be satisfied by all $\epsilon \in (0,\epsilon^*)$, where 
\begin{equation}
\hspace{-2mm}
    \epsilon^* = \left(\tfrac{a_{\psi_f}}{a_f d}\big(\tfrac{(b_1 + db_2)^2 a_{\psi_s}^2 a_{\psi_f}^2}{4(a_s a_{\psi_s}^2 - \mu)}+\mu d\big) + \tfrac{b_3 a_{\psi_f}^2}{a_f} \right)^{-1}.
    \label{eqn: epsilon star}
\end{equation}
Then we have that for all $|\xi^y|_{\mathcal{E}^y} \leq \Delta_1$, 
$
U^\circ(\xi^y, F^y(\xi^y, \epsilon)) \leq -\mu (U_s(\xi_s) + d U_f(\xi_s,\xi_f)) + 2 \nu_1 \leq - \mu U(\xi^y) + 2 \nu_1
$.
% \begin{equation*}
%     \begin{aligned}
%         U^\circ(\xi^y, F^y(\xi^y, \epsilon)) &\leq -\mu (U_s(\xi_s) + d U_f(\xi_s,\xi_f)) + 2 \nu_1 \\
%         & \leq - \mu U(\xi^y) + 2 \nu_1.
%     \end{aligned}
% \end{equation*}
%
%
By definition of $\Delta_1$ in \eqref{eqn: Delta_1}, we have $U(\xi^y) \leq \Delta^U$ implies $|\xi^y|_{\mathcal{E}^y} \leq \Delta_1$. Then by definition of $\mu_1$ in \eqref{eqn: mu} and $\nu_2$ in \eqref{eqn: nu_2}, we have
\begin{small}
\begin{equation}
\setlength\abovedisplayskip{-5pt}%shrink space
\setlength\belowdisplayskip{5pt}
    U^\circ(\xi^y, F^y(\xi^y, \epsilon)) \leq -\mu_1 U(\xi^y), \ \forall U(\xi^y) \in [\nu_2, \Delta^U].
    \label{eqn: U derivative}
\end{equation}
\end{small}

\noindent\emph{\underline{During jumps}:} 
When slow dynamics update, (\ref{eqn: Us jump}) implies that 
$
U(G_s^y(\xi^y))= {U_s}((x,h_s(\kappa_s, e_s), 0, \kappa_s+1)) + d{U_f}(G_s^y(\xi^y))\leq {U_s}(\xi_s)+ d{U_f}(G_s^y(\xi^y))
$ for all $\xi^y \in \mathcal{D}_s^{y,\epsilon}$, where $\mathcal{D}_s^{y,\epsilon}$ is the set that a slow transmission can be triggered.
% \begin{equation*}
%     \begin{aligned}
%     &U(G_s^y(\xi^y)) \\
%     &= {U_s}((x,h_s(\kappa_s, e_s), 0, \kappa_s+1)) + d{U_f}(G_s^y(\xi^y)) \\
%     &\leq {U_s}(\xi_s)+ d{U_f}(G_s^y(\xi^y))
%     \end{aligned}
% \end{equation*}
%
By definitions of $U_f$ and $G_s^y$, we have
$
U(G_s^y(\xi^y)) \leq U_s(\xi_s)+ d\big( V_f(x,h_y(\kappa_s,x,e_s,y)) + \gamma_f \phi_f(\tau) W_f^2(\kappa_f, e_f) \big).
$
% \begin{equation*}
%     \begin{aligned}
%     &U(G_s^y(\xi^y)) \\
%     &\leq U_s(\xi_s)+ d\big( V_f(x,h_y(x,e_s,y)) + \gamma_f \phi_f(\tau) W_f^2(\kappa_f, e_f) \big)
%     \end{aligned}
% \end{equation*}
By adding and subtracting the term $U_f(\xi_s,\xi_f)$ and in the view of Assumption \ref{Assumption Vf at slow transmission}, we have
$
U(G_s^y(\xi^y)) \leq  U(\xi^y) + d\Big(V_f\big(x, h_y(\kappa_s,x,e_s,y)\big) - V_f(x,y)\Big)
\leq  U(\xi^y) + d\Big( \lambda_1 W_s^2(\kappa_s,e_s) + \lambda_2 \sqrt{W_s^2(\kappa_s,e_s) V_f(x,y)} \Big)
$.
% \begin{equation*}
%     \begin{aligned}
%     U(G_s^y(\xi^y)) \leq & U(\xi^y) + d\Big(V_f\big(x, h_y(x,e_s,y)\big) - V_f(x,y)\Big) \\
%     \leq & U(\xi^y) + d\Big( \lambda_1 W_s^2(\kappa_s,e_s) \\
%         &+ \lambda_2 \sqrt{W_s^2(\kappa_s,e_s) V_f(x,y)} \Big).
%                \end{aligned} %\label{eqn: U slow jump part 2}
% \end{equation*}
By completion of square, we have $\sqrt{W_s^2(\kappa_s,e_s) V_f(x,y)} \leq \tfrac{1}{2\sqrt{d}}W_s^2(\kappa_s,e_s) + \tfrac{\sqrt{d}}{2}V_f(x,y)$. At the same time, $W_s^2(\kappa_s,e_s) \leq \tfrac{1}{\gamma_s \lambda_s^*}U(\xi^y)$ and $d V_f(x,y) \leq U(\xi^y)$ by definition of $U(\xi^y)$. Then we have
\begin{equation}
    \begin{aligned}
        U(G_s^y(\xi^y)) \leq & U(\xi^y) + d \lambda_1 W_s^2(\kappa_s,e_s)  + \tfrac{1}{2\sqrt{d}} \lambda_2 W_s^2(\kappa_s,e_s)\\ 
            &  + \lambda_2 \tfrac{\sqrt{d}}{2} \big(d V_f(x,y) \big)
        \\
         \leq & \left(1 + \tfrac{\lambda_1}{\gamma_s \lambda_s^*} d + \tfrac{1}{2}\left(\tfrac{\lambda_2}{\gamma_s \lambda_s^*} + \lambda_2 \right)\sqrt{d}  \right) U(\xi^y)
         \\
         = & a_d U(\xi^y),
    \end{aligned}
    \label{eqn: U slow jump}
\end{equation}
where $a_d$ is defined in \eqref{eqn: a_d}.
For fast dynamics updates, (\ref{eqn: Uf jump}) implies for all $\xi^y \in \mathcal{D}_f^{y,\epsilon}$,
\begin{equation}
    \begin{aligned}
    U(G_f^y(\xi^y)) &= {U_s}(\xi_s) + d{U_f}(G_f^y(\xi^y))\\
    &\hspace{-5mm} \leq {U_s}(\xi_s) + d{U_f}(\xi_s, \xi_f) = U(\xi^y) ,
    \end{aligned}
    \label{eqn: U fast update}
\end{equation}
where we can see $U$ is non-increasing at fast transmissions.
Let $\xi^y(t,j)$ be a solution to $\mathcal{H}_1^y$, $(t, j)\in\text{dom}\, \xi^y$ and $0 = t_0 \leq t_1 \leq \cdots \leq t_{j+1} = t$ satisfying 
$$
\text{dom} \ \xi^y \cap ([0,t]\times \{0,\cdots,j \}) = \bigcup_{i\in \{0,\cdots,j\}  } [t_i,t_{i+1}] \times \{ i\}.
$$
%We note the sequence $\{t_1, \cdots, t_j \}$ are made of two disjoint sub-sequences, which correspond to transmission times of slow and fast variables, respectively, and are denoted by $\{t_{j_k^s}\}$ and $\{t_k^f\}$, respectively. 
The sequence $\{1,\cdots,j \}$ is divided into two disjoint sub-sequences representing slow and fast transmissions, respectively, and are denoted by $\mathcal{J}^s \coloneqq \{j_1^s, j_2^s, \cdots ,j_m^s \}$ and $\{j_1^f, j_2^f, \cdots ,j_n^f \}$, respectively.
We will first focus on the trajectory of $U(\xi^y(s,i))$, at each $i \in \mathcal{J}^s$ and all $s\in t_i $. Note that by abuse of notation, we write $U(s,i) \coloneqq U(\xi^y(s,i))$.
%
%For a sequence $\left\{j  \right\}_{j=0}^{\infty}$ of both slow and fast jumps, we pick all slow jumps to form a subsequence, denoted as $\left\{j_k^s  \right\}_{k=0}^{\infty}$ where $j_k^s$ denotes $k^{\text{th}}$ slow jump at time $t_{j_k^s} \in \mathcal{T}^s$. 
The notation $(s_1,i_1)\preceq (s_2,i_2)$ indicates $s_1+i_1 \leq s_2+i_2$, while $(s_1,i_1)\prec (s_2,i_2)$ implies $s_1+i_1 < s_2+i_2$. 
\begin{claim}
    %If $(t_a, j_a)\in \text{dom} \ \xi^y$ satisfies $U(t_a,j_a)\leq \nu_2$ and $(t_{j_k^s}, j_k^s) \preceq (t_a,j_a) \preceq (t_{j_{k+1}^s}, j_{k+1}^s-1)$, then $U(s,i)\leq \nu_2$ for all $(t_a, j_a) \preceq (s,i) \preceq (t_{j_{k+1}^s}, j_{k+1}^s-1)$.
    Suppose $(t_a, j_a) \preceq (t,j) \in \text{dom} \ \xi^y$ satisfies $U(t_a,j_a)\leq \nu_2$. Let $\mathcal{J}^s_{>a} \coloneqq \{i \in \mathcal{J}^s | i > j_a \}$.
    If $\mathcal{J}^s_{>a} = \emptyset$, then $U(s,i)\leq \nu_2$ for all $(t_a, j_a) \preceq (s,i) \preceq (t, j)$. Otherwise, let $j_{k}^s = \min \mathcal{J}^s_{>a}  $, then $U(s,i)\leq \nu_2$ for all $(t_a, j_a) \preceq (s,i) \preceq (t_{j_{k}^s}, j_{k}^s-1)$.
    \label{claim: invariant during flow}
\end{claim}
\textbf{Proof:} We first proof Claim \ref{claim: invariant during flow} for the case $\mathcal{J}^s_{>a} \neq \emptyset$ by contradiction. Assume there exist $(t_b, j_b)$ such that $(t_a, j_a)\preceq (t_b, j_b) \preceq (t_{j_{k}^s}, j_{k}^s-1)$ and $U(t_b,j_b) > \nu_2$. Then by continuity of $U$, there exist $(t_c,j_c)$ such that $(t_a, j_a)\preceq (t_c, j_c) \preceq (t_b, j_b)$, $U(t_c,j_c) = \nu_2$ and $U(s,i) \geq \nu_2 $ for all $(t_c, j_c)\preceq (s, i) \preceq (t_b, j_b)$. 
Then by \eqref{eqn: U derivative}, we have $ U(t_b,j_b) \leq U(t_c,j_c) + \int_{t_c}^{t_b} U^\circ(\xi^y, F^y(\xi^y, \epsilon)) dt \leq \nu_2 -\mu_1 \nu_2  (t_b-t_c) \leq \nu_2$,
% \begin{equation*}
%     \begin{aligned}
%         U(t_b,j_b) &= U(t_c,j_c) + \int_{t_c}^{t_b} U^\circ(\xi^y, F^y(\xi^y, \epsilon)) dt \\
%                     &\leq \nu_2 -\mu_1 \nu_2  (t_b-t_c) \leq \nu_2,
%     \end{aligned}
% \end{equation*}
which conflict our assumption that $U(t_b,j_b) > \nu_2$, and we proved Claim \ref{claim: invariant during flow} for the case $\mathcal{J}^s_{>a} \neq \emptyset$. The case $\mathcal{J}^s_{>a} = \emptyset$ can be proved similarly.

%%%%%%%%%%%%%%%%%%%%

Suppose $j_k^s$, $j_{k+1}^s \in\mathcal{J}^s$. By \eqref{eqn: U derivative}, \cite{teel2000assigning}, comparison principle \cite[Lemma 3.4]{nonlinear_systems_Khalil} and the fact that $U$ is non-increasing at fast transmissions, we have that if $U(t_{j_k^s}, j_k^s) \in [\nu_2, \Delta^U]$, then
\begin{equation}
    U(s,i) \leq U(t_{j_k^s}, j_k^s) \exp \! \big(-\mu_1(s-t_{j_k^s}) \big) \label{eqn: Exponential U flow}, 
\end{equation}
whenever $U(s,i) \geq \nu_2$, as well as $(t_{j_k^s}, j_k^s)\! \preceq \! (s,i) \! \preceq \!(t_{j_{k+1}^s}, j_{k+1}^s - 1)$ or $(t_{j_m^s}, j_m^s)\! \preceq \! (s,i) \! \preceq \!(t, j)$, where $j_m^s$ is the last element in $\mathcal{J}^s$ and $m$ corresponds to the amount of slow transmissions in trajectory $\xi^y(t,j)$.

Since $t_{j_{k+1}^s} - t_{j_k^s} \geq \tau_{\text{miati}}^s$ for all $j_k^s,j_{k+1}^s\in\mathcal{J}^s$, we have that $U(t^s_{k+1},j^s_{k+1}-1) \leq U(t^s_k, j^s_k) \exp \! \big(-\mu_1 \tau_{\text{miati}}^s \big)$
% \begin{equation}
%     U(t^s_{k+1},j^s_{k+1}-1) \leq U(t^s_k, j^s_k) \exp \! \big(-\mu_1 \tau_{\text{miati}}^s \big), 
% \end{equation}
if $U(s,i)\geq \nu_2$ for all $(t_{j_k^s}, j_k^s) \preceq (s,i) \preceq (t_{j_{k+1}^s}, j_{k+1}^s - 1)$ and $U(t_{j_k^s}, j_k^s) \in [\nu_2, \Delta^U]$. Then we will experience a slow transmission, by \eqref{eqn: U slow jump}, we have for all $j_k^s,j_{k+1}^s\in\mathcal{J}^s$ and $U(t_{j_k^s}, j_k^s) \in [\nu_2, \Delta^U]$, if $U(s,i)\geq \nu_2$ for all $(t_{j_k^s}, j_k^s) \preceq (s,i) \preceq (t_{j_{k+1}^s}, j_{k+1}^s - 1)$, then
\begin{equation}
    \begin{aligned}
        U(t^s_{k+1}, j^s_{k+1}) &\leq a_d U(t^s_{k+1},j^s_{k+1} - 1)
        \\
        &\leq a_d U(t_{j_k^s}, j_k^s) \exp (-\mu_1\tau_{\text{miati}}^s ).
    \end{aligned} 
    \label{eqn: flow then jump}
\end{equation}
By definition of $\lambda$, $d$ and $a_d$ in \eqref{eqn: lambda}, \eqref{eqn: d} and \eqref{eqn: a_d}, we have for all $U(t_{j_k^s}, j_k^s) \in [\nu_2, \Delta^U]$
\begin{equation}
\begin{aligned}
    U(t^s_{k+1}, j^s_{k+1}) &\leq a_d U(t_{j_k^s}, j_k^s) \exp (-\mu_1\tau_{\text{miati}}^s ) \\
    &= \lambda U(t^s_k,j^s_k)
\end{aligned}
    \label{eqn: U slow jump with lambda}
\end{equation}
if $U(s,i)\geq \nu_2$ for all $(t_{j_k^s}, j_k^s) \preceq (s,i) \preceq (t_{j_{k+1}^s}, j_{k+1}^s - 1)$.

The following claim provides an upper bound of $U(\xi^y(t,j))$ for all $(t, j) \in \text{dom}\ \xi^y$. 
% \begin{claim}
%     Given $|\xi(0,0)|_{\mathcal{E}} \leq \Delta$, we have 
%     \begin{equation}
%         U(t,j) \leq \frac{a_d}{\lambda}U(t_{j_k^s},j_k^s)  \exp \!\left(-\tfrac{\ln{(\nicefrac{1}{\lambda})}}{\tau_{\text{mati}}^s}(t-t_{j_k^s}) \right) + a_d \nu_2,
%         \label{eqn: upper bound of U}
%     \end{equation}
%     for all $(t,j) \in \text{dom}\ \xi$, $k \in \mathbb{Z} and \red{(t_k)}$.
%     \label{Claim U upper bound}
% \end{claim}
\begin{claim}
    Given $|\xi(0,0)|_{\mathcal{E}} \leq \Delta$, we have 
    \begin{equation}
        U(t, j) \leq \max \left\{ \frac{a_d}{\lambda}U(0,0)  \exp \!\left(-\tfrac{\ln{(\nicefrac{1}{\lambda})}}{\tau_{\text{mati}}^s}(t) \right), a_d \nu_2\right\},
        \label{eqn: upper bound of U}
    \end{equation}
    for all $(t, j) \in \text{dom}\ \xi^y$.
    \label{Claim U upper bound}
\end{claim}
\textbf{Proof of Claim \ref{Claim U upper bound}:}
We first consider the case that $\mathcal{J}^s$ is non-empty, where the scenario with $\mathcal{J}^s = \emptyset$ can be deduced from step 1 below directly.
The proof has two steps, in the first step, we check the upper bound of $U(s,i)$ for all $(s,i) \preceq (t_{j_1^s}, j_1^s)$ and $(s,i) \in \text{dom}\,\xi^y $, then we verify the upper bound of $U(t,j)$ in the second step.
Here we remind that  $|\xi(0,0)|_{\mathcal{E}} \leq \Delta$ implies $U(0,0) \leq \overline{\alpha}_U(\Delta^y) < \Delta^U$ and $\Delta^U = a_d \overline{\alpha}_U(\Delta^y)$.

\emph{\underline{Step 1} }: Since the initial condition of $\tau_s$ is not guaranteed to be zero, $t_{j_1^s}$ (the time of first slow transmission) might be smaller than $\tau_{\text{miati}}^s$. In this step, we focus on the interval $(s,i) \preceq (t_{j_1^s}, j_1^s)$, $(s,i) \in \text{dom}\,\xi^y $.  
We divide into two cases based on the initial condition.

\noindent \emph{\underline{Case 1-1}}: If $U(0,0) \leq \nu_2$, by Claim \ref{claim: invariant during flow}, we know $U(s,i) \leq \nu_2$ for all $(s,i) \preceq (t_{j_1^s}, j_1^s-1)$. Then by \eqref{eqn: U slow jump}, $U(t_{j_1^s}, j_1^s) \leq a_d \nu_2$.

\noindent \emph{\underline{Case 1-2} }: If $U(0,0) \in (\nu_2, \overline{\alpha}_U(\Delta^y)]$, by \eqref{eqn: U derivative}, \eqref{eqn: U fast update} and Claim \ref{claim: invariant during flow}, we have $U(s,i) \leq U(0,0)$ for all $(s,i) \preceq (t_{j_1^s}, j_1^s - 1)$. Then by \eqref{eqn: U slow jump} we have $U(s,i) \leq a_d U(0,0)$ for all $(s,i) \preceq (t_{j_1^s}, j_1^s)$.

For any $(t,j)\in \text{dom}\ \xi^y$ with $\mathcal{J}^s = \emptyset$, Claim \ref{claim: invariant during flow} implies $U(t,j) \leq \nu_2$ if $U(0,0)\leq \nu_2$, and Case 1-2 shows $U(t,j)\leq U(0,0)$ if $U(0,0) \in (\nu_2, \overline{\alpha}_U(\Delta^y)]$. Since $t \leq \tau_{\text{mati}}^s$, we have $ \frac{a_d}{\lambda}U(0,0)  \exp \!\left(-\tfrac{\ln{(\nicefrac{1}{\lambda})}}{\tau_{\text{mati}}^s}(t) \right)   \geq U(0,0)$, which shows \eqref{eqn: upper bound of U} hold for all $(t,j)\in \text{dom}\ \xi^y$ when $\mathcal{J}^s = \emptyset$.
%
% we can show \eqref{Claim U upper bound} hold. Since $t_{j_1^s} \leq \tau_{\text{mati}}^s$, $t \leq \tau_{\text{mati}}^s$ we have 
% $ \frac{a_d}{\lambda}U(0,0)  \exp \!\left(-\tfrac{\ln{(\nicefrac{1}{\lambda})}}{\tau_{\text{mati}}^s}(s) \right)   \geq a_d U(0,0)$ for all $s \in [0, t_{j_1^s}]$. By considering both cases, we can see inequality \eqref{eqn: upper bound of U} is satisfied for all $(s,i) \preceq (t_{j_1^s}, j_1^s)$. 
At the same time, from Cases 1-1 and 1-2, we can see $|\xi(0,0)|_{\mathcal{E}} \leq \Delta$ implies $U(t_{j_1^s}, j_1^s) \leq a_d \overline{\alpha}_U = \Delta^U $.

\emph{\underline{Step 2} }: In this step, we verify \eqref{eqn: upper bound of U} for $(t,j)$ with $\mathcal{J}^s \neq \emptyset$. We also divide into two cases based on $U(t_{j_1^s}, j_1^s)$.

\noindent \emph{\underline{Case 2-1} }:We consider the case that $U(t_{j_1^s}, j_1^s) \leq a_d \nu_2$.
We remind that $m$ is the amount of slow transmissions in the trajectory $\xi^y(t,j)$.
If $m=1$, by \eqref{eqn: U derivative} and Claim \eqref{claim: invariant during flow}, we have $U(t,j)\leq a_d \nu_2$.
Next we consider the scenario that $m\geq 2$.
Let $\Omega_{\nu_2} \coloneqq \{\xi^y \in \mathbb{X} | U(\xi^y) \leq \nu_2\}$.
By \eqref{eqn: Exponential U flow}, definition of $\lambda$, $d$ and $a_d$ in \eqref{eqn: lambda}, \eqref{eqn: d} and \eqref{eqn: a_d}, we can see the time that takes $\xi^y(t_{j_1^s}, j_1^s)$ to enter the set $\Omega_{\nu_2}$ is less than or equal to $\tau_{\text{miati}}^s$. Then by Claim \ref{claim: invariant during flow}, we have $U(t_2^s, j_2^s-1) \leq \nu_2$. By \eqref{eqn: U slow jump}, we have $U(t_{j_2^s}, j_2^s) \leq a_d \nu_2$. Finally, by concatenation, we can prove $U(t,j)\leq a_d \nu_2$.
%
% which implies $U(t_2^s, j_2^s-1) \leq \nu_2$ if $m\geq 2$. By \eqref{eqn: U derivative} and and Claim \eqref{claim: invariant during flow}, we have $U(t,j)\leq a_d \nu_2$ if $m=1$. If $m\geq 2$, by \eqref{eqn: U slow jump}, we have $U(t_{j_2^s}, j_2^s) \leq a_d \nu_2$. Finally, by concatenation, we can prove $U(t,j)\leq a_d \nu_2$.

\noindent \emph{\underline{Case 2-2} }:
We consider the case that $U(t_{j_1^s}, j_1^s) \in (a_d \nu_2, \Delta^U] $, which can exist only if $U(0,0) \in (\nu_2, \overline{\alpha}_U(\Delta^y)]$. Then by Cases 1-2, we know $U(t_{j_1^s}, j_1^s) \leq a_d U(0,0)$. We consider two possible scenarios.

Firstly, if $U(s,i) > \nu_2$ for all $(t_{j_1}^s, j_1^s) \preceq (s,i) \preceq (t,j)$, then by \eqref{eqn: U slow jump with lambda} and concatenation, we have $U(t_{j_m^s}, j_m^s) \leq \lambda^{m-1} U(t_{j_1^s},j_1^s) \leq \lambda^{m-1} a_d U(0,0)$,
% \begin{equation}
%     U(t_{j_m^s}, j_m^s) \leq \lambda^{m-1} U(t_{j_1^s},j_1^s) \leq \lambda^{m-1} a_d U(0,0) , \label{eqn: Discrete model is UGES}
% \end{equation}
where $j_m^s$ is the last element in $\mathcal{J}^s$.
Then by \eqref{eqn: Exponential U flow}, we have 
\begin{equation}
\setlength\belowdisplayskip{3pt}
U(t,j) \leq \lambda^{m-1} a_d U(0,0) \text{exp}\big(-\mu_1(t-t_{j_m^s})\big).
\label{eqn: U(t,j) upperbournd case 2-2}
\end{equation}

The first component on the right-hand side of \eqref{eqn: upper bound of U} can be written as
\begin{equation}
    \begin{aligned}
        &\frac{a_d}{\lambda}U(0,0)  \exp \!\left(-\tfrac{\ln{(\nicefrac{1}{\lambda})}}{\tau_{\text{mati}}^s}t \right) \\
        =& \frac{a_d}{\lambda}U(0,0)  \exp \!\left(-\tfrac{\ln{(\nicefrac{1}{\lambda})}}{\tau_{\text{mati}}^s}t_{j_m^s} \right) \exp \!\left(-\tfrac{\ln{(\nicefrac{1}{\lambda})}}{\tau_{\text{mati}}^s} (t-t_{j_m^s}) \right).
    \end{aligned}
    \label{eqn: verify upper bound part 1}
\end{equation}
Since $\tfrac{t_{j_m^s}}{\tau_{\text{mati}}^s} \leq m$, we have 
$\frac{a_d}{\lambda}U(0,0)  \exp \!\left(-\tfrac{\ln{(\nicefrac{1}{\lambda})}}{\tau_{\text{mati}}^s}t_{j_m^s} \right) \geq \frac{a_d}{\lambda}U(0,0)  \exp \!\left(-\ln{(\nicefrac{1}{\lambda})} m \right) =  a_d U(0,0) \lambda^{m-1}$.
%
% \begin{equation}
%     \begin{aligned}
%         & \frac{a_d}{\lambda}U(0,0)  \exp \!\left(-\tfrac{\ln{(\nicefrac{1}{\lambda})}}{\tau_{\text{mati}}^s}t_{j_m^s} \right) \\
%         \geq & \frac{a_d}{\lambda}U(0,0)  \exp \!\left(-\ln{(\nicefrac{1}{\lambda})} m \right) 
%         =  a_d U(0,0) \lambda^{m-1}.
%     \end{aligned} \label{eqn: verify discrete upper bound}
% \end{equation}
%
By definition of $\lambda$ in \eqref{eqn: lambda}, and since $a_d >1$, we have
$\exp \!\left(-\tfrac{\ln{(\nicefrac{1}{\lambda})}}{\tau_{\text{mati}}^s} (t-t_{j_m^s}) \right) 
=  \lambda^{\tfrac{t-t_{j_m^s}}{\tau_{\text{mati}}^s}} 
\geq   \big(a_d \exp (-\mu_1  \tau_{\text{miati}}^s $ $) \big)^{\tfrac{t-t_{j_m^s}}{\tau_{\text{mati}}^s}}
\geq  \exp \! \big(-\mu_1 \tfrac{\tau_{\text{miati}}^s}{\tau_{\text{mati}}^s} (t - t_{j_m^s})\big)
\geq  \exp \! \big(-\mu_1 (t - t_{j_m^s})\big)
$.
%
% \begin{equation}
%     \begin{aligned}
%         & \exp \!\left(-\tfrac{\ln{(\nicefrac{1}{\lambda})}}{\tau_{\text{mati}}^s} (t-t_{j_m^s}) \right) \\
%         %=&\exp \!\left( -\ln{(\nicefrac{1}{\lambda})}\right)^{\tfrac{t-t_{j_k^s}}{\tau_{\text{mati}}^s}} 
%         = & \lambda^{\tfrac{t-t_{j_m^s}}{\tau_{\text{mati}}^s}} 
%         \geq  \big(a_d \exp (-\mu_1 \tau_{\text{miati}}^s) \big)^{\tfrac{t-t_{j_m^s}}{\tau_{\text{mati}}^s}} \\
%         \geq & \exp \! \big(-\mu_1 \tfrac{\tau_{\text{miati}}^s}{\tau_{\text{mati}}^s} (t - t_{j_m^s})\big) \\
%         \geq & \exp \! \big(-\mu_1 (t - t_{j_m^s})\big). \\
%     \end{aligned}
%     \label{eqn: verify upper bound part 2}
% \end{equation}
Then by \eqref{eqn: U(t,j) upperbournd case 2-2} and \eqref{eqn: verify upper bound part 1}, we show that if $U(s,i) > \nu_2$ for all $(t_{j_1}^s, j_1^s) \preceq (s,i) \preceq (t,j)$, \eqref{eqn: upper bound of U} hold. 

Secondly, we consider the scenario that there exist $(t_a,j_a) \preceq (t,j)$, $(t_a,j_a) \in \text{dom} \ \xi^y$ such that $U(t_a,j_a) \leq \nu_2$. Following the similar procedure as in Case 2-1, we can prove $U(t,j)\leq a_d \nu_2$. Then we have proved Claim \ref{Claim U upper bound}.   $\hfill\square$

Now we convert the upper bound on $U(\xi^y(t,j))$ to the upper bound on $|\xi|_{\mathcal{E}}$.  By \eqref{eqn: U less than nu2 implies xi less than nu}, we know $U(\xi^y) \leq a_d \nu_2$ implies $|\xi|_{\mathcal{E}} \leq \nu$. Additionally, since $|\xi|_{\mathcal{E}} \leq \zeta_2(|\xi^y|_{\mathcal{E}^y})$, $|\xi^y|_{\mathcal{E}^y} \leq \zeta_2(|\xi|_{\mathcal{E}})$ and by sandwich bound \eqref{eqn: No disturbance U sandwich bound}, we have $U(t,j) \leq  \frac{a_d}{\lambda}U(0,0)  \exp \!\left(-\tfrac{\ln{(\nicefrac{1}{\lambda})}}{\tau_{\text{mati}}^s}t \right)$ implies $|\xi(t,j)| \leq \beta_1(\xi(0,0), t)$ where 
$
\beta_1(s, t) \coloneqq \zeta_2(\underline{\alpha}_{U}^{-1}(\frac{a_d}{\lambda}\overline{\alpha}_U(\zeta_2(s))  \exp \!\big(-\tfrac{\ln{(\nicefrac{1}{\lambda})}}{\tau_{\text{mati}}^s} t \big)))
$
% \begin{equation*}
%     \beta_1(s, t) \coloneqq \zeta_2(\underline{\alpha}_{U}^{-1}(\frac{a_d}{\lambda}\overline{\alpha}_U(\zeta_2(s))  \exp \!\left(-\tfrac{\ln{(\nicefrac{1}{\lambda})}}{\tau_{\text{mati}}^s}(t) \right))),
% \end{equation*}
and $\beta_1 \in \mathcal{KL}$. Then by Claim \ref{Claim U upper bound}, we know that $|\xi(t,j)|_\mathcal{E} \leq \beta_1(|\xi(0,0)|_\mathcal{E}, t) + \nu$ for all $|\xi(0,0)|_\mathcal{E}\leq \Delta$ .
% \begin{equation}
%     |\xi(t,j)|_\mathcal{E} \leq \beta_1(|\xi(0,0)|_\mathcal{E}, t) + \nu.
% \end{equation}
%
Additionally, since $t \geq \tau_{\text{miati}}^s (j-1)$, we have 
%
% $
% \exp \!\left(-\tfrac{\ln{(\nicefrac{1}{\lambda})}}{ \tau_{\text{mati}}^s}(t) \right)
% \leq  \exp \!\left(-\tfrac{\ln{(\nicefrac{1}{\lambda})}}{\tau_{\text{mati}}^s}(\tfrac{t}{2}+\tfrac{\tau_{\text{miati}}^s}{2}(j-1))) \right)
% \leq  \exp \!\left( \tfrac{\ln{(\nicefrac{1}{\lambda})}\tau_{\text{miati}}^s}{2 \tau_{\text{mati}}^s}\right)
% $
% $
% \exp \!\left(-\tfrac{\ln{(\nicefrac{1}{\lambda})}}{2 \tau_{\text{mati}}^s}  \min\{1,\tau_{\text{miati}}^s \} (t+j) \right)
% \eqqcolon  \alpha_1(t+j)           
% $.
\begin{equation}
    \begin{aligned}
        &\exp \!\left(-\tfrac{\ln{(\nicefrac{1}{\lambda})}}{ \tau_{\text{mati}}^s}(t) \right) \\
        %=&  \exp \!\left(-\tfrac{\ln{(\nicefrac{1}{\lambda})}}{\tau_{\text{mati}}^s}(\tfrac{t}{2}+\tfrac{t}{2})) \right) \\
        \leq & \exp \!\left(-\tfrac{\ln{(\nicefrac{1}{\lambda})}}{\tau_{\text{mati}}^s}(\tfrac{t}{2}+\tfrac{\tau_{\text{miati}}^s}{2}(j-1))) \right) \\
        % =& \exp \!\left( \tfrac{\ln{(\nicefrac{1}{\lambda})}\tau_{\text{miati}}^s}{2\tau_{\text{mati}}^s}\right)   
        %     \exp \!\left(-\tfrac{\ln{(\nicefrac{1}{\lambda})}}{2 \tau_{\text{mati}}^s}(t+\tau_{\text{miati}}^sj) \right) \\
        \leq & \exp \!\left( \tfrac{\ln{(\nicefrac{1}{\lambda})}\tau_{\text{miati}}^s}{2 \tau_{\text{mati}}^s}\right)   
            \exp \!\left(-\tfrac{\ln{(\nicefrac{1}{\lambda})}}{2 \tau_{\text{mati}}^s}  \min\{1,\tau_{\text{miati}}^s \} (t+j) \right)
            \\
            \eqqcolon & \alpha_1(t+j).
    \end{aligned}
    \label{eqn: change t to t+j}
\end{equation}
Then we have for all $|\xi(0,0)|_\mathcal{E} \leq \Delta$, there exist $\beta \in \mathcal{KL}$ such that $|\xi(t,j)|_\mathcal{E} \leq \beta(|\xi(0,0)|_\mathcal{E}, t+j) + \nu$,
% \begin{equation*}
%      |\xi(t,j)|_\mathcal{E} \leq \beta(|\xi(0,0)|_\mathcal{E}, t+j) + \nu,
% \end{equation*}
where $\beta(s,t+j) = \zeta_2(\underline{\alpha}_{U}^{-1}(\frac{a_d}{\lambda}\overline{\alpha}_U(\zeta_2(s))  \alpha_1 (t+j)))$.

\section{Proof of Theorem \ref{Theorem Exponential decay}} 
The first step is to show the $\psi_s$ and $\psi_f$ in Lemma \ref{Lemma MATI} are linear. By \eqref{eqn: NCS Ws dot}, \eqref{eqn: NCS Vs flow} and the definition of $U_s$ in \eqref{eqn: definition of U_s}, we can show
$
U_s^\circ(\xi_s; F_s^y(x,0,e_s, 0)) \leq - \rho_s(|x|) - \rho_s\left(W_s(\kappa_d, e_s)\right)
$
along the same line as \cite[(27)]{carnevale_stability}.
%
%
% \begin{equation*}
%     \begin{aligned}
%         & U_s^\circ(\xi_s; F_s^y(x,0,e_s, 0)) \\
%         \leq &\left< \tfrac{\partial V_s}{\partial x}, f_x(x,\overline{H}(x,e_s),e_s,0)\right> \\
%             & + \gamma_d \Big( -2 L_s \phi_s(\tau_s) - \gamma_d \big(\phi_s^2(\tau_s)+1\big) \Big) W_s^2(\kappa_d, e_s) \\
%             &+ 2 \gamma_d\phi_s W_s(\kappa_d,e_s) \left< \tfrac{\partial {W_s}(\kappa_d,{e_s})}{\partial {e_s}}, f_{e_s}(x,\overline{H}(x,e_s),e_s, 0)\right> \\
%         \leq & - \rho_s(|x|) - \rho_s\left(W_s(\kappa_d, e_s)\right) - H_s^2(x,e_s) + \gamma_d^2 W_s^2(\kappa_d, e_s) \\
%             &+ \gamma_d \Big( -2 L_s \phi_s(\tau_s) - \gamma_d \big(\phi_s^2(\tau_s)+1\big) \Big) W_s^2(\kappa_d, e_s) \\
%             &+ 2 \gamma_d \phi_s(\tau_s) W_s(\kappa_d, e_s)\big(L_s {W_s}(\kappa_d, e_s)  + H_s(x,e_s) \big)  \\
%         \leq & - \rho_s(|x|) - \rho_s\left(W_s(\kappa_d, e_s)\right) \\
%             &- \big(H_s(x,e_s) - \gamma_d \phi_s(\tau_s)W_s(\kappa_d, e_s)\big)^2 \\
%         \leq & - \rho_s(|x|) - \rho_s\left(W_s(\kappa_d, e_s)\right).
%     \end{aligned}
% \end{equation*}
Additionally, since $a_{\rho_s} s^2 \leq \rho_s(s)$ for all $s \in \mathbb{R}$, we have $U_s^\circ(\xi_s; F_s^y(x,0,e_s, 0)) \leq -a_{\rho_s} |x|^2 - a_{\rho_s} W_s^2(\kappa_d, e_s)$.
% \begin{equation}
%     U_s^\circ(\xi_s; F_s^y(x,0,e_s, 0)) \leq -a_{\rho_s} |x|^2 - a_{\rho_s} W_s^2(\kappa_d, e_s)
% \end{equation}
Then by \eqref{eqn: Ws exponential sandwich bound}, we have
$
U_s^\circ(\xi_s; F_s^y(x,0,e_s, 0)) \leq  - \big(a_{\rho_s} |x|^2 + a_{\rho_s} W_s^2(\kappa_d,e_s)\big) \leq  - a_{\rho_s} \min \{ 1, \underline{a}_{W_s}^2 \} (|x|^2+|e_s|^2) \eqqcolon  - a_s \psi_s^2(|(x,e_s)|)
$,
% \begin{equation*}
%     \begin{aligned}
%         &U_s^\circ(\xi_s; F_s^y(x,0,e_s, 0)) \\
%         \leq & - \big(a_{\rho_s} |x|^2 + a_{\rho_s} W_s^2(\kappa_d,e_s)\big) \\
%         %\leq & - \big(a_{\rho_s} |x|^2 + a_{\rho_s} \underline{a}_{W_s}^2 |e_s|^2\big) \\
%         \leq & - a_{\rho_s} \min \{ 1, \underline{a}_{W_s}^2 \} (|x|^2+|e_s|^2) \\
%         \eqqcolon & - a_s \psi_s^2(|(x,e_s)|),
%     \end{aligned}
% \end{equation*}
which implies \eqref{eqn: Us flow} is satisfied with $a_s \coloneqq a_{\rho_s} \min \{ 1, \underline{a}_{W_s}^2 \}$ and $\psi_s(|(x,e_s)|) \coloneqq |(x,e_s)|$. 
Moreover, we have $\underline{a}_{U_s}|(x,e_s)|^2 \leq U_s(\xi_s) \leq \overline{a}_{U_s}|(x,e_s)|^2,$ where $\underline{a}_{U_s} \coloneqq \min\{\underline{a}_{V_s}, \gamma_s \lambda_s^* \underline{a}_{W_s}^2 \}$ and $\overline{a}_{U_s} \coloneqq \max\{\overline{a}_{V_s}, \gamma_s \tfrac{1}{\lambda_s^*} \overline{a}_{W_s}^2 \}$.

Along the same line as $U_s$, we can proof that \eqref{eqn: Uf flow} is satisfied with $a_f \! \coloneqq \! a_{\rho_f} \min \{ 1, \underline{a}_{W_f}^2 \}$ and $\psi_f(|(y,e_f)|) \coloneqq |(y,e_f)|$. Moreover, we have $\underline{a}_{U_f}|(y,e_f)|^2 \leq U_f(\xi_s, \xi_f) \leq \overline{a}_{U_f}|(y,e_f)|^2,$ where $\underline{a}_{U_f} \coloneqq \min\{\underline{a}_{V_f}, \gamma_f \lambda_f^* \underline{a}_{W_f}^2 \}$ and $\overline{a}_{U_f} \coloneqq \max\{\overline{a}_{V_f}, \gamma_f \tfrac{1}{\lambda_f^*} \overline{a}_{W_f}^2 \}$. Then we satisfy Assumption \ref{Assumption Extra 2} 
% \begin{equation}
% \begin{aligned}
%     \psi_s(|(x,e_s)|) &\leq a_{\psi_s}\sqrt{U_s(\xi_s)}, \\
%     \psi_f(|(y,e_f)|) &\leq a_{\psi_f}\sqrt{U_f(\xi_s,\xi_f)} ,
% \end{aligned}
% \end{equation} \todo{} % Can be replaced by Assumption 5
with $a_{\psi_s} = \underline{a}_{U_s}^{-\frac{1}{2}}$ and $a_{\psi_f} = \underline{a}_{U_f}^{-\frac{1}{2}}$.
Same as the proof of Theorem \ref{Theorem H_1}, we define composite Lyapunov function $U$ as $U(\xi_s, \xi_f) \coloneqq U_s(\xi) + d U_f(\xi_s,\xi_f)$, where $d \in (0,1)$. Then $U$ has sandwich bound 
\begin{equation}
    \underline{a}_{U} |\xi^y|_{\mathcal{E}^y}^2 \leq U(\xi^y) \leq \overline{a}_{U} |\xi^y|_{\mathcal{E}^y}^2, \label{eqn: Exponential U sandwich bound}
\end{equation}
where $\underline{a}_{U} \coloneqq \min \{\underline{a}_{U_s}, d \underline{a}_{U_f} \}$ and $\overline{a}_{U} \coloneqq \max \{\overline{a}_{U_s}, d \overline{a}_{U_f} \}$.

\noindent\emph{\underline{During flow}:} 
We can obtain \eqref{eqn: U derivative during flow eqn1}, as well as
$
U^\circ(\xi^y, 
$
$
F^y(\xi^y, \epsilon)) \leq - 
\left[ \begin{smallmatrix}
            \sqrt{U_s(\xi_s)} \\ \sqrt{U_f(\xi_s, \xi_f)}
        \end{smallmatrix} \right]^T
        \Lambda
        \left[ \begin{smallmatrix}
            \sqrt{U_s(\xi_s)} \\ \sqrt{U_f(\xi_s, \xi_f)}
        \end{smallmatrix} \right]
$,
% \begin{equation*}
%     \begin{aligned}
%         U^\circ(\xi^y, F^y(\xi^y, \epsilon)) \leq - 
%         \left[ \begin{smallmatrix}
%             \sqrt{U_s(\xi_s)} \\ \sqrt{U_f(\xi_s, \xi_f)}
%         \end{smallmatrix} \right]^T
%         \Lambda
%         \left[ \begin{smallmatrix}
%             \sqrt{U_s(\xi_s)} \\ \sqrt{U_f(\xi_s, \xi_f)}
%         \end{smallmatrix} \right],
%     \end{aligned}
% \end{equation*}
where $\Lambda$ is defined in \eqref{eqn: Lambda}, along the same line as the proof of Theorem \ref{Theorem H_1} by setting $\nu_1$ to be zero.
% $\Lambda \coloneqq 
%     \left[ \begin{smallmatrix}
%         a_s a_{\psi_s}^2 & -\tfrac{1}{2}(b_1 + d b_2)a_{\psi_s}a_{\psi_f} \\
%         -\tfrac{1}{2}(b_1 + d b_2)a_{\psi_s}a_{\psi_f} & d (\tfrac{a_f}{\epsilon} - b_3)a_{\psi_f}^2
%     \end{smallmatrix} \right]$. \todo{} % Lambda is same in the proof Theorem 1
In order to satisfy $\Lambda \geq \mu
    \left[ \begin{smallmatrix}
        1 & 0 \\ 0 & d
    \end{smallmatrix} \right]$, where $\mu$ is defined in \eqref{eqn: mu}, we need to satisfy inequality \eqref{eqn: inequality of epsilon Exponential} by having $\epsilon \in (0,\epsilon^*]$, where $\epsilon^*$ is defined by \eqref{eqn: epsilon star} and $d$ in \eqref{eqn: epsilon star} is given later.
% \begin{subequations}
% \begin{align}
%     a_s a_{\psi_s}^2 > \mu \\
%     (a_s a_{\psi_s}^2-\mu) \big( d (\tfrac{a_f}{\epsilon} - b_3)a_{\psi_f}^2 - \mu d \big) &\geq \tfrac{1}{4}(b_1 + db_3)^2a_{\psi_s}^2a_{\psi_f}^2, \label{eqn: inequality of epsilon Exponential}
% \end{align}
% \end{subequations} \todo{}
where the first inequality is satisfied by the definition of $\mu$, and the second inequality can be satisfied by taking $\epsilon$ sufficiently small.
Then we have 
\begin{equation}
    \begin{aligned}
        U^\circ(\xi^y, F^y(\xi^y, \epsilon)) &\leq -\mu (U_s(\xi_s) + d U_f(\xi_s,\xi_f)) \\
        & \leq - \mu U(\xi^y).
        \label{eqn: U dot Exponential}
    \end{aligned}
\end{equation}
\noindent\emph{\underline{During jumps}:} 
Same as the proof of Theorem \ref{Theorem H_1}, we have $U(G_s^y(\xi^y)) \leq a_d U(\xi^y)$ at slow transmissions and $U(G_f^y(\xi^y)) \leq U(\xi^y)$ at fast transmissions.
Suppose $j_k^s, j_{k+1}^s \in \mathcal{J}^s$.
By \eqref{eqn: U dot Exponential}, the fact that $U$ is non-increasing at fast transmissions and comparison principle, we have
\begin{equation}
    U(s,i) \leq U(t_{j_k^s}, j_k^s) \exp \! \big(-\mu(t-t_{j_k^s}) \big) \label{eqn: Exponential U flow - Exponential}, 
\end{equation}
for all $(t_{j_k^s}, j_k^s) \preceq (s,i) \preceq (t_{j_{k+1}^s}, j_{k+1}^s - 1)$ and $(s,i)\in \text{dom}\,\xi^y$.
Along the same line as deriving \eqref{eqn: flow then jump}, we have $U(t_{j_{k+1}^s}, j^s_{k+1}) \leq a_d U(t_{j_k^s}, j_k^s) \exp (-\mu\tau_{\text{miati}}^s )$.
% \begin{equation*}
%     \begin{aligned}
%         U(t_{j_{k+1}^s}, j^s_{k+1}) &\leq a_d U(t_{j_{k+1}^s},j^s_{k+1} - 1)
%         \\
%         &\leq a_d U(t_{j_k^s}, j_k^s) \exp (-\mu\tau_{\text{miati}}^s ).
%     \end{aligned} %
% \end{equation*}
By definition of $a_d$ in \eqref{eqn: a_d}, we have that for any $\tau_{\text{miati}}^{s} \leq T(L_s, \gamma_s, \lambda_s^*)$, $\lambda \in (\exp (-\mu\tau_{\text{miati}}^s ), 1)$, there exist 
\begin{equation}
d^* = \tfrac{-b+\sqrt{b^2-4a \tilde{c}}}{2a},
\label{eqn: d star exponential}
\end{equation}
where $a = \tfrac{\lambda_1}{\gamma_s \lambda_s^*}$, $b= \tfrac{1}{2}( \tfrac{\lambda_1}{\gamma_s \lambda_s^*} + \lambda_2)$ and $\tilde{c}= 1 - \lambda e^{\mu \tau_{\text{miati}}^s}$, such that by taking $d =d^*$, we have $U(t_{j_{k+1}^s}, j^s_{k+1}) \leq \lambda U(t_{j_k^s},j^s_k)$.
% \begin{equation*}
%     U(t_{j_{k+1}^s}, j^s_{k+1}) \leq \lambda U(t_{j_k^s},j^s_k)
% \end{equation*}
Then the inequality \eqref{eqn: inequality of epsilon Exponential} is satisfied by all $\epsilon \in (0, \epsilon^*)$, where $\epsilon^*$ is defined in \eqref{eqn: epsilon star} with $d = d^*$.
By concatenation, we have $U(t_{j_k^s}, j^s_k) \leq  \lambda^{k-1}U(t_{j_1^s}, j_1^s)$.
% \begin{equation*}
%     \begin{aligned}
%         U(t_{j_k^s}, j^s_k) \leq & \lambda^{k-1}U(t_{j_1^s}, j_1^s),
%     \end{aligned}
% \end{equation*}
%
Moreover, since $U$ is non-increasing during flow and upper bounded by $U(G_s^y(\xi^y)) \leq a_d U(\xi^y)$ at slow transmission, we have $ U(t_{j_1^s},j_1^s) \leq  a_d U(0,0) $. Then we have $U(t_{j_k^s}, j_k^s) \leq a_d \lambda^{k-1} U(0,0)$.
% \begin{equation}
%     U(t_{j_k^s}, j_k^s) \leq a_d \lambda^{k-1} U(0,0). \label{eqn: Exponential U slow jump decay}
% \end{equation}
%
Now we have obtained the upper bound of trajectory during the interval between slow transmissions (i.e., \eqref{eqn: Exponential U flow - Exponential}) and the upper bound at each slow transmission. 
Then along the same line as the proof of Claim \ref{Claim U upper bound}, by setting $\Delta$ to infinity and $\nu$ to zero, we can show $U(\xi^y(t,j)) \leq \frac{a_d}{\lambda}U(0,0)  \exp \!\left(-\tfrac{\ln{(\nicefrac{1}{\lambda})}}{\tau_{\text{mati}}^s} t \right)$,
%
% \begin{equation*}
%     U(t,j) \leq \frac{a_d}{\lambda}U(0,0)  \exp \!\left(-\tfrac{\ln{(\nicefrac{1}{\lambda})}}{\tau_{\text{mati}}^s}(t) \right),
% \end{equation*}
for all $(t,j) \in  \text{dom} \ \xi^y$.
 By \eqref{eqn: Exponential U sandwich bound}, we have
 $
 |\xi^y(t,j)|_{\mathcal{E}^y} 
\leq  \big(\tfrac{1}{\underline{a}_U} U(t,j)\big)^{\nicefrac{1}{2}}
\leq  \left(\tfrac{a_d}{\lambda \underline{a}_U}U(0,0)  \exp \!\left(-\tfrac{\ln{(\nicefrac{1}{\lambda})}}{\tau_{\text{mati}}^s}t \right) \right)^{\nicefrac{1}{2}} 
=  \left(\tfrac{a_d \overline{a}_U}{\lambda \underline{a}_U} \right)^{\nicefrac{1}{2}}  |\xi^y(0,0)|_{\mathcal{E}^y}  \exp \!\left(-\tfrac{\ln{(\nicefrac{1}{\lambda})}}{ \tau_{\text{mati}}^s}t \right)^{\nicefrac{1}{2}}
$.
 % \begin{equation*}
 %     \begin{aligned}
 %         &|\xi^y(t,j)|_{\mathcal{E}^y}  
 %         \leq  \big(\tfrac{1}{\underline{a}_U} U(t,j)\big)^{\nicefrac{1}{2}} 
 %         \\
 %         \leq & \left(\tfrac{a_d}{\lambda \underline{a}_U}U(0,0)  \exp \!\left(-\tfrac{\ln{(\nicefrac{1}{\lambda})}}{\tau_{\text{mati}}^s}t \right) \right)^{\nicefrac{1}{2}} 
 %         \\
 %         \leq & \left(\tfrac{a_d}{\lambda \underline{a}_U} \overline{a}_U |\xi^y(0,0)|_{\mathcal{E}^y}^2  \exp \!\left(-\tfrac{\ln{(\nicefrac{1}{\lambda})}}{\tau_{\text{mati}}^s}t \right) \right)^{\nicefrac{1}{2}} \\    
 %         = & \left(\tfrac{a_d \overline{a}_U}{\lambda \underline{a}_U} \right)^{\nicefrac{1}{2}}  |\xi^y(0,0)|_{\mathcal{E}^y}  \exp \!\left(-\tfrac{\ln{(\nicefrac{1}{\lambda})}}{ \tau_{\text{mati}}^s}t \right)^\frac{1}{2}.
 %     \end{aligned}
 % \end{equation*} 
%
Since $\overline{H}$ is globally Lipschitz and $\overline{H}(0,0) = 0$, we have $\overline{H}(x,e_s) \leq L|(x,e_s)|$, where $L$ is the Lipschitz constant. 
 Then by $y = z - \overline{H}(x,e_s)$, there exist $h_1 = 1 + L$ such that
 $|\xi(t,j)|_{\mathcal{E}} \leq h_1 |\xi^y(t,j)|_{\mathcal{E}^y} $ and $|\xi^y(t,j)|_{\mathcal{E}^y} \leq h_1 |\xi(t,j)|_{\mathcal{E}} $. Then the upper bound of $|\xi(t,j)|_{\mathcal{E}} $ is
 $
 |\xi(t,j)|_{\mathcal{E}} \leq h_1^2 \left(\tfrac{a_d \overline{a}_U}{\lambda \underline{a}_U} \right)^{\nicefrac{1}{2}}  |\xi(0,0)|_{\mathcal{E}}  \exp \!\left(-\tfrac{\ln{(\nicefrac{1}{\lambda})}}{ \tau_{\text{mati}}^s}t \right)^{\nicefrac{1}{2}}
 $.
 % \begin{equation*}
 %     |\xi(t,j)|_{\mathcal{E}} \leq h_1^2 \left(\tfrac{a_d \overline{a}_U}{\lambda \underline{a}_U} \right)^{\nicefrac{1}{2}}  |\xi(0,0)|_{\mathcal{E}}  \exp \!\left(-\tfrac{\ln{(\nicefrac{1}{\lambda})}}{ \tau_{\text{mati}}^s}t \right)^\frac{1}{2}.
 % \end{equation*}
% Additionally, since $t \geq \tau_{\text{miati}} (j-1)$, we have
% \begin{equation*}
%     \begin{aligned}
%         &\exp \!\left(-\tfrac{\ln{(\nicefrac{1}{\lambda})}}{2 \tau_{\text{mati}}^s}(t) \right) \\
%         %=&  \exp \!\left(-\tfrac{\ln{(\nicefrac{1}{\lambda})}}{2 \tau_{\text{mati}}^s}(\tfrac{t}{2}+\tfrac{t}{2})) \right) \\
%         \leq & \exp \!\left(-\tfrac{\ln{(\nicefrac{1}{\lambda})}}{2 \tau_{\text{mati}}^s}(\tfrac{t}{2}+\tfrac{\tau_{\text{miati}}^s}{2}(j-1))) \right) \\
%         %=& \exp \!\left( \tfrac{\ln{(\nicefrac{1}{\lambda})}\tau_{\text{miati}}^s}{4 \tau_{\text{mati}}^s}\right)   
%             %\exp \!\left(-\tfrac{\ln{(\nicefrac{1}{\lambda})}}{4 \tau_{\text{mati}}^s}(t+\tau_{\text{miati}}^sj) \right) \\
%         \leq & \exp \!\left( \tfrac{\ln{(\nicefrac{1}{\lambda})}\tau_{\text{miati}}^s}{4 \tau_{\text{mati}}^s}\right)   
%             \exp \!\left(-\tfrac{\ln{(\nicefrac{1}{\lambda})}}{4 \tau_{\text{mati}}^s}  \min\{1,\tau_{\text{miati}}^s \} (t+j) \right).
%     \end{aligned}
% \end{equation*}
By \eqref{eqn: change t to t+j}, we have 
$
|\xi(t,j)|_{\mathcal{E}} \leq c_1 |\xi(0,0)|_{\mathcal{E}}\exp \! \big(- c_2 (t+j)\big)
$,
% \begin{equation*}
%     |\xi(t,j)|_{\mathcal{E}} \leq c_1 |\xi(0,0)|_{\mathcal{E}}\exp \! \big(- c_2 (t+j)\big),
% \end{equation*}
where $c_1 = h_1^2 \left(\tfrac{a_d \overline{a}_U}{\lambda \underline{a}_U} \right)^{\nicefrac{1}{2}} \exp \!\left( \tfrac{\ln{(\nicefrac{1}{\lambda})}\tau_{\text{miati}}^s}{4 \tau_{\text{mati}}^s}\right)$ and $c_2 = \tfrac{\ln{(\nicefrac{1}{\lambda})}}{4 \tau_{\text{mati}}^s}  \min\{1,\tau_{\text{miati}}^s \}$.
%Now we have shown $\mathcal{H}_1$ is uniformly globally pre-exponentially stable w.r.t $\mathcal{E}$, and we can prove $\mathcal{H}_1$ is UGES along the same line as the proof of Theorem \ref{Theorem H_1} (i.e. completeness of solution). $\hfill\blacksquare$

%\input{Chapters/Appendix_Stable fast subsystem}     % Sections and subsections are supported  

\section{Proof of Proposition \ref{Proposition LTI}}

Claim \ref{Claim for LTI section} has shown that \eqref{eqn: NCS assumption Ws sandwich bound}-\eqref{eqn: NCS Ws dot} in Assumption \ref{Assumption reduced model} hold and \eqref{eqn: Ws exponential sandwich bound} in Assumption \ref{Assumption Exponential} holds. 
Next, we show \eqref{eqn: NCS Vs flow} in Assumption \ref{Assumption reduced model}, as well as \eqref{eqn: Vs exponential sandwich bound} and \eqref{eqn: a_{rho_s}} in Assumption \ref{Assumption Exponential} hold.
Let $P_s = \left[\begin{smallmatrix}
    p_{11}^s  & \bigstar \\ {p_{12}^{s\top}} & p_{22}^s
\end{smallmatrix} \right] > 0$, where $p_{11}^s$ is a $n_{x_p} $ by $ n_{x_p}$ symmetric matrix, $p_{12}^s$ is a $n_{x_p} $ by $ n_{x_c}$ matrix and $p_{22}^s$ is a $n_{x_c} $ by $ n_{x_c}$ symmetric matrix. Let $V_s = x^\top P_s x$, then \eqref{eqn: Vs exponential sandwich bound} is satisfied with $\underline{a}_{V_s} = \lambda_{\text{min}}(P_s)$ and $\overline{a}_{V_s} = \lambda_{\text{max}}(P_s)$. Moreover, we have
\begin{equation}
    \begin{aligned}
        &\left< \tfrac{\partial {V_s}(x)}{\partial x},f_x(x,\overline{H}(x,e_s),e_s, 0) \right>   \\
        =& x^\top (P_s A_{11}^s + A_{11}^{s\top} P_s) x + x^\top P_s A_{12}^s e_s + e_s^\top A_{12}^{s\top} P_s x .
    \end{aligned}
        \label{eqn: Linear case Vs dot}
\end{equation}
Inequalities \eqref{eqn: NCS Vs flow} and \eqref{eqn: a_{rho_s}} are satisfied if
\eqref{eqn: Linear case Vs dot inequality} holds.
\begin{equation}
    \begin{aligned}
        &\left< \tfrac{\partial {V_s}(x)}{\partial x},f_x(x,\overline{H}(x,e_s),e_s, 0) \right>   \\
        \leq & -a_{\rho_s} x^\top x - a_{\rho_s} e_s^\top e_s - x^\top A_{H_s}^\top A_{H_s} x  + \gamma_s^2 \underline{a}_{W_s}^2 e_s^\top e_s.
    \end{aligned}
    \label{eqn: Linear case Vs dot inequality}
\end{equation}
By substituting \eqref{eqn: Linear case Vs dot} into \eqref{eqn: Linear case Vs dot inequality}, we show that \eqref{eqn: NCS Vs flow} in Assumption \ref{Assumption reduced model} and \eqref{eqn: a_{rho_s}} in Assumption \ref{Assumption Exponential} are satisfied if \eqref{eqn: LMIs} with $\ell = s$ holds.
% \begin{equation}
%     \left[\begin{smallmatrix}
%         A_{11}^s P_s + P_s A_{11}^{s\top} + a_{\rho_s} I + A_{H_s}^\top A_{H_s} &  \bigstar  \\
%         A_{12}^{s\top} P_s & a_{\rho_s} I - \gamma_s^2 \underline{a}_{W_s}^2 I
%     \end{smallmatrix}\right]
%     \leq 0.
%     \label{eqn: LMI slow}
% \end{equation}

Similarly, Claim \ref{Claim for LTI section} show \eqref{eqn: NCS assumption Wf sandwich bound}-\eqref{eqn: NCS Wf dot} and \eqref{eqn: Ws exponential sandwich bound} are satisfied. 
Let $P_f = \left[\begin{smallmatrix}
    p_{11}^f  & \bigstar \\ {p_{12}^{f\top}} & p_{22}^f
\end{smallmatrix} \right] > 0$, where $p_{11}^f$ is a $n_{z_p} $ by $ n_{z_p}$ symmetric matrix, $p_{12}^f$ is a $n_{z_p} $ by $ n_{z_c}$ matrix and $p_{22}^f$ is a $n_{z_c} $ by $ n_{z_c}$ symmetric matrix. Let $V_f = y^\top P_f y$, then \eqref{eqn: Vs exponential sandwich bound} is satisfied with $\underline{a}_{V_f} = \lambda_{\text{min}}(P_f)$ and $\overline{a}_{V_f} = \lambda_{\text{max}}(P_f)$.
Moreover, we can show \eqref{eqn: NCS Vf flow} in Assumption \ref{Assumption boundary layer system} and \eqref{eqn: a_{rho_f}} in Assumption \ref{Assumption Exponential} hold if LMI \eqref{eqn: LMIs} with $\ell = f$ is satisfied.
% \begin{equation}
%     \left[\begin{smallmatrix}
%         A_{11}^f P_f + P_f A_{11}^{f\top} + a_{\rho_f} I + A_{H_f}^\top A_{H_f} & \bigstar \\
%         A_{12}^{f\top} P_f & a_{\rho_f} I - \gamma_f^2 \underline{a}_{W_f}^2 I
%     \end{smallmatrix}\right] \leq 0.
%     \label{eqn: LMI fast}
% \end{equation}
At this point, we show Assumptions \ref{Assumption reduced model}, \ref{Assumption boundary layer system} and \ref{Assumption Exponential} hold if the LMI \eqref{eqn: LMIs} with $\ell \in \{s,f\}$ is satisfied.

We then Verify Assumptions \ref{Assumption Vf at slow transmission} and \ref{Assumption interconnection Exponential}. By definition of $h_y(\kappa_s,x,e_s,y)$ in \eqref{eqn: Jump of y at slow transmission} and $\overline{H}$ in \eqref{eqn: H bar linear}, we have $h_y(\kappa_s,x,e_s,y) =  y -A_{33}^{-1} A_{32} (e_s - h_s(\kappa_s, e_s)) $.
% \begin{equation}
%     \begin{aligned}
%         h_y(\cyan{\kappa_s,}x,e_s,y)
%         % =& y + \overline{H}(x,e_s) - \overline{H}(x,h_s(\kappa_s, e_s))
%         % \\
%         % =& y +  (- A_{33}^{-1} A_{31} x - A_{33}^{-1} A_{32} e_s) - 
%         %     \\
%         %     & ( - A_{33}^{-1} A_{31} x - A_{33}^{-1} A_{32} h_s(\kappa_s, e_s))
%         % \\
%         =  y -A_{33}^{-1} A_{32} (e_s - h_s(\kappa_s, e_s)).
%     \end{aligned}
%     \label{eqn: h_y linear}
% \end{equation}
Since we assumed when a slow node gets access to the network, some elements of $e_s$ reset to zero, we have $|e_s - h_s(\kappa_s, e_s)| \leq |e_s|$.
% \begin{equation*}
%     |e_s - h_s(\kappa_s, e_s)| \leq |e_s|.
% \end{equation*}
Then by definition of $h_y$, we have 
% $
% V_f(x, h_y(x,e_s,y)) - V_f(x,y)
% %
% =  (y -A_{33}^{-1} A_{32} (e_s - h_s(\kappa_s, e_s)))^\top P_f (y -A_{33}^{-1} A_{32} (e_s - h_s(\kappa_s, e_s))) - y^\top P_f y
% %
% \leq 2 | P_f A_{33}^{-1} A_{32}| |y| |e_s| + |A_{32}^\top A_{33}^{-1\top} P_f A_{33}^{-1} A_{32}| |e_s|^2
% %
% \leq  \lambda_1 W_s^2(\kappa_s, e_s) + \lambda_2 \sqrt{W_s^2(\kappa_s, e_s) V_f(x,y)}
% $,
\begin{small}
\begin{equation}
\setlength\abovedisplayskip{-5pt}%shrink space
    \begin{aligned}
        &V_f(x, h_y(\kappa_s,x,e_s,y)) - V_f(x,y) \\
        %= & h_y^\top(x,e_s, y) P_f h_y - y^\top P_f y \\
        = & (y -A_{33}^{-1} A_{32} (e_s - h_s(\kappa_s, e_s)))^\top P_f \\
            &(y -A_{33}^{-1} A_{32} (e_s - h_s(\kappa_s, e_s))) - y^\top P_f y  \\
        \leq & 2 | P_f A_{33}^{-1} A_{32}| |y| |e_s| + |A_{32}^\top A_{33}^{-1\top} P_f A_{33}^{-1} A_{32}| |e_s|^2 \\
        \leq & \lambda_1 W_s^2(\kappa_s, e_s) + \lambda_2 \sqrt{W_s^2(\kappa_s, e_s) V_f(x,y)} ,
    \end{aligned}
    \label{eqn: lambda_1 and lambda_2}
\end{equation}
\end{small}
where $\lambda_1 = \tfrac{1}{\underline{a}_{W_s}^2}  |A_{32}^\top A_{33}^{-1\top} P_f A_{33}^{-1} A_{32}|$ and $\lambda_2 = \tfrac{2}{\underline{a}_{W_s}  \sqrt{\underline{a}_{V_f}}  } | P_f A_{33}^{-1} A_{32}| $. We have shown that we satisfy Assumption \ref{Assumption Vf at slow transmission}. Next, we show that Assumption \ref{Assumption interconnection Exponential} always hold. We first verify inequality \eqref{eqn: SPNCS interconnection Exponential 1}. We have 
$
\tfrac{\partial U_s}{\partial \xi_s} \!= \!
        \left[\! \begin{smallmatrix}
        2 x^\top \! P_s &
        2\gamma_s \phi_s(\tau_s) W_s(\kappa_s, e_s) \tfrac{\partial W_s}{\partial e_s} &
        -\gamma_s^2(\phi_s^2(\tau_s) + 1 ) W_s(\kappa_s, e_s)^2 &
        0
    \end{smallmatrix}\!\right]
$.
% \begin{equation*}
%     \begin{aligned}
%     \tfrac{\partial U_s}{\partial \xi_s} 
%     &= \left[ \begin{smallmatrix} \tfrac{\partial U_s}{\partial x} &\tfrac{\partial U_s}{\partial e_s} &\tfrac{\partial U_s}{\partial \tau_s} &\tfrac{\partial U_s}{\partial \kappa_s}\end{smallmatrix} \right]
%     \\
%     &=\left[ \begin{smallmatrix}
%         (2 x^\top P_s)^\top \\
%         (2\gamma_s \phi_s(\tau_s) W_s(\kappa_s, e_s) \tfrac{\partial W_s}{\partial e_s})^\top \\
%         \gamma_s(-\gamma_s(\phi_s^2(\tau_s) + 1 )) W_s(\kappa_s, e_s)^2 \\
%         0
%     \end{smallmatrix} \right]^\top.
%     \end{aligned}
% \end{equation*}
Additionally, we have 
$   F_s^y(x, y, e_s, e_f) \!= \!\!
    \left[ \begin{smallmatrix}
        A_{11}^s & A_{12}^s & A_{13} & A_{14} \\
        A_{21}^s & A_{22}^s & A_{23} & A_{24} \\
        0 & 0& 0 & 0 \\
        0 & 0& 0 & 0
    \end{smallmatrix} \right] 
$ 
$
    \left[ \begin{smallmatrix}
        x \\ e_s \\ y \\ e_f
    \end{smallmatrix} \right]
    +
    \left[ \begin{smallmatrix}
        0 \\ 0 \\ 1 \\ 0
    \end{smallmatrix} \right]
$,
%
% \begin{equation*}
%     F_s^y(x, y, e_s, e_f) = 
%     \left[ \begin{smallmatrix}
%         A_{11}^s & A_{12}^s & A_{13} & A_{14} \\
%         A_{21}^s & A_{22}^s & A_{23} & A_{24} \\
%         0 & 0& 0 & 0 \\
%         0 & 0& 0 & 0
%     \end{smallmatrix} \right]
%     \left[ \begin{smallmatrix}
%         x \\ e_s \\ y \\ e_f
%     \end{smallmatrix} \right]
%     +
%     \left[ \begin{smallmatrix}
%         0 \\ 0 \\ 1 \\ 0
%     \end{smallmatrix} \right],
% \end{equation*}
which implies
$$
F_s^y(x,y,e_s,e_f) - F_s^y(x,0,e_s,0) = 
    \left[ \begin{smallmatrix}
        A_{13}y + A_{14}e_f \\ A_{23}y + A_{24}e_f \\ 0 \\ 0
    \end{smallmatrix} \right].
$$
%
%
% \begin{equation*}
%     F_s^y(x,y,e_s,e_f) - F_s^y(x,0,e_s,0) = 
%     \left[ \begin{smallmatrix}
%         A_{13}y + A_{14}e_f \\ A_{23}y + A_{24}e_f \\ 0 \\ 0
%     \end{smallmatrix} \right].
% \end{equation*}
By \cite[Remark 11]{dragan_stability}, there exist $L_1 \geq 0$ such that $\left|\tfrac{\partial W_s(\kappa_s,e_s)}{\partial e_s} \right| \leq L_1$, then \eqref{eqn: SPNCS interconnection Exponential 1} is satisfied by
\begin{small}
\begin{equation}
\setlength\abovedisplayskip{-4pt}%shrink space
    \begin{aligned}
        &\Big < \tfrac{\partial {U_s}}{\partial \xi_s}, F_s^y(x,y,e_s,e_f) - F_s^y(x,0,e_s,0)  \Big>
        \\ 
        = & 2 x^\top P_s (A_{13}y + A_{14}e_f) + 
            \\
            & 2 \gamma_s \phi_s(\tau_s)W_s(\kappa_s,e_s) \tfrac{\partial W_s}{\partial e_s}(A_{23} y + A_{24}e_f)
        \\
        \leq & \left[ \begin{smallmatrix}
        |x| \\ |e_s|
    \end{smallmatrix} \right]^\top
    \Lambda_{b_1}
    \left[ \begin{smallmatrix}
        |y| \\ |e_f|
    \end{smallmatrix} \right]
    \leq  b_1 \psi_s(|(x,e_s)|) \psi_f(|(y,e_f)|),
    \end{aligned}
    \label{eqn: Lambda_b1}
\end{equation}
\end{small}
where 
$\Lambda_{b_1} = \left[\begin{smallmatrix}
    |P_s A_{13}| & |P_s A_{14}| \\ \tfrac{\gamma_s}{\lambda_s^*}\overline{a}_{W_s} L_1 |A_{22}| & \tfrac{\gamma_s}{\lambda_s^*}\overline{a}_{W_s} L_1 |A_{24}|
\end{smallmatrix}\right]$, $b_1 = \sqrt{\lambda_{\text{max}}(\Lambda_{b_1}^\top \Lambda_{b_1})}$ and $\psi_s(s) = \psi_f(s) = s$ for all $s \in \mathbb{R}_{\geq 0}$.
Finally, we validate the inequality \eqref{eqn: SPNCS interconnection Exponential 2} in Assumption \ref{Assumption interconnection Exponential}. By definition of $U_f$ in \eqref{eqn: definition of U_f}, we have 
$\tfrac{\partial U_f}{\partial \xi_s} = 0$,
$\tfrac{\partial U_f}{\partial y} = 2 y^\top P_f$, 
$\tfrac{\partial \overline{H}}{\partial \xi_s} = \left[ \begin{smallmatrix}
            -A_{33}^{-1} A_{31} & -A_{33}^{-1} A_{32} & 0 &0
\end{smallmatrix} \right] $,
$\tfrac{\partial U_f}{\partial e_f} = 2 \gamma_f \phi_f(\tau_f)W_f(\kappa_f, e_f) \tfrac{\partial W_f}{\partial e_f}$,
and
$\tfrac{\partial \tilde{k}}{\partial \xi_s} = \left[ \begin{smallmatrix}
            \left[\begin{smallmatrix}
                A_x^{p_f} & 0 \\ 0 & A_x^{c_f}
            \end{smallmatrix}\right] & 0 & 0 & 0\end{smallmatrix} \right]$.
%
% \begin{equation*}
%     \begin{aligned}
%         \tfrac{\partial U_f}{\partial \xi_s} &= 0, \qquad \tfrac{\partial U_f}{\partial y} = 2 y^\top P_f \\
%         \tfrac{\partial \overline{H}}{\partial \xi_s} &= \left[ \begin{smallmatrix}
%             -A_{33}^{-1} A_{31} & -A_{33}^{-1} A_{32} & 0 &0
%         \end{smallmatrix} \right] ,
%         \\
%         \tfrac{\partial U_f}{\partial e_f} &= 2 \gamma_f \phi_f(\tau_f)W_f(\kappa_f, e_f) \tfrac{\partial W_f}{\partial e_f},
%         \\
%         \tfrac{\partial \tilde{k}}{\partial \xi_s} &= \left[ \begin{smallmatrix}
%             \left[\begin{smallmatrix}
%                 A_x^{p_f} & 0 \\ 0 & A_x^{c_f}
%             \end{smallmatrix}\right] & 0 & 0 & 0
%         \end{smallmatrix} \right].
%     \end{aligned}
% \end{equation*}
Then along the same line as \eqref{eqn: Lambda_b1}, we can show that there exist a matrix $\Lambda_{b_2}$, a symmetric matrix $\Lambda_{b_3}$ and $b_2$, $b_3 \geq 0$ such that
% $
% \Big< \tfrac{\partial {U_f}}{\partial \xi_s} - \tfrac{\partial {U_f}}{\partial y} \tfrac{\partial \overline{H}}{\partial \xi_s} - \tfrac{\partial {U_f}}{\partial e_f} \tfrac{\partial \tilde k}{\partial \xi_s} ,  F_s^y(x,y,e_s,e_f) \Big>
% %
% \leq  \left[ \begin{smallmatrix}
%             |x| \\ |e_s|
%         \end{smallmatrix} \right]^\top
%         \Lambda_{b_2}
%         \left[ \begin{smallmatrix}
%             |y| \\ |e_f|
%         \end{smallmatrix} \right]
%         + 
%         \left[ \begin{smallmatrix}
%         |y| \\ |e_f|
%         \end{smallmatrix} \right]^\top
%         \Lambda_{b_3}
%         \left[ \begin{smallmatrix}
%             |y| \\ |e_f|
%         \end{smallmatrix} \right]
% %
% \leq    b_2 \psi_s\left(\left| (x, e_s) \right|\right) $ $ \psi_f\left(\left| (y, e_f) \right|\right) + b_3 \psi_f^2\left(\left| (y, e_f) \right|\right)
% $,
\begin{equation}
    \begin{aligned}
        &\Big< \tfrac{\partial {U_f}}{\partial \xi_s} - \tfrac{\partial {U_f}}{\partial y} \tfrac{\partial \overline{H}}{\partial \xi_s} - \tfrac{\partial {U_f}}{\partial e_f} \tfrac{\partial \tilde k}{\partial \xi_s} ,  F_s^y(x,y,e_s,e_f) \Big>
        \\
        \leq & \left[ \begin{smallmatrix}
            |x| \\ |e_s|
        \end{smallmatrix} \right]^\top
        \Lambda_{b_2}
        \left[ \begin{smallmatrix}
            |y| \\ |e_f|
        \end{smallmatrix} \right]
        + 
        \left[ \begin{smallmatrix}
        |y| \\ |e_f|
        \end{smallmatrix} \right]^\top
        \Lambda_{b_3}
        \left[ \begin{smallmatrix}
            |y| \\ |e_f|
        \end{smallmatrix} \right]
        \\
        \leq &  b_2 \psi_s\left(\left| (x, e_s) \right|\right) \psi_f\left(\left| (y, e_f) \right|\right) + b_3 \psi_f^2\left(\left| (y, e_f) \right|\right),
    \end{aligned}
    \label{eqn: Lambda_b2 and Lambda_b3}
\end{equation}
where $b_2 \!= \! \sqrt{\lambda_{\text{max}}(\Lambda_{b_2}^\top \Lambda_{b_2}) }$, $b_3 = \lambda_{\text{max}}(\Lambda_{b_3})$, which implies the inequality \eqref{eqn: SPNCS interconnection Exponential 2} is satisfied.

% \bibliographystyle{plain}        % Include this if you use bibtex 
% \bibliography{autosam}           % and a bib file to produce the 

% \begin{wrapfigure}{l}{25mm} 
%     \includegraphics[width=1in,height=1.25in,clip,keepaspectratio]{Biography/romain.jpg}
%   \end{wrapfigure}\par
%   \textbf{Romain Postoyan} received the ``Ing\'enieur'' degree in Electrical and Control Engineering from ENSEEIHT (France) in 2005. He obtained the M.Sc. by Research in Control Theory \& Application from Coventry University (United Kingdom) in 2006 and the Ph.D. in Control Engineering from Universit\'e Paris-Sud (France) in 2009. In 2010, he was a research assistant at the University of Melbourne (Australia). Since 2011, he is a CNRS researcher at the ``Centre de Recherche en Automatique de Nancy'' (France). He received the `Habilitation à Diriger des Recherches (HDR)'' in 2019 from Université de Lorraine (Nancy, France). He serves/served as an associate editor for the journals: IEEE Transactions on Automatic Control, Automatica, IEEE Control Systems Letters and IMA Journal of Mathematical Control and Information.\par

                                        % in the appendices.
\end{document}